\documentclass[10pt]{article}
\usepackage[T1]{fontenc}
\usepackage[latin9]{inputenc}
\usepackage{color}
\usepackage[english]{babel}

\usepackage{textcomp}
\usepackage{amsthm}
\usepackage{amsmath}
\usepackage{amssymb}
\usepackage[all]{xy}
\usepackage[unicode=true, pdfusetitle,
 bookmarks=true,bookmarksnumbered=false,bookmarksopen=false,
 breaklinks=false,pdfborder={0 0 1},backref=false,colorlinks=true]
 {hyperref}

\theoremstyle{plain}

\usepackage{mathrsfs}
\usepackage{hyperref}
\usepackage{xy}
\hypersetup{colorlinks=true,citecolor=blue}
\begin{document}

\title{Universal coefficients for overconvergent cohomology \\
and the geometry of eigenvarieties}

\author{David Hansen%
\thanks{Department of Mathematics, Boston College, Chestnut Hill MA; dvd.hnsn@gmail.com%
}\\
with an appendix by James Newton}
\maketitle
\begin{abstract}
We prove a universal coefficients theorem for the overconvergent cohomology
modules introduced by Ash and Stevens, and give several applications.
In particular, we sketch a very simple construction of eigenvarieties
using overconvergent cohomology and prove many instances of a conjecture
of Urban on the dimensions of these spaces. For example, when the
underlying reductive group is an inner form of $\mathrm{GL}_{2}$
over a quadratic imaginary extension of $\mathbf{Q}$, the cuspidal
component of the eigenvariety is a rigid analytic curve.

\tableofcontents{}
\end{abstract}

\section{Introduction}

\subsection{Background and results}

Since the pioneering works of Serre \cite{SeAntwerp}, Katz \cite{Katzmodschemes},
and especially Hida \cite{Hidordinv,Hidtotreal,HidordGL2} and Coleman
 \cite{Colemanclassicaloc,Colemanbanachfamilies}, \emph{p-}adic families
of modular forms have become a major topic in modern number theory.
Aside from their intrinsic beauty, these families have found applications
towards Iwasawa theory, the Bloch-Kato conjecture, modularity lifting
theorems, and the local and global Langlands correspondences \cite{BCnontemp,BCast,BCsign,BLGGT,Sknote,Wilesiwasawatotreal}.
One of the guiding examples in the field is Coleman and Mazur's\emph{
eigencurve} \cite{CMeigencurve,BuEigen}, a universal object parametrizing
all overconvergent \emph{p-}adic modular forms of fixed tame level
and finite slope. Concurrently with their work, Stevens introduced
his beautifully simple idea of \emph{overconvergent cohomology} \cite{Strigidsymbs},
a group-cohomological avatar of overconvergent \emph{p}-adic modular
forms\emph{.} Ash and Stevens developed these ideas much further in \cite{AS}:
as conceived there, overconvergent cohomology works for any connected
reductive $\mathbf{Q}$-group $G$ split at $p$, and leads to natural
candidates for quite general eigenvarieties. When the group $G^{\mathrm{der}}(\mathbf{R})$
possesses discrete series representations, Urban \cite{UrEigen} used
overconvergent cohomology to construct eigenvarieties interpolating
classical forms with nonzero Euler-Poincaré multiplicities, showing
that his construction yields spaces which are equidimensional of the
same dimension as weight space. In this article we develop new tools
to analyze the situation for general groups.

To describe our results, we first introduce some notation. Fix a reductive
$\mathbf{Q}$-group scheme $G$ with $G^{\mathrm{der}}(\mathbf{R})$
noncompact, and fix a prime $p$ such that $G$ is split over $\mathbf{Q}_{p}$.
Choose a Borel subgroup $B=TN$ and Iwahori subgroup $I$ of $G(\mathbf{Q}_{p})$.
Let $\mathcal{W}=\mathrm{Hom}_{\mathrm{cts}}(T(\mathbf{Z}_{p}),\mathbf{G}_{m})^{\mathrm{an}}$
be the rigid analytic space of \emph{p-}adically continuous characters
of $T(\mathbf{Z}_{p})$; as a rigid space, $\mathcal{W}$ is a disjoint
union of open balls, each of dimension equal to the rank of $G$.%
\footnote{By {}``rank'' we shall always mean {}``absolute rank'', i.e. the
dimension of any maximal torus, split or otherwise.%
} Multiplication of characters gives $\mathcal{W}$ the structure of
a rigid analytic group. Given a continuous character $\lambda:T(\mathbf{Z}_{p})\to\overline{\mathbf{Q}_{p}}^{\times}$,
we also denote the corresponding point of $\mathcal{W}(\overline{\mathbf{Q}_{p}})$
by $\lambda$. A weight $\lambda$ is \emph{arithmetic }if $\lambda=\lambda^{\mathrm{alg}}\varepsilon$
where $\lambda^{\mathrm{alg}}$ is a $B$-dominant algebraic character
of $T$ and $\varepsilon$ is a character of finite order. In §2.1
we define, given an arbitrary affinoid open subset $\Omega\subset\mathcal{W}$
and an integer $s\gg0$, an orthonormalizable Banach $A(\Omega)$-module
$\mathbf{A}_{\Omega}^{s}$ and its dual module $\mathbf{D}_{\Omega}^{s}\simeq\mathrm{Hom}_{A(\Omega)}^{\mathrm{cts}}(\mathbf{A}_{\Omega}^{s},A(\Omega))$.%
\footnote{If $X$ is an affinoid variety, $A(X)$ denotes the underlying affinoid
algebra.%
} The module $\mathbf{A}_{\Omega}^{s}$ admits a canonical continuous
left $A(\Omega)$-linear $I$-action, and $\mathbf{D}_{\Omega}^{s}$
inherits a corresponding right action. Furthermore, there are canonical
$A(\Omega)[I]$-equivariant transition maps $\mathbf{A}_{\Omega}^{s}\to\mathbf{A}_{\Omega}^{s+1}$
which are injective and compact; we set\[
\mathcal{A}_{\Omega}=\lim_{s\to\infty}\mathbf{A}_{\Omega}^{s}\]
and\[
\mathcal{D}_{\Omega}=\lim_{\infty\leftarrow s}\mathbf{D}_{\Omega}^{s}.\]
For any rigid Zariski closed subspace $\Sigma\subset\Omega$, let
$A(\Omega)\to A(\Sigma)$ denote the corresponding surjection on coordinate
rings, and set $\mathcal{D}_{\Sigma}=\mathcal{D}_{\Omega}\otimes_{A(\Omega)}A(\Sigma)$.
If $\lambda$ is a classical dominant weight for $B$, with $V_{\lambda}$
the corresponding irreducible right $G(\mathbf{Q}_{p})$-representation
of highest weight $\lambda$, there is a canonical continuous, surjective
$I$-equivariant {}``integration'' map $i_{\lambda}:\mathcal{D}_{\Omega}\to V_{\lambda}$
which factors through the map $\mathcal{D}_{\Omega}\to\mathcal{D}_{\lambda}$.
There's no hope of finding a coherent sheaf of modules over $\Omega$
with the $V_{\lambda}$'s as fibers, since the dimensions of the latter
spaces vary unboundedly; in this light, $\mathcal{D}_{\Omega}$ \emph{p-}adically
interpolates the $V_{\lambda}$'s in about the most straightforward
way possible.

Let $K_{\infty}$ be a maximal compact subgroup of $G(\mathbf{R})$,
and let $Z_{\infty}$ be the real points of the center of $G$. For
any open compact subgroup $K^{p}\subset G(\mathbf{A}_{f}^{p})$ in
the prime-to-$p$ finite adeles of $G$, the modules $\mathcal{A}_{\Omega}$
and $\mathcal{D}_{\Omega}$ give rise to local systems on the Shimura
manifold $Y(K^{p}I)=G(\mathbf{Q})\backslash G(\mathbf{A})/K^{p}IK_{\infty}Z_{\infty}$.
If $M$ is any $I$-module, we write $H_{\ast}(K^{p},M)$ and $H^{\ast}(K^{p},M)$
for the homology and cohomology of the local system induced by $M$
on $Y(K^{p}I)$. Ash and Stevens showed that the maps $\mathcal{D}_{\lambda}\to V_{\lambda}$
give rise to degree-preserving Hecke-equivariant morphisms\[
H^{\ast}(K^{p},\mathcal{D}_{\lambda})\to H^{\ast}(K^{p},V_{\lambda})\]
with cokernel contained in the space of forms of {}``critical slope''.
The target, by Matsushima's formula and its generalizations, is isomorphic
as a Hecke module to a finite-dimensional space of classical automorphic
forms; the source, on the other hand, is much larger. In the case
of $G=\mathrm{GL_{2}}/\mathbf{Q}$, fundamental unpublished work
of Stevens asserts that the modules $H^{1}(K^{p},\mathcal{D}_{\lambda})$
contain exactly the same finite-slope Hecke data as spaces of \emph{p}-adic
overconvergent modular forms, and that $H^{1}(K^{p},\mathcal{D}_{\Omega})$
interpolates the $H^{1}(K^{p},\mathcal{D}_{\lambda})$'s in a natural
way. In general, the overconvergent cohomology modules $H^{\ast}(K^{p},\mathcal{D}_{\lambda})$
and their interpolations $H^{\ast}(K^{p},\mathcal{D}_{\Omega})$ seem
to be an excellent surrogate for spaces of overconvergent \emph{p-}adic
modular forms.

Our first main result is an analogue of the universal coefficients
theorem for overconvergent cohomology. As the duality functor $\mathrm{Hom}_{A(\Omega)}(-,A(\Omega))$
is far from exact, this result naturally takes the form of a spectral
sequence. Before stating it, though, we need to explain one key difficulty.
The structure of the individual $A(\Omega)$-modules $H_{i}(K^{p},\mathcal{A}_{\Omega})$
and $H^{i}(K^{p},\mathcal{D}_{\Omega})$ is very mysterious; in general
they are not flat or finitely generated over $A(\Omega)$, and it
seems unclear whether they're even Hausdorff in their natural topology.\emph{
}More precisely, $H_{\ast}(K^{p},\mathcal{A}_{\Omega})$ is the homology
of a non-canonical chain complex $C_{\bullet}(K^{p},\mathcal{A}_{\Omega})$
of topological $A(\Omega)$-modules, and there's no particular reason
for the boundaries of this complex to form a closed subspace of the
cycles. To avoid this difficulty, as well the fact that the derived
functors of $\mathrm{Hom}_{A(\Omega)}^{\mathrm{cts}}(-,A(\Omega))$
only make sense in the setting of \emph{relative} homological algebra,
we appeal to the theory of slope decompositions. More precisely, we
define in §2.1 a monoid $\Delta\subset G(\mathbf{Q}_{p})$ containing
$I$ which acts on the modules described above and extends their $I$-actions,
and a submonoid $\Delta^{+}\subset\Delta$ which acts completely continuously.
The Hecke algebra $\mathcal{A}_{p}=\mathcal{H}_{\mathbf{Q}_{p}}(I\backslash\Delta/I)$
is commutative, and we define \[
\mathbf{T}(K^{p})=\mathcal{A}_{p}\otimes_{\mathbf{Q}_{p}}\mathcal{H}_{\mathbf{Q}_{p}}(K^{p}\backslash G(\mathbf{A}_{f}^{p})/K^{p})^{\mathrm{sph}};\]
this is a commutative, non-Noetherian $\mathbf{Q}_{p}$-algebra, which
acts on homology and cohomology in the usual way. For any element
$t\in\Delta^{+}$, we set $U_{t}:=[ItI]\in\mathbf{T}(K^{p})$ and
refer to $U_{t}$ as a \emph{controlling operator, }in analogy with
the classical $U_{p}$-operator of Atkin and Lehner and its fundamental
role in the theory of overconvergent modular forms. In §2.2 we explain
how to lift the operator $U_{t}$ from homology to a compact $A(\Omega)$-linear
operator $\tilde{U}_{t}$ on the chain complex $C_{\bullet}(K^{p},\mathcal{A}_{\Omega})$.
By combining some remarkable results of Buzzard and Ash-Stevens, we
find that for any given rational number $h=a/b\in\mathbf{Q}_{\geq0}$
and any point $x\in\mathcal{W}$, there is an affinoid $\Omega\subset\mathcal{W}$
containing $x$ together with a \emph{slope-$\leq h$ decomposition}
\[
C_{\bullet}(K^{p},\mathcal{A}_{\Omega})\simeq C_{\bullet}(K^{p},\mathcal{A}_{\Omega})_{\leq h}\oplus C_{\bullet}(K^{p},\mathcal{A}_{\Omega})_{>h}\]
into closed $A(\Omega)[\tilde{U_{t}}]$-stable subcomplexes, where
$C_{\bullet}(K^{p},\mathcal{A}_{\Omega})_{>h}$ is roughly the maximal
$A(\Omega)$-subcomplex on which $p^{-a}\tilde{U}_{t}^{b}$ acts topologically
nilpotently and the complex $C_{\bullet}(K^{p},\mathcal{A}_{\Omega})_{\leq h}$
is a complex of finite flat $A(\Omega)$-modules (we give a precise,
much more general definition of slope-$\leq h$ decompositions in
§2.3). Passing to homology yields a decomposition of Hecke modules\[
H_{\ast}(K^{p},\mathcal{A}_{\Omega})\simeq H_{\ast}(K^{p},\mathcal{A}_{\Omega})_{\leq h}\oplus H_{\ast}(K^{p},\mathcal{A}_{\Omega})_{>h},\]
with $H_{\ast}(K^{p},\mathcal{A}_{\Omega})_{\leq h}$ finitely presented
as an $A(\Omega)$-module. We define a \emph{slope datum }as a triple
$(U_{t},\Omega,h)$ where $U_{t}$ is a controlling operator, $\Omega\subset\mathcal{W}$
is a connected affinoid open subset, and $h\in\mathbf{Q}_{\geq0}$
is such that $C_{\bullet}(K^{p},\mathcal{A}_{\Omega})$ admits a slope-$\leq h$
decomposition for the $\tilde{U}_{t}$-action. For the actual construction
of eigenvarieties the choice of a controlling operator $U_{t}$ is
immaterial (and for many groups there is a canonical choice), so we
shall sometimes implicitly fix a controlling operator and then refer
to the pair $(\Omega,h)$ as a slope datum\emph{.}

\textbf{Theorem 1.1. }\emph{Fix a slope datum $(U_{t},\Omega,h)$,
and let $\Sigma\subseteq\Omega$ be an arbitrary rigid Zariski closed
subspace. Then $H^{\ast}(K^{p},\mathcal{D}_{\Sigma})$ admits a slope-$\leq h$
decomposition, and there is a convergent first quadrant spectral sequence\[
E_{2}^{i,j}=\mathrm{Ext}_{A(\Omega)}^{i}(H_{j}(K^{p},\mathcal{A}_{\Omega})_{\leq h},A(\Sigma))\Rightarrow H^{i+j}(K^{p},\mathcal{D}_{\Sigma})_{\leq h}.\]
Furthermore, there is a convergent second quadrant spectral sequence\[
E_{2}^{i,j}=\mathrm{Tor}_{-i}^{A(\Omega)}(H^{j}(K^{p},\mathcal{D}_{\Omega})_{\leq h},A(\Sigma))\Rightarrow H^{i+j}(K^{p},\mathcal{D}_{\Sigma})_{\leq h}.\]
In addition, there are analogous spectral sequences relating Borel-Moore
homology with compactly supported cohomology, and boundary homology
with boundary cohomology, and there are morphisms between the spectral
sequences compatible with the morphisms between these different cohomology
theories. Finally, the spectral sequences and the morphisms between
them are equivariant for the natural Hecke actions on their $E_{2}$
pages and abutments; more succinctly, they are spectral sequences
of $\mathbf{T}(K^{p})$-modules.}

When no ambiguity is likely, we will refer to the two spectral sequences
of Theorem 1.1 as {}``the Ext spectral sequence'' and {}``the Tor
spectral sequence.'' We wish to emphasize that, unlike most familiar
spectral sequences of universal coefficients type, these sequences
have no particular reason to degenerate. However, the payoff of understanding
their $E_{2}$ pages is very great: this knowledge is enough to calculate
the dimensions of eigenvarieties. With the spectral sequences in hand,
we prove the following basic result.

\textbf{Theorem 1.2.}\emph{ Fix a slope datum $(U_{t},\Omega,h)$.}

\emph{i. For any $i$, $H_{i}(K^{p},\mathcal{A}_{\Omega})_{\leq h}$
is a faithful $A(\Omega)$-module if and only if $H^{i}(K^{p},\mathcal{D}_{\Omega})_{\leq h}$
is faithful.}

\emph{ii. If $G$ has $\mathbf{Q}$-anisotropic derived group, $H_{i}(K^{p},\mathcal{A}_{\Omega})_{\leq h}$
and $H^{i}(K^{p},\mathcal{D}_{\Omega})_{\leq h}$ are torsion $A(\Omega)$-modules
for all $i$, unless $G^{\mathrm{der}}(\mathbf{R})$ has a discrete
series, in which case they are torsion for all $i\neq\frac{1}{2}\mathrm{dim}G(\mathbf{R})/K_{\infty}Z_{\infty}$.}

Some precursors of Theorem 1.2.ii were discovered by Hida (Theorem
5.2 of \cite{HidordGL2}, about which more below) and Ash-Pollack-Stevens \cite{APSrigid}.
We should note that in proving the non-faithfulness of $H^{i}(K^{p},\mathcal{D}_{\Omega})_{\leq h}$
in the discrete series case for $i\neq\frac{1}{2}\mathrm{dim}G(\mathbf{R})/K_{\infty}Z_{\infty}$,
we make crucial use of the homology module $H_{\ast}(K^{p},\mathcal{A}_{\Omega})_{\leq h}$
and the Ext spectral sequence, and \emph{we don't know how to prove
this result by working solely with $\mathcal{D}_{\Omega}$.} This
is a basic example of the benefit of treating the modules $\mathcal{A}_{\Omega}$
and $\mathcal{D}_{\Omega}$ on an equal footing.

Before turning to explain our motivations, let us say a word about
the proof of the Ext spectral sequence of Theorem 1.1, the Tor spectral
sequence being very similar. Given a commutative ring $R$, an ideal
$\mathfrak{a}\subset R$, an $R$-module $M$, and a chain complex
$C_{\bullet}$ of projective $R$-modules, there is a quite general
convergent spectral sequence\[
E_{2}^{i,j}=\mathrm{Ext}_{R}^{i}(H_{j}(C_{\bullet}),R/\mathfrak{a})\Rightarrow H^{i+j}(\mathrm{Hom}_{R}(C_{\bullet},R/\mathfrak{a})).\]
This result, however, isn't directly applicable towards proving Theorem
1.1 because we don't know whether or not $\mathbf{A}_{\Omega}^{s}$
is a projective $A(\Omega)$-module.%
\footnote{This boils down to the question, seemingly open, of whether $\mathbf{Q}_{p}\left\langle X,Y\right\rangle $
is projective over $\mathbf{Q}_{p}\left\langle X\right\rangle $.%
} Fortunately, the theory of slope decompositions implies that $C_{\bullet}(K^{p},\mathcal{A}_{\Omega})_{\leq h}$
is a complex of finite flat and hence projective $A(\Omega)$-modules,
from which Theorem 1.1 follows easily \emph{except }for the assertion
regarding Hecke equivariance. This latter property is absolutely crucial
for applications, since it allows one to localize the entire spectral
sequence at any ideal in the Hecke algebra. We proceed by constructing
a homotopy action of $\mathbf{T}(K^{p})$ directly on the complex
$C_{\bullet}(K^{p},\mathcal{A}_{\Omega})_{\leq h}$ and making use
of a derived-categorical point of view. Granted the preliminary materials
of §2, the proof of Theorem 1.1 is actually rather short, and we refer
the reader to §3 for details.

\subsection{The geometry of eigenvarieties}

Our main motivation for Theorem 1.1 is the problem, broadly speaking,
of analyzing the geometry of eigenvarieties. To describe the spaces
we have in mind, fix a controlling operator $U_{t}$. For any slope
datum $(\Omega,h)$ we define $\mathbf{T}_{\Omega,h}(K^{p})$ as the
finite commutative $A(\Omega)$-algebra generated by the image of
$\mathbf{T}(K^{p})\otimes_{\mathbf{Q}_{p}}A(\Omega)$ in $\mathrm{End}_{A(\Omega)}\left(H^{\ast}(K^{p},\mathcal{D}_{\Omega})_{\leq h}\right)$.
Given a Noetherian $\mathbf{Q}_{p}$-algebra $A$, an $A$-module
$M$, and a $\mathbf{Q}_{p}$-algebra homomorphism $\phi:A\to\overline{\mathbf{Q}_{p}}$,
we define the $\phi$-generalized eigenspace $M_{(\phi)}$ of $M$
as \[
M_{(\phi)}=\lim_{n\to\infty}M[\mathfrak{M}^{n}]\]
where $\mathfrak{M}=\ker\phi$.

\textbf{Theorem 1.3. }\emph{There is a separated $\mathbf{Q}_{p}$-rigid
analytic space $\mathscr{X}=\mathscr{X}(K^{p})$ equipped with a (not
necessarily surjective) morphism $w:\mathscr{X}\to\mathcal{W}$, locally
finite in the domain over its image, with the following properties:}

\emph{i. There is a canonical algebra homomorphism $\psi:\mathbf{T}(K^{p})\to\mathcal{O}(\mathscr{X})$,
and the image of any double coset operator $[K^{p}gK^{p}]$ lies in
the subring $\mathcal{O}(\mathscr{X})^{\leq1}$ of power-bounded functions.}

\emph{ii. Given any $\mathbf{Q}_{p}$-algebra homomorphism $\phi:\mathbf{T}(K^{p})\to\overline{\mathbf{Q}_{p}}$,
the $\phi$-generalized eigenspace of $H^{\ast}(K^{p},\mathcal{D}_{\lambda})_{<\infty}$
is nonzero if and only if there is a point $x\in\mathscr{X}$ with
$w(x)=\lambda$ such that the diagram\[
\xymatrix{\mathbf{T}(K^{p})\ar[rr]^{\phi}\ar[dr]_{\psi(x)} &  & \overline{\mathbf{Q}_{p}}\\
 & \kappa(x)\ar[ur]_{\iota}}
\]
commutes for some embedding $\iota:\kappa(x)\to\overline{\mathbf{Q}_{p}}$.
(Here $\kappa(x)$ is the residue field of $ $$\mathcal{O}_{\mathscr{X},x}$,
and $\psi(x)$ is the image of $\psi$ under $\mathcal{O}(\mathscr{X})\to\mathcal{O}_{\mathscr{X},x}\to\kappa(x)$.)}

\emph{iii. The affinoid rigid spaces $\mathscr{X}_{\Omega,h}=\mathrm{Sp}\mathbf{T}_{\Omega,h}(K^{p})$
admissibly cover $\mathscr{X}$ as $(\Omega,h)$ runs over all slope
data.}

\emph{iv. For each degree $n$, the $\mathbf{T}_{\Omega,h}(K^{p})$-modules
$H^{n}(K^{p},\mathcal{D}_{\Omega})_{\leq h}$ glue together over the
aforementioned admissible covering of $\mathscr{X}(K^{p})$ into a
coherent $\mathcal{O}_{\mathscr{X}}$-module sheaf $\mathscr{M}^{n}$.}

\emph{v. The space $\mathscr{X}$ is independent of the choice of
controlling operator $U_{t}$ used in its construction.}

In other words, $\mathscr{X}$ analytically parametrizes the finite-slope
eigenpackets occuring in $H^{\ast}(K^{p},\mathcal{D}_{\lambda})$
as $\lambda$ varies over weight space. We sketch the construction
of $\mathscr{X}$ in §3.2 - we stress, however, that \emph{none} of
the key ideas in this construction are due to us: they are due to
Coleman-Mazur, Buzzard, and Ash-Stevens. We also note that for some
groups, we can describe a canonical rigid subgroup $\mathcal{W}_{0}\subset\mathcal{W}$,
the \emph{space of null weights}, together with a group-theoretic
splitting $\mathcal{W}\to\mathcal{W}_{0}$ and a universal character
$\mathbf{T}(K^{p})\to\mathcal{O}(\mathcal{W}/\mathcal{W}_{0})$ such
that $\mathscr{X}$ is a $\mathcal{W}/\mathcal{W}_{0}$-bundle over
$\mathscr{X}_{0}=\mathscr{X}\cap w^{-1}(\mathcal{W}_{0})$. Restricting
to the space $\mathscr{X}_{0}$ amounts to factoring out wild twists,
and is analogous to the Sen null space in Galois deformation theory.
When $G=\mathrm{GL}_{n}/\mathbf{Q}$ and $T=\mathrm{diag}(t_{1},\dots,t_{n})$
is the standard maximal torus, $\mathcal{W}_{0}$ is simply the space
of characters trivial on the one-parameter subgroup $\mathrm{diag}(1,\dots,1,t_{n})$.
The first sign the construction described in Theorem 1.3 is worth
consideration is probably the following theorem, which is due to Stevens \cite{Stfamilies}
and Bellaïche \cite{BeCrit}.

\textbf{Proposition 1.4. }\emph{When $G=\mathrm{GL}_{2}$ and $K^{p}=\mathrm{GL}_{2}(\mathbf{A}_{f}^{p})\cap K_{1}(N)$,
$\mathscr{X}_{0}$ is isomorphic to the Coleman-Mazur-Buzzard eigencurve
$\mathscr{C}(N)$ of tame level $N$.}

In general, when $G^{\mathrm{der}}(\mathbf{R})$ has a discrete series,
standard limit multiplicity results yield an abundance of classical
automorphic forms of essentially every arithmetic weight, and one
expects that correspondingly every irreducible component $\mathscr{X}_{i}$
of the eigenvariety $\mathscr{X}$ which contains a {}``suitably
general'' classical point has maximal dimension, namely $\mathrm{dim}\mathscr{X}_{i}=\mathrm{dim}\mathcal{W}=\mathrm{rank}G$.
This numerical coincidence is\emph{ characteristic }of the groups
for which $G^{\mathrm{der}}(\mathbf{R})$ has a discrete series. More
precisely, define the \emph{defect }and the \emph{amplitude }of $G$,
respectively, as the integers $l(G)=\mathrm{rank}G-\mathrm{rank}K_{\infty}Z_{\infty}$
and $q(G)=\frac{1}{2}(\mathrm{dim}(G(\mathbf{R})/K_{\infty}Z_{\infty})-l(G))$.
Note that $l(G)$ is zero if and only if $G^{\mathrm{der}}(\mathbf{R})$
has a discrete series, and that algebraic representations with regular
highest weight contribute to $(\mathfrak{g},K_{\infty})$-cohomology
exactly in the unbroken range of degrees $[q(G),q(G)+l(G)]$. We say
a point $x\in\mathscr{X}(K^{p})$ is \emph{classical }if the induced
eigenpacket $\psi(x):\mathbf{T}(K^{p})\to\kappa(x)$ matches the Hecke
data of an algebraic automorphic representation $\pi$ of $G(\mathbf{A_{Q}})$,
and $x$ is \emph{regular }if $\pi_{\infty}$ contributes to $(\mathfrak{g},K_{\infty})$-cohomology
with coefficients in an irreducible algebraic representation with
regular highest weight. The definition of a\emph{ non-critical }classical
point is rather more subtle and we defer it until §2.4. The following
conjecture is a special case of a conjecture of Urban (Conjecture
5.7.3 of \cite{UrEigen}).

\textbf{Conjecture 1.5. }\emph{Every irreducible component $\mathscr{X}_{i}$
of $\mathscr{X}(K^{p})$ containing a cuspidal non-critical regular
classical point has dimension $\mathrm{dim}\mathcal{W}-l(G)$.}

To put this conjecture in context, we note that the inequality $\mathrm{dim}\mathscr{X}_{i}\geq\mathrm{dim}\mathcal{W}-l(G)$
is known. The idea for such a result is independently due to Stevens
and Urban, whose proofs remain unpublished. On reading an earlier
version of this paper, James Newton discovered a very short and elegant
proof of this inequality, which is given here in Appendix B.

In any case, using the spectral sequences, we verify Conjecture 1.5
in many cases.

\textbf{Theorem 1.6.}
\begin{description}
\item [{i.}] \emph{If $l(G)=0$, then Conjecture 1.5 is true, and if $x\in\mathscr{X}(K^{p})$
is a cuspidal non-critical regular classical point, then $\mathscr{X}$
is smooth at $x$ and the weight morphism $w$ is étale at $x$.}
\item [{ii.}] \emph{If $l(G)=1$, then Conjecture 1.5 is true. Furthermore,
if $l(G)\geq1$, then every irreducible component of $\mathscr{X}(K^{p})$
containing a cuspidal non-critical regular classical point has dimension
at most $\mathrm{dim}\mathcal{W}-1$.}
\end{description}

Some basic examples of groups with $l(G)=1$ include $\mathrm{GL}_{3}/\mathbf{Q}$
and $\mathrm{Res}_{F/\mathbf{Q}}H$ where $F$ is a number field with
exactly one complex embedding and $H$ is an $F$-inner form of $\mathrm{GL}_{2}/F$.
In particular, we have the following corollary, which was our original
motivation for this project.

\textbf{Corollary 1.7. }\emph{Let $F$ be an imaginary quadratic extension
of $\mathbf{Q}$, and let $G=\mathrm{Res}_{F/\mathbf{Q}}H$ where
$H$ is an $F$-inner form of $\mathrm{GL}_{2}/F$ (possibly the split
form). Then there is a two-dimensional space of null weights $\mathcal{W}_{0}$,
and any component of the eigenvariety $\mathscr{X}_{0}$ containing
a cuspidal non-critical regular classical point is a rigid analytic
curve.}

In the ordinary case, Corollary 1.7 is a result of Hida (§5-6 of \cite{HidordGL2}).
In fact, we are able to give a much more precise description of the
geometry of the eigenvariety in this case.

\textbf{Theorem 1.8. }\emph{Maintaining the notation and assumptions
of Corollary 1.7, the eigenvariety $\mathscr{X}_{0}$ is naturally
a union of subspaces $ $$\mathscr{X}_{0}^{\mathrm{punc}}\cup\mathscr{X}_{0}^{\mathrm{Eis}}\cup\mathscr{X}_{0}^{\mathrm{cusp}}$
where:}

\emph{i. $\mathscr{X}_{0}^{\mathrm{punc}}$ is zero-dimensional, supported
over the trivial weight in $\mathcal{W}_{0}$. The eigenpackets carried
by $\mathscr{X}_{0}^{\mathrm{punc}}$ are of the form $T_{\mathfrak{q}}\mapsto(1+\mathbf{N}\mathfrak{q})\eta(\mathfrak{q})$,
where $\eta$ is an everywhere-unramified Hecke character of $F$.}

\emph{ii. $\mathscr{X}_{0}^{\mathrm{Eis}}$ is empty unless $H$ is
$F$-split, in which case $\mathscr{X}_{0}^{\mathrm{Eis}}$ is finite
and flat over $\mathcal{W}_{0}$. The fiber of $\mathscr{X}_{0}^{\mathrm{Eis}}$
over $\lambda\in\mathcal{W}_{0}$ carries the eigenpackets $T_{\mathfrak{q}}\mapsto\mathbf{N}\mathfrak{q}\cdot\lambda(\mathfrak{q})\eta_{1}(\mathfrak{q})+\eta_{2}(\mathfrak{q})$,
where $\eta_{1}$ and $\eta_{2}$ are certain finite-order Hecke characters
of $F$ of conductors dividing $\mathfrak{n}$.}

\emph{iii. $\mathscr{X}_{0}^{\mathrm{cusp}}$ is equidimensional of
dimension one.}

The proofs of Theorem 1.2 and Theorem 1.6 systematically play the
two spectral sequences against each other, and are curiously intertwined
with the details of the \emph{construction} of $\mathscr{X}(K^{p})$;
in particular, we make crucial use of a non-canonical intermediate
rigid analytic space $\mathcal{Z}$ which depends heavily on the chain
complex $C_{\bullet}(K^{p},-)$. We also appeal to a simple but powerful
analytic continuation principle: if $ $$\mathscr{X}$ is an equidimensional
rigid analytic space and $\mathscr{F}$ is a coherent analytic sheaf
over $\mathscr{X}$ such that $\mathrm{Supp}\mathscr{F}$ contains
\emph{some }admissible open subset, then $\mathrm{Supp}\mathscr{F}$
contains an entire irreducible component of $\mathscr{X}$.

\subsection{The relation to other incarnations of \emph{p}-adic automorphic forms}

It seems pertinent for us to make some remarks on the relation between
overconvergent cohomology and other constructions of \emph{p-}adic
automorphic forms and eigenvarieties.

Aside from overconvergent cohomology, the other main approach to a
{}``cohomological'' construction of eigenvarieties is Emerton's
completed cohomology \cite{Eminterpolate,CEcompletedsurvey}. We briefly
review the main definitions. Fix a tame level $K^{p}$ and an open
compact subgroup $K_{p}\subset G(\mathbf{Q}_{p})$, and choose a filtration
$K_{p}\supset K_{p}^{1}\supset K_{p}^{2}\supset\dots\supset K_{p}^{i}\supset\dots$
of $K_{p}$ by open normal subgroups such that $\cap_{i=1}^{\infty}K_{p}^{i}=\{1\}$.
We make the definitions\[
\widetilde{H}_{i}=\widetilde{H}_{i}(K^{p})=\lim_{\infty\leftarrow j}H_{i}(G(\mathbf{Q})\backslash G(\mathbf{A})/K_{p}^{j}K^{p}K_{\infty}Z_{\infty},\mathbf{Z}_{p})\]
and\[
\widetilde{H}^{i}=\widetilde{H}^{i}(K^{p})=\lim_{\infty\leftarrow k}\lim_{j\to\infty}H^{i}(G(\mathbf{Q})\backslash G(\mathbf{A})/K_{p}^{j}K^{p}K_{\infty}Z_{\infty},\mathbf{Z}/p^{k}).\]
We write $\widetilde{H}_{i}^{BM}$ and $\widetilde{H}_{c}^{i}$ for
the obvious Borel-Moore and compactly supported variants. The natural
$K_{p}$-action on $\widetilde{H}^{i}$ and $\widetilde{H}_{i}$ extends
to a continuous $G(\mathbf{Q}_{p})$-action which is \emph{independent
}of the choice of $K_{p}$ (explaining our omission of $K_{p}$ from
the notation). For each $n$, there are short exact sequences of continuous
$G(\mathbf{Q}_{p})$-modules \[
0\to\mathrm{Hom}^{\mathrm{cts}}(\widetilde{H}_{n-1},\mathbf{Q}_{p}/\mathbf{Z}_{p})\to\widetilde{H}^{n}\to\mathrm{Hom}^{\mathrm{cts}}(\widetilde{H}_{n},\mathbf{Z}_{p})\to0\]
and\[
0\to\mathrm{Hom}(\widetilde{H}^{n+1}[p^{\infty}],\mathbf{Q}_{p}/\mathbf{Z}_{p})\to\widetilde{H}_{n}\to\mathrm{Hom}^{\mathrm{cts}}(\widetilde{H}^{n},\mathbf{Z}_{p})\to0.\]
In addition, for any fixed choice of $K_{p}$, the $K_{p}$-action
on $\widetilde{H}_{i}$ extends to a left action of the (Noetherian)
completed group ring $\Lambda=\mathbf{Z}_{p}[[K_{p}]]$ which gives
$\widetilde{H}_{i}$ the structure of a finitely presented $\Lambda$-module.
Finally, there is a spectral sequence of the form\[
E_{2}^{i,j}=\mathrm{Ext}_{\Lambda}^{i}(\widetilde{H}_{j},\Lambda)\Rightarrow\widetilde{H}_{d-i-j}^{BM}\]
where $d=\mathrm{dim}G(\mathbf{R})/K_{\infty}Z_{\infty}$.

The following theorem and conjecture are due to Calegari and Emerton.

\textbf{Theorem 1.9. }\emph{Suppose $G$ is semisimple. Then $\widetilde{H}_{i}$
is a torsion $\Lambda$-module unless $G(\mathbf{R})$ has a discrete
series and $i=\frac{1}{2}\mathrm{dim}G(\mathbf{R})/K_{\infty}Z_{\infty}$,
in which case $\widetilde{H}_{i}$ is a faithful $\Lambda$-module.}

\textbf{Conjecture 1.10. }\emph{Define the }codimension \emph{of a
$\Lambda$-module $M$, written $\mathrm{codim}(M)$ as the least
integer $i$ such that $\mathrm{Ext}_{\Lambda}^{i}(M,\Lambda)\neq0$.
Then $\mathrm{codim}(\widetilde{H}_{q(G)})=l(G)$ and $\mathrm{codim}(\tilde{H}_{i})>l(G)$
for $i\neq q(G)$.}

The reader will not fail to notice the very strong formal similarities
between Theorem 1.2 and Theorem 1.9, and between Conjecture 1.5 and
Conjecture 1.10. It seems very likely, to borrow a phrase from \cite{CEcompletedsurvey},
that this relation {}``is more than one of mere analogy.'' Some
evidence for this has begun to accumulate: in a forthcoming paper \cite{Hcomp}
we will construct Hecke-equivariant morphisms \[
\phi^{n}:H^{n}(K^{p},\mathcal{A}_{\lambda})_{<\infty}\to J_{B}^{T(\mathbf{Z}_{p})=\lambda}(\widetilde{H}^{n}(K^{p})_{\mathrm{la}}),\]
where $\lambda\in\mathcal{W}(\overline{\mathbf{Q}_{p}})$ is arbitrary,
$\mathcal{A}_{\lambda}=\mathcal{A}_{\Omega}\otimes_{A(\Omega)}A(\Omega)/\mathfrak{m}_{\lambda}$,
$J_{B}^{T(\mathbf{Z}_{p})=\lambda}(-)$ denotes the $\lambda$-isotypic
subspace of Emerton's locally analytic Jacquet module, and $(-)_{\mathrm{la}}$
denotes the functor of passage to $\mathbf{Q}_{p}$-locally analytic
vectors defined in \cite{STalgadm}. In many cases the morphisms $\phi^{n}$
are injective. Furthermore, the $\phi^{n}$'s are defined as the edge
maps of a spectral sequence - very different from those considered
in this paper - which computes $H^{\ast}(K^{p},\mathcal{A}_{\lambda})$
using the continuous cohomology of the modules $\widetilde{H}^{\ast}(K^{p})_{\mathrm{la}}$.

On the other hand, the completed cohomology groups should contain
much more data than overconvergent cohomology: completed cohomology
is expected to pick up essentially every Galois representation, while
overconvergent cohomology should only pick up the trianguline Galois
representations. Thus it is not obvious, at least to this author,
how direct a connection between Conjectures 1.5 and 1.10 might be
expected. In any case, the two theories seem to play complementary
roles: completed cohomology is amenable to the powerful methods of
\emph{p}-adic analytic representation theory, while overconvergent
cohomology with its theory of slope decompositions is planted firmly
in the world of finitely generated modules over Noetherian rings.

Finally, despite the title of this paper, the reader may have noticed
a conspicuous absence of overconvergent modular forms in our discussion.
Let us simply say that when $G$ gives rise to Shimura varieties,
we expect that the Hecke data occuring in spaces of overconvergent
modular forms for $G$ will also occur in overconvergent cohomology,
and that Chenevier's interpolation theorem (Theorémè 1 of \cite{CHjacquetlanglands})
will allow a fairly straightforward verification of this expectation
whenever a {}``small slope forms are classical''-type theorem is
available.

\subsection*{Notation and terminology}

Our notation and terminology is mostly standard. In nonarchimedian
functional analysis and rigid analytic geometry we essentially follow \cite{BGR}.
In the body of the paper, $k$ denotes an extension field of $\mathbf{Q}_{p}$,
complete for its norm $\left|\bullet\right|_{k}$. If $M$ and $N$
are $k$-Banach spaces, we write $\mathcal{L}_{k}(M,N)$ for the space
of continuous $k$-linear maps between $M$ and $N$; the operator
norm\[
|f|=\sup_{m\in M,\,|m|_{M}\leq1}|f(m)|_{N}\]
makes $\mathcal{L}_{k}(M,N)$ into a $k$-Banach space. If $(A,|\bullet|_{A})$
is a $k$-Banach space which furthermore is a commutative Noetherian
$k$-algebra whose multiplication map is (jointly) continuous, we
say $A$ is a $k$-\emph{Banach algebra}. An $A$-module $M$ which
is also a $k$-Banach space is a \emph{Banach $A$-module }if the
structure map $A\times M\to M$ extends to a continuous map $A\widehat{\otimes}_{k}M\to M$,
or equivalently if the norm on $M$ satisfies $\left|am\right|_{M}\leq C|a|_{A}|m|_{M}$
for all $a\in A$ and $m\in M$ with some fixed constant $C$. For
a topological ring $R$ and topological $R$-modules $M,N$, we write
$\mathcal{L}_{R}(M,N)$ for the $R$-module of continuous $R$-linear
maps $f:M\to N$. When $A$ is a $k$-Banach algebra and $M,N$ are
Banach $A$-modules, we topologize $\mathcal{L}_{A}(M,N)$ via its
natural Banach $A$-module structure. We write $\mathrm{Ban}_{A}$
for the category whose objects are Banach $A$-modules and whose morphisms
are elements of $\mathcal{L}_{A}(-,-)$. If $I$ is any set and $A$
is a $k$-Banach algebra, we write $c_{I}(A)$ for the module of sequences
$\mathbf{a}=(a_{i})_{i\in I}$ with $|a_{i}|_{A}\to0$; the norm $|\mathbf{a}|=\sup_{i\in I}|a_{i}|_{A}$
gives $c_{I}(A)$ the structure of a Banach $A$-module. If $M$ is
any Banach $A$-module, we say $M$ is \emph{orthonormalizable }if
$M$ is \emph{isomorphic} to $c_{I}(A)$ for some $I$ (such modules
are called {}``potentially orthonormalizable'' in  \cite{BuEigen}).

If $A$ is an affinoid algebra, then $\mathrm{Sp}A$, the \emph{affinoid
space }associated with $A$, denotes the locally G-ringed space $(\mathrm{Max}A,\mathcal{O}_{A})$
where $\mathrm{Max}A$ is the set of maximal ideals of $A$ endowed
with the Tate topology and $\mathcal{O}_{A}$ is the extension of
the assignment $U\mapsto A_{U}$, for affinoid subdomains $U\subset\mathrm{Max}A$
with representing algebras $A_{U}$, to a structure sheaf on $\mathrm{Max}A$.
If $X$ is an affinoid space, we write $A(X)$ for the coordinate
ring of $X$. If $A$ is reduced we equip $A$ with the canonical
supremum norm. If $X$ is a rigid analytic space, we write $\mathcal{O}_{X}(X)$
or $\mathcal{O}(X)$ for the ring of global sections of the structure
sheaf on $X$. Given a point $x\in X$, we write $\mathfrak{m}_{x}$
for the corresponding maximal ideal in $\mathcal{O}_{X}(U)$ for any
admissible affinoid open $U\subset X$ containing $x$, and $\kappa(x)$
for the residue field $\mathcal{O}_{X}(U)/\mathfrak{m}_{x}$; $\mathcal{O}_{X,x}$
denotes the local ring of $\mathcal{O}_{X}$ at $x$ in the Tate topology,
and $\widehat{\mathcal{O}_{X,x}}$ denotes the $\mathfrak{m}_{x}$-adic
completion of $\mathcal{O}_{X,x}$.

In homological algebra our conventions follow \cite{Weibelhomalg}.
If $R$ is a ring, we write $\mathbf{K}^{?}(R)$, $?\in\{+,-,b,\emptyset\}$
for the homotopy category of $?$-bounded $R$-module complexes and
$\mathbf{D}^{?}(R)$ for its derived category.

\subsection*{Acknowledgments}

This paper represents a portion of my forthcoming Ph.D. dissertation
at Boston College. First and foremost, I'm grateful to my advisor,
Avner Ash, for giving very generously of his time, insight, and enthusiasm.
During the development of these ideas, I also enjoyed stimulating
conversations with Joël Bellaïche, Jay Pottharst, Glenn Stevens, and
Jack Thorne; it's a pleasure to acknowledge their help and influence
here. I'm particularly grateful to Glenn Stevens for emphasizing the
utility of slope decompositions to me, and for his kind encouragement.
Kevin Buzzard, James Newton, and Jack Thorne offered helpful comments
and corrections on earlier drafts, for which I thank them sincerely.

\section{Background material}

In §2.1-§2.3 we lay down some foundational notation and definitions;
the ideas in these sections are entirely due to Ash and Stevens \cite{AS},
though our presentation differs somewhat in insignificant details.

\subsection{Algebraic groups and overconvergent coefficient modules}

\subsubsection*{Lie theoretic data}

Fix a prime $p$ and a connected, reductive $\mathbf{Q}$-group $G$;
we suppose $G/\mathbf{Q}_{p}$ is split and is the generic fiber of
a smooth group scheme $\mathcal{G}/\mathbf{Z}_{p}$. Fix a Borel $B=TN$,
an opposite Borel $B^{\mathrm{opp}}=TN^{\mathrm{opp}}$, and a compatible
Iwahori subgroup $I\subset\mathcal{G}(\mathbf{Z}_{p})\subset G(\mathbf{Q}_{p})$.
Set $X^{\ast}=\mathrm{Hom}(T,\mathbf{G}_{m})$ and $X_{\ast}=\mathrm{Hom}(\mathbf{G}_{m},T)$,
and let $\Phi$ and $\Phi^{+}$ be the sets of roots and positive
roots, respectively, for the Borel $B$. We write $X_{+}^{\ast}$
for the cone of $B$-dominant weights; $\rho\in X^{\ast}\otimes_{\mathbf{Z}}\frac{1}{2}\mathbf{Z}$
denotes half the sum of the positive roots.

Set $N^{\circ}=\left\{ n\in N^{\mathrm{opp}}(\mathbf{Z}_{p}),\, n\equiv1\,\mathrm{in}\, G(\mathbf{Z}/p\mathbf{Z})\right\} $,
so the Iwahori decomposition reads $I=N^{\circ}\cdot T(\mathbf{Z}_{p})\cdot N(\mathbf{Z}_{p})$.
For any integer $c\geq1$, we set \[
I_{0}^{c}=\left\{ g\in I,\, g\,\mathrm{mod}\, p^{c}\in B(\mathbf{Z}/p^{c}\mathbf{Z})\right\} \]
and\[
I_{1}^{c}=\left\{ g\in I,\, g\,\mathrm{mod}\, p^{c}\in N(\mathbf{Z}/p^{c}\mathbf{Z})\right\} .\]
Note that $I_{1}^{c}$ is normal in $I_{0}^{c}$, with quotient $T(\mathbf{Z}/p^{c}\mathbf{Z})$.
For $s\geq1$ a positive integer, define\[
I^{s}=\left\{ g\in I|g\,\equiv1\,\mathrm{in}\, G(\mathbf{Z}/p^{s}\mathbf{Z})\right\} \]
and $T^{s}=T(\mathbf{Z}_{p})\cap I^{s}$, $N^{s}=N(\mathbf{Z}_{p})\cap I^{s}$.
Note that $N^{s}$ is normal in $B(\mathbf{Z}_{p})$ and in $N(\mathbf{Z}_{p})$. 

We define semigroups of $T(\mathbf{Q}_{p})$ by\[
\Lambda=\left\{ t\in T(\mathbf{Q}_{p}),\, t^{-1}N(\mathbf{Z}_{p})t\subseteq N(\mathbf{Z}_{p})\right\} \]
and \[
\Lambda^{+}=\left\{ t\in T(\mathbf{Q}_{p}),\,\bigcap_{i=1}^{\infty}t^{-i}N(\mathbf{Z}_{p})t^{i}=\{1\}\right\} .\]
A simple calculation shows that $t\in T(\mathbf{Q}_{p})$ is contained
in $\Lambda$ (resp. $\Lambda^{+}$) if $v_{p}(\alpha(t))\leq0$ (resp.
$v_{p}(\alpha(t))<0$) for all $\alpha\in\Phi^{+}$. Using these,
we define semigroups of $G(\mathbf{Q}_{p})$ by $\Delta=I\Lambda I,\;\Delta^{+}=I\Lambda^{+}I$.
Note the inclusions $I\subset\Delta\supset\Delta^{+}$. The Iwahori
decompsition extends to $\Delta$: any element $g\in\Delta$ has a
unique decomposition $g=\mathrm{n}^{\circ}(g)\mathrm{t}(g)\mathrm{n}(g)$
with $\mathrm{n}^{\circ}\in N^{\circ}$, $\mathrm{t}\in\Lambda$,
$\mathrm{n}\in N(\mathbf{Z}_{p})$. Fix once and for all a group homomorphism
$\sigma:T(\mathbf{Q}_{p})\to T(\mathbf{Z}_{p})$ which splits the
inclusion $T(\mathbf{Z}_{p})\subset T(\mathbf{Q}_{p})$. Note that
$\sigma$ splits $\Lambda$ as the direct product $T(\mathbf{Z}_{p})\cdot\left(\Lambda\cap\ker\sigma\right)$
via $\delta\mapsto(\sigma(\delta),\delta\cdot\sigma(\delta)^{-1})$.

We fix an analytic isomorphism $\psi:N(\mathbf{Z}_{p})\simeq\mathbf{Z}_{p}^{d}$,
$d=\mathrm{dim}N$, which identifies cosets of $N^{s}$ with additive
cosets $\mathbf{a}+p^{s}\mathbf{Z}_{p}^{d}\subset\mathbf{Z}_{p}^{d}$,
$\mathbf{a}\in\mathbf{Z}_{p}^{d}$.

\textbf{Definition. }\emph{If $R$ is any $\mathbf{Q}_{p}$-Banach
algebra and $s$ is a nonnegative integer, the module $\mathbf{A}(N(\mathbf{Z}_{p}),R)^{s}$
of} $s$-locally analytic $R$-valued functions on $N(\mathbf{Z}_{p})$\emph{
is the $R$-module of continuous functions $f:N(\mathbf{Z}_{p})\to R$
such that \[
f_{\mathbf{a}}=f\left(\psi^{-1}\left(\mathbf{a}+p^{s}(x_{1},\dots,x_{d})\right)\right):\mathbf{Z}_{p}^{d}\to R\]
is given by an element of the $d$-variable Tate algebra $R\left\langle x_{1},\dots,x_{d}\right\rangle $
for any fixed $\mathbf{a}\in\mathbf{Z}_{p}^{d}$.}

Given $f\in\mathbf{A}(N(\mathbf{Z}_{p}),R)^{s}$ and $\mathbf{i}=(i_{1},\dots,i_{d})\in\mathbf{N}^{d}$,
define coefficients $c(f_{\mathbf{a}},\mathbf{i})\in R$ by $f_{\mathbf{a}}(x_{1},\dots,x_{d})=\sum_{\mathbf{i}\in\mathbb{N}^{d}}c(f_{\mathbf{a}},\mathbf{i})x_{1}^{i_{1}}\cdots x_{d}^{i_{d}}$.
The norm\[
\left\Vert f\right\Vert =\mathrm{sup}_{\mathbf{a}\in\mathbf{Z}_{p}^{d}}\mathrm{sup}_{\mathbf{i}\in\mathbf{N}^{d}}|c(f_{\mathbf{a}},\mathbf{i})|_{R}\]
defines a Banach $R$-module structure on $\mathbf{A}(N(\mathbf{Z}_{p}),R)^{s}$,
with respect to which the canonical inclusion $\mathbf{A}(N(\mathbf{Z}_{p}),R)^{s}\subset\mathbf{A}(N(\mathbf{Z}_{p}),R)^{s+1}$
is compact.

\subsubsection*{Weights and Modules}

The space of weights is the space $\mathcal{W}=\mathrm{Hom}_{\mathrm{cts}}(T(\mathbf{Z}_{p}),\mathbf{G}_{m})$;
we briefly recall the rigid analytic structure on $\mathcal{W}$,
following the discussion in §3.4-3.5 of \cite{AS}. For any rational
number $r$, we define\[
\mathbf{Q}_{p}\left\langle p^{r}X\right\rangle =\left\{ \sum_{n=0}^{\infty}a_{n}X^{n}\in\mathbf{Q}_{p}[[X]],\,\lim_{n\to\infty}p^{rn}|a_{n}|\to0\right\} ;\]
this is a $\mathbf{Q}_{p}$-affinoid algebra, and $\mathcal{B}[p^{r}]=\mathrm{Sp}\mathbf{Q}_{p}\left\langle p^{r}X\right\rangle $
is the $\mathbf{Q}_{p}$-rigid analytic disk of radius $p^{r}$. Choose
a monotone increasing sequence of rational numbers $\mathbf{r}=\{r_{i}\}_{i\in\mathbf{N}}$
with $r_{i}<0$ and $\lim_{i\to\infty}r_{i}=0$; the natural morphisms
$\mathbf{Q}_{p}\left\langle p^{r_{i+1}}X\right\rangle \to\mathbf{Q}_{p}\left\langle p^{r_{i}}X\right\rangle $
dualize to maps $\mathcal{B}[p^{r_{i}}]\to\mathcal{B}[p^{r_{i+1}}]$,
and we define \emph{the} $\mathbf{Q}_{p}$-\emph{rigid analytic open
unit disk $\mathcal{B}$ }as the natural gluing of the $\mathcal{B}[p^{r_{i}}]$'s
along these maps. The rigid structure on $\mathcal{B}$ is independent
of the choice of rational sequence $\mathbf{r}$. There is a natural
isomorphism $\mathcal{W}\simeq\widehat{T(\mathbf{Z}_{p})_{\mathrm{tors}}}\times\mathcal{B}^{d}$,
$d=\mathrm{dim}\mathbf{T}$, and we equip $\mathcal{W}$ with the
unique rigid structure for which this is an isomorphism of rigid analytic
spaces. 

Let $\Omega\subset\mathcal{W}$ be an admissible open affinoid subset.

\textbf{Lemma 2.1.1. }\emph{The ring $A(\Omega)$ is a regular ring.}

\emph{Proof. }By construction, $\mathcal{W}$ is a disjoint union
of finitely many open polydisks $\mathcal{B}^{d}$; write $\phi:\Omega\to\mathcal{B}^{d}\subset\mathcal{W}$
for the open immersion of $\Omega$ into whichever polydisk contains
it. For any point $x\in\mathcal{B}^{d}$, the complete local ring
$\widehat{\mathcal{O}_{\mathcal{B}^{d},x}}$ is a power series ring
in $d$ variables over $k_{x}$, so is regular. For any $\omega\in\Omega$
the open immersion $\phi$ induces an isomorphism $\phi^{\ast}:\widehat{\mathcal{O}_{\mathcal{B}^{d},\phi(\omega)}}\simeq\widehat{\mathcal{O}_{\Omega,\omega}}$,
so the completed local rings of $A(\Omega)$ at all maximal ideals
are regular. By Proposition 7.3.2/8 of \cite{BGR}, this implies that
each algebraic local ring $A(\Omega)_{\mathfrak{m}},\,\mathfrak{m}\in\mathrm{Max}A(\Omega)$
is regular, as desired. $\square$

We turn now to the key definitions of this section. Given $\Omega\subset\mathcal{W}$
admissible open, we write $\chi_{\Omega}:T(\mathbf{Z}_{p})\to A(\Omega)^{\times}$
for the unique character it determines. We define $s[\Omega]$ as
the minimal integer such that $\chi_{\Omega}|_{T^{s[\Omega]}}$ is
analytic. For any integer $s\geq s[\Omega]$, we make the definition\[
\mathbf{A}_{\Omega}^{s}=\left\{ f:I\to A(\Omega),\, f\,\mathrm{analytic\, on\, each\,}I^{s}-\mathrm{coset},\, f(n^{\circ}tg)=\chi_{\Omega}(t)f(g)\,\forall n^{\circ}\in N^{\circ},\, t\in T(\mathbf{Z}_{p}),\, g\in I\right\} .\]
By the Iwahori decomposition, restricting an element $f\in\mathbf{A}_{\Omega}^{s}$
to $N(\mathbf{Z}_{p})$ induces an isomorphism\begin{eqnarray*}
\mathbf{A}_{\Omega}^{s} & \simeq & \mathbf{A}(N(\mathbf{Z}_{p}),A(\Omega))^{s}\\
f & \mapsto & f|_{N(\mathbf{Z}_{p})},\end{eqnarray*}
and we regard $\mathbf{A}_{\Omega}^{s}$ as a Banach $A(\Omega)$-module
via pulling back the Banach module structure on $\mathbf{A}(N(\mathbf{Z}_{p}),A(\Omega))^{s}$
under this isomorphism. Right translation gives $\mathbf{A}_{\Omega}^{s}$
the structure of a continuous left $A(\Omega)[I]$-module. More generally,
the formula\[
(N^{\circ}b)\star\delta=N^{\circ}\sigma(\delta)\delta^{-1}b\delta,\, b\in B(\mathbf{Z}_{p})\simeq N^{\circ}\backslash I\,\mathrm{and}\,\delta\in\Lambda\]
yields a right action of $\Delta$ on $N^{\circ}\backslash I$ which
extends the natural right translation action by $I$ (cf.\emph{ }§2.5
of \cite{AS}) and induces a left $\Delta$-action on $\mathbf{A}_{\Omega}^{s}$
which we denote by $\delta\star f,\, f\in\mathbf{A}_{\Omega}^{s}$.
For any $\delta\in\Delta^{+}$, the image of the operator $\delta\star-\in\mathcal{L}_{A(\Omega)}(\mathbf{A}_{\Omega}^{s},\mathbf{A}_{\Omega}^{s})$
factors through the inclusion $\mathbf{A}_{\Omega}^{s-1}\hookrightarrow\mathbf{A}_{\Omega}^{s}$,
and so defines a completely continuous operator on $\mathbf{A}_{\Omega}^{s}$.
The Banach dual \begin{eqnarray*}
\mathbf{D}_{\Omega}^{s} & = & \mathcal{L}_{A(\Omega)}(\mathbf{A}_{\Omega}^{s},A(\Omega))\\
 & \simeq & \mathcal{L}_{A(\Omega)}(\mathbf{A}(N(\mathbf{Z}_{p}),k)^{s}\widehat{\otimes}_{k}A(\Omega),A(\Omega))\\
 & \simeq & \mathcal{L}_{k}(\mathbf{A}(N(\mathbf{Z}_{p}),k)^{s},A(\Omega))\end{eqnarray*}
inherts a dual right action of $\Delta$, and the operator $\delta\star-$
for $\delta\in\Delta^{+}$ likewise factors through the inclusion
$\mathbf{D}_{\Omega}^{s+1}\hookrightarrow\mathbf{D}_{\Omega}^{s}$.

We define \[
\mathcal{A}_{\Omega}=\lim_{s\to\infty}\mathbf{A}_{\Omega}^{s}\]
where the direct limit is taken with respect to the natural compact,
injective transition maps $\mathbf{A}_{\Omega}^{s}\to\mathbf{A}_{\Omega}^{s+1}$.
Note that $\mathcal{A}_{\Omega}$ is topologically isomorphic to the
module of $A(\Omega)$-valued locally analytic functions on $N(\mathbf{Z}_{p})$,
equipped with the finest locally convex topology for which the natural
maps $\mathbf{A}_{\Omega}^{s}\hookrightarrow\mathcal{A}_{\Omega}$
are continuous. The $\Delta$-actions on $\mathbf{A}_{\Omega}^{s}$
induce a continuous $\Delta$-action on $\mathcal{A}_{\Omega}$. Set\[
\mathcal{D}_{\Omega}=\left\{ \mu:\mathcal{A}_{\Omega}\to A(\Omega),\,\mu\,\mathrm{is}\, A(\Omega)-\mathrm{linear\, and\, continuous}\right\} ,\]
and topologize $\mathcal{D}_{\Omega}$ via the coarsest locally convex
topology for which the natural maps $\mathcal{D}_{\Omega}\to\mathbf{D}_{\Omega}^{s}$
are continuous. In particular, the canonical map\[
\mathcal{D}_{\Omega}\to\lim_{\infty\leftarrow s}\mathbf{D}_{\Omega}^{s}\]
is a topological isomorphism of locally convex $A(\Omega)$-modules,
and $\mathcal{D}_{\Omega}$ is compact and Fréchet. Note that the
transition maps $\mathbf{D}_{\Omega}^{s+1}\to\mathbf{D}_{\Omega}^{s}$
are \emph{injective}, so $\mathcal{D}_{\Omega}=\cap_{s\gg0}\mathbf{D}_{\Omega}^{s}$.

Suppose $\Sigma\subset\Omega$ is a Zariski closed subspace; by Corollary
9.5.2/8 of \cite{BGR}, $\Sigma$ arises from a surjection $A(\Omega)\twoheadrightarrow A(\Sigma)$
with $A(\Sigma)$ an affinoid algebra. We make the definitions $\mathbf{D}_{\Sigma}^{s}=\mathbf{D}_{\Omega}^{s}\otimes_{A(\Omega)}A(\Sigma)$
and $\mathcal{D}_{\Sigma}=\mathcal{D}_{\Omega}\otimes_{A(\Omega)}A(\Sigma)$.

\textbf{Proposition 2.1.2. }\emph{There are canonical topological
isomorphisms $\mathbf{D}_{\Sigma}^{s}\simeq\mathcal{L}_{A(\Omega)}(\mathbf{A}_{\Omega}^{s},A(\Sigma))$
and $\mathcal{D}_{\Sigma}\simeq\mathcal{L}_{A(\Omega)}(\mathcal{A}_{\Omega},A(\Sigma))$.}

\emph{Proof. }Set $\mathfrak{a}_{\Sigma}=\ker(A(\Omega)\to A(\Sigma))$,
so $A(\Sigma)\simeq A(\Omega)/\mathfrak{a}_{\Sigma}$. The definitions
immediately imply isomorphisms\begin{eqnarray*}
\mathbf{D}_{\Sigma}^{s} & \simeq & \mathcal{L}_{A(\Omega)}(\mathbf{A}_{\Omega}^{s},A(\Omega))/\mathfrak{a}_{\Sigma}\mathcal{L}_{A(\Omega)}(\mathbf{A}_{\Omega}^{s},A(\Omega))\\
 & \simeq & \mathcal{L}_{A(\Omega)}(\mathbf{A}_{\Omega}^{s},A(\Omega))/\mathcal{L}_{A(\Omega)}(\mathbf{A}_{\Omega}^{s},\mathfrak{a}_{\Sigma}),\end{eqnarray*}
so the first isomorphism will follow if we can verify that the sequence
\[
0\to\mathcal{L}_{A(\Omega)}(\mathbf{A}_{\Omega}^{s},\mathfrak{a}_{\Sigma})\to\mathcal{L}_{A(\Omega)}(\mathbf{A}_{\Omega}^{s},A(\Omega))\to\mathcal{L}_{A(\Omega)}(\mathbf{A}_{\Omega}^{s},A(\Sigma))\]
is exact on the right. Given a $k$-Banach space $E$, write $b(E)$
for the Banach space of bounded sequences $\{(e_{i})_{i\in\mathbb{N}},\mathrm{sup}_{i\in\mathbb{N}}\left|e_{i}\right|_{E}<\infty\}$.
Choosing an orthonormal basis of $\mathbf{A}(N(\mathbf{Z}_{p}),A(\Omega))^{s}$
gives rise to an isometry $\mathcal{L}_{A(\Omega)}(\mathbf{A}_{\Omega}^{s},E)\simeq b(E)$
for $E$ any Banach $A(\Omega)$-module. Thus we need to show the
surjectivity of the reduction map $b(A(\Omega))\to b(A(\Sigma))$.
Choose a presentation $A(\Omega)=T_{n}/\mathfrak{b}_{\Omega}$, so
$A(\Sigma)=T_{n}/\mathfrak{b}_{\Sigma}$ with $\mathfrak{b}_{\Omega}\subseteq\mathfrak{b}_{\Sigma}$.
Quite generally for any $\mathfrak{b}\subset T_{n}$, the function\[
f\in T_{n}/\mathfrak{b}\mapsto\left\Vert f\right\Vert _{\mathfrak{b}}=\mathrm{inf}_{\tilde{f}\in f+\mathfrak{b}}\left\Vert \tilde{f}\right\Vert _{T_{n}}\]
defines a norm on $T_{n}/\mathfrak{b}$. By Proposition 3.7.5/3 of \cite{BGR},
there is a unique Banach algebra structure on any affinoid algebra.
Hence for any sequence $(f_{i})_{i\in\mathbb{N}}\in b(A(\Sigma))$,
we may choose a bounded sequence of lifts $(\widetilde{f}_{i})_{i\in\mathbb{N}}\in b(T_{n})$;
reducing the latter sequence modulo $\mathfrak{b}_{\Omega}$, we are
done.

Taking inverse limits in the sequence we just proved to be exact,
the second isomorphism follows. $\square$

For any complete extension $k/\mathbf{Q}_{p}$, let $\mathscr{A}_{G}(k)$
denote the ring of $k$-valued algebraic functions on $G(\mathbf{Q}_{p})$;
we write $\mathscr{A}_{G}=\mathscr{A}_{G}(\mathbf{Q}_{p})$. $ $
We regard $\mathscr{A}_{G}(k)$ as a right $G(\mathbf{Q}_{p})$-module
via left translation of functions. $ $Suppose $\lambda\in X_{+}^{\ast}\subset\mathcal{W}(\mathbf{Q}_{p})$
is a dominant weight for $\mathbf{B}$, with $V_{\lambda}$ the corresponding
irreducible right $G(\mathbf{Q}_{p})$-representation of highest weight
$\lambda$. The function $v_{\lambda}$ defined on the big cell of
the Bruhat decomposition by\[
v_{\lambda}(n'tn)=\lambda(t),\,(n',t,n)\in N^{opp}(\mathbf{Q}_{p})\times T(\mathbf{Q}_{p})\times N(\mathbf{Q}_{p})\]
extends to a well-defined algebraic function $v_{\lambda}\in\mathscr{A}_{G}$.
By the Borel-Weil-Bott theorem, the $\mathbf{Q}_{p}[N(\mathbf{Q}_{p})]$-orbit
of $v_{\lambda}$ spans a canonical copy of $V_{\lambda}$ inside
$\mathscr{A}_{G}$ with $v_{\lambda}$ a highest weight vector.

More generally, suppose $\lambda\in\mathcal{W}(k)$ is an arithmetic
weight, with $\lambda=\lambda^{\mathrm{alg}}\cdot\epsilon$. Let $s=s[\lambda]$
be the smallest integer such that $T^{s}\subseteq\ker\epsilon$. Let
$\Lambda^{s}\subset\Lambda$ be the monoid generated by $\Lambda\cap\ker\sigma$
and $T^{s}$, and set $\Delta^{s}=I_{1}^{s}\Lambda^{s}I_{1}^{s}$.
We define $f_{\lambda}$ by \[
f_{\lambda}(n'tn)=\lambda^{\mathrm{alg}}(t)\epsilon(\sigma(t)),\,(n',t,n)\in N^{opp}(\mathbf{Q}_{p})\times T(\mathbf{Q}_{p})\times N(\mathbf{Q}_{p}).\]
Note that $f_{\lambda}|_{\Delta^{s[\lambda]}}=v_{\lambda^{\mathrm{alg}}}|_{\Delta^{s[\lambda]}}$.
The function $h\mapsto f_{\lambda}(hg),\, h\in I_{1}^{s[\lambda]}\,\mathrm{and}\, g\in\Delta^{s[\lambda]}$
defines an element of $\mathbf{A}_{\lambda}^{s}\otimes_{k}V_{\lambda^{\mathrm{alg}}}$
for any $s\geq s[\lambda]$.

\textbf{Proposition 2.1.3. }\emph{The formula \[
v\star t=\lambda^{\mathrm{alg}}(t^{-1}\sigma(t))v\cdot t,\, t\in\Lambda\,\mathrm{and}\, v\in V_{\lambda^{\mathrm{alg}}}\subset\mathcal{A}_{G}\]
extends to a well-defined right action of $\Delta$ on $V_{\lambda^{\mathrm{alg}}}$,
and the map}\begin{eqnarray*}
i_{\lambda}:\mathbf{D}_{\lambda}^{s} & \to & V_{\lambda^{\mathrm{alg}}}\\
\mu & \mapsto & i_{\lambda}(\mu)=\int f_{\lambda}(ng)\mu(n)\end{eqnarray*}
\emph{is equivariant for the $\star$-action of $\Delta^{s[\lambda]}$.}

For any $\delta\in\Lambda^{s[\lambda]}$, we calculate \begin{eqnarray*}
i_{\lambda}(\mu\star\delta) & = & \int f_{\lambda}(ng)(\mu\star\delta)(n)\\
 & = & \int f_{\lambda}(\sigma(\delta)\delta^{-1}n\delta g)\mu(n)\\
 & = & \lambda^{\mathrm{alg}}(\delta^{-1}\sigma(\delta))\epsilon(\sigma(\delta^{-1}\sigma(\delta)))\int f_{\lambda}(n\delta g)\mu(n)\\
 & = & \lambda^{\mathrm{alg}}(\delta^{-1}\sigma(\delta))\int f_{\lambda}(n\delta g)\mu(n)\\
 & = & i_{\lambda}(\mu)\star\delta\end{eqnarray*}
If $\gamma\in I_{1}^{s[\lambda]}$, the $\star$-actions are simply
the usual actions, and we have\begin{eqnarray*}
i_{\lambda}(\mu\star\gamma) & = & \int f_{\lambda}(hg)(\mu\star\gamma)(h)\\
 & = & \int f_{\lambda}(h\gamma g)\mu(h)\\
 & = & i_{\lambda}(\mu)\star\gamma.\end{eqnarray*}
$\square$ $ $

\subsubsection*{The case of $\mathrm{GL}_{n}/\mathbf{Q}_{p}$}

We examine the case of $\mathrm{GL}_{n}$ in detail. We choose $B$
and $B^{\mathrm{opp}}$ as the upper and lower triangular Borel subgroups,
respectively, and we identify $T$ with diagonal matrices. The splitting
$\sigma$ is canonically induced from the homomorphism\begin{eqnarray*}
\mathbf{Q}_{p}^{\times} & \to & \mathbf{Z}_{p}^{\times}\\
x & \mapsto & p^{-v_{p}(x)}x.\end{eqnarray*}
Since $T(\mathbf{Z}_{p})\simeq(\mathbf{Z}_{p}^{\times})^{n}$, we
canonically identify a character $\lambda:T(\mathbf{Z}_{p})\to R^{\times}$
with the $n$-tuple of characters $(\lambda_{1},\dots,\lambda_{n})$
where \begin{eqnarray*}
\lambda_{i}:\mathbf{Z}_{p}^{\times} & \to & R^{\times}\\
x & \mapsto & \lambda\circ\mathrm{diag}(\underset{i-1}{\underbrace{1,\dots,1}},x,1,\dots,1).\end{eqnarray*}
Dominant weights are identified with characters $\lambda=(\lambda_{1},\dots,\lambda_{n})$
with $\lambda_{i}(x)=x^{k_{i}}$ for integers $k_{1}\geq k_{2}\geq\dots\geq k_{n}$. 

We want to explain how to {}``twist away'' one dimension's worth
of weights in a canonical fashion. For any $\Omega\subset\mathcal{W}$,
a simple calculation shows that the $\star$-action of $\Delta$ on
$\mathbf{A}_{\Omega}^{s}$ is given explicitly by the formula\[
(\delta\star f)(x)=\chi_{\Omega}(\sigma(\mathrm{t}(\delta))\mathrm{t}(\delta)^{-1}\mathrm{t}(x\delta))f(\mathrm{n}(x\delta)),\,\delta\in\Delta,\, x\in N(\mathbf{Z}_{p}),\, f\in\mathbf{A}(N(\mathbf{Z}_{p}),A(\Omega))^{s}.\]
Given $1\leq i\leq n$, let $m_{i}(g)$ denote the determinant of
the upper-left $i$-by-$i$ block of $g\in\mathrm{GL}_{n}$. For any
$g\in\Delta$, a pleasant calculation left to the reader shows that\[
\mathrm{t}(g)=\mathrm{diag}(m_{1}(g),m_{1}(g)^{-1}m_{2}(g),\dots,m_{i}^{-1}(g)m_{i+1}(g),\dots,m_{n-1}(g)^{-1}\mathrm{det}g).\]
In particular, writing $\lambda^{0}=(\lambda_{1}\lambda_{n}^{-1},\lambda_{2}\lambda_{n}^{-1},\dots,\lambda_{n-1}\lambda_{n}^{-1},1)$
yields a canonical isomorphism\[
\mathbf{A}_{\lambda}^{s}\simeq\mathbf{A}_{\lambda^{0}}^{s}\otimes\lambda_{n}(\det\cdot|\det|_{p})\]
of $\Delta$-modules. If $\Gamma\subset\mathrm{GL}_{n}(\mathbf{Q})$
is a congruence lattice which satisfies $\det\Gamma=1$, this descends
to an isomorphism of Hecke modules\[
H_{\ast}(\Gamma,\mathbf{A}_{\lambda}^{s})\simeq H_{\ast}(\Gamma,\mathbf{A}_{\lambda^{0}}^{s})\otimes\lambda_{n}(\det\cdot|\det|_{p}).\]
As such we define the \emph{null space }$\mathcal{W}^{0}\subset\mathcal{W}$
as the subspace of weights which are trivial on the one-parameter
subgroup $\mathrm{diag}(1,\dots,1,x)$, with its induced rigid analytic
structure. Restricting our attention to weights in the null space
amounts to factoring out central twists by {}``wild'' characters.

In the case of $\mathrm{GL}_{2}$ we can be even more explicit. Here
$\Delta$ is generated by the center of $G(\mathbf{Q}_{p})$ and by
the monoid\[
\Sigma_{0}(p)=\left\{ g=\left(\begin{array}{cc}
a & b\\
c & d\end{array}\right)\in M_{2}(\mathbf{Z}_{p}),\, c\in p\mathbf{Z}_{p},\, a\in\mathbf{Z}_{p}^{\times},\, ad-bc\neq0\right\} .\]
Another simple calculation shows that the center of $G(\mathbf{Q}_{p})$
acts on $\mathbf{A}_{\lambda}^{s}$ through the character $z\mapsto\lambda(\sigma(z)),$
while the monoid $\Sigma_{0}(p)$ acts via\[
(g\cdot f)(x)=(\lambda_{1}\lambda_{2}^{-1})(a+cx)\lambda_{2}(\det g|\det g|_{p})f\left(\frac{b+dx}{a+cx}\right),\, f\in\mathbf{A}(N(\mathbf{Z}_{p}),k)^{s},\,\left(\begin{array}{cc}
1 & x\\
 & 1\end{array}\right)\in N(\mathbf{Z}_{p}),\]
exactly as in \cite{Strigidsymbs}.

\paragraph*{Remarks. }

There are some subtle differences between the different modules we
have defined. The assignment $\Omega\mapsto\mathcal{A}_{\Omega}$
describes a presheaf over $\mathcal{W}$, and the modules $\mathbf{A}_{\Omega}^{s}$
are orthonormalizable. On the other hand, the modules $\mathbf{D}_{\Omega}^{s}$
are \emph{not }obviously orthonormalizable, and $\mathcal{D}_{\Omega}$
doesn't obviously form a presheaf. There are alternate modules of
distributions available, namely $\tilde{\mathbf{D}}_{\Omega}^{s}=\mathcal{L}_{k}(\mathbf{A}(N(\mathbf{Z}_{p}),k)^{s},k)\widehat{\otimes}_{k}A(\Omega)$
and $\tilde{\mathcal{D}}_{\Omega}=\lim_{\leftarrow s}\tilde{\mathcal{D}}_{\Omega}^{s}$,
equipped with suitable actions such that there is a natural $A(\Omega)[\Delta]$-equivariant
embedding $\tilde{\mathcal{D}}_{\Omega}\hookrightarrow\mathcal{D}_{\Omega}$.
The modules $\tilde{\mathbf{D}}_{\Omega}^{s}$ are orthonormalizable,
and $\tilde{\mathcal{D}}_{\Omega}$ forms a presheaf over weight space,
but of course is not the continuous dual of $\mathcal{A}_{\Omega}$.
Despite these differences, the practical choice to work with one module
or the other is really a matter of taste: for any $\lambda\in\Omega(\overline{\mathbf{Q}_{p}})$,
there are isomorphisms \[
\mathcal{D}_{\Omega}\otimes_{A(\Omega)}A(\Omega)/\mathfrak{m}_{\lambda}\simeq\tilde{\mathcal{D}}_{\Omega}\otimes_{A(\Omega)}A(\Omega)/\mathfrak{m}_{\lambda}\simeq\mathcal{D}_{\lambda},\]
and in point of fact the slope-$\leq h$ subspaces of $H^{\ast}(K^{p},\mathcal{D}_{\Omega})$
and $H^{\ast}(K^{p},\tilde{\mathcal{D}}_{\Omega})$, when they are
defined, are canonically isomorphic as Hecke modules, as we show in
Proposition 2.4.6 below. One of our implicit goals in this paper is
to demonstrate the feasibility of working successfully with the modules
$\mathcal{D}_{\Omega}$ by treating the dual modules $\mathcal{A}_{\Omega}$
on an equal footing.

\subsection{Shimura manifolds and cohomology of local systems}

In this section we set up our conventions for the homology and cohomology
of local systems on Shimura manifolds. Following \cite{AS}, we compute
homology and cohomology using two different families of resolutions:
some extremely\emph{ }large {}``adelic'' resolutions which have
the advantage of making the Hecke action transparent, and resolutions
with good finiteness properties constructed from simplicial decompositions
of the Borel-Serre compactifications of locally symmetric spaces.

\subsubsection*{Resolutions and complexes}

Let $G/\mathbf{Q}$ be a connected reductive group. Let $G(\mathbf{R})^{\circ}$
denote the connected component of $G(\mathbf{R})$ containing the
identity element, with $G(\mathbf{Q})^{\circ}=G(\mathbf{Q})\cap G(\mathbf{R})^{\circ}$.
Fix a maximal compact subgroup $K_{\infty}\subset G(\mathbf{R})$
with $K_{\infty}^{\circ}$ the connected component containing the
identity, and let $Z_{\infty}$ denote the real points of a maximal
$\mathbf{Q}$-split torus contained in the center of $G$. Given an
open compact subgroup $K_{f}\subset G(\mathbf{A}_{f})$, we define
the \emph{Shimura manifold of level $K_{f}$ }by\begin{eqnarray*}
Y(K_{f}) & = & G(\mathbf{Q})\backslash G(\mathbf{A})/K_{f}K_{\infty}^{\circ}Z_{\infty}.\end{eqnarray*}
This is a possibly disconnected Riemannian orbifold. By strong approximation
there is a finite set of elements $\gamma(K_{f})=\{x_{i},x_{i}\in G(\mathbf{A}_{f})\}$
with\[
G(\mathbf{A})=\coprod_{x_{i}\in\gamma(K_{f})}G(\mathbf{Q})^{\circ}G(\mathbf{R})^{\circ}x_{i}K_{f}.\]
Defining $\Gamma(x_{i})=G(\mathbf{Q})^{\circ}\cap x_{i}K_{f}x_{i}^{-1}$,
we have a decomposition\[
Y(K_{f})=G(\mathbf{Q})\backslash G(\mathbf{A})/K_{f}K_{\infty}^{\circ}Z_{\infty}\simeq\coprod_{x_{i}\in\gamma(K_{f})}\Gamma(x_{i})\backslash D_{\infty},\]
where $D_{\infty}=G(\mathbf{R})^{\circ}/K_{\infty}^{\circ}Z_{\infty}$
is the symmetric space associated with $G$. We shall assume for the
remainder of this paper that $K_{f}$ is \emph{neat:} the image of
each $\Gamma(x_{i})$ in $G^{\mathrm{ad}}(\mathbf{R})$ is torsion-free,
in which case $Y(K_{f})$ is a smooth manifold. If $N$ is any right
$K_{f}$-module, the double quotient\[
\widetilde{N}=G(\mathbf{Q})\backslash\left(D_{\infty}\times G(\mathbf{A}_{f})\times N\right)/K_{f}\]
naturally gives rise to a local system on $Y(K_{f})$. Set $D_{\mathbf{A}}=D_{\infty}\times G(\mathbf{A}_{f})$,
and let $C_{\bullet}(D_{\mathbf{A}})$ denote the complex of singular
chains on $D_{\mathbf{A}}$ endowed with the natural action of $G(\mathbf{Q})\times G(\mathbf{A}_{f})$.
If $M$ and $N$ are left and right $K_{f}$-modules, respectively,
we define the complexes of \emph{adelic chains }and \emph{adelic cochains
}by\[
C_{\bullet}^{ad}(K_{f},M)=C_{\bullet}(D_{\mathbf{A}})\otimes_{\mathbf{Z}[G(\mathbf{Q})\times K_{f}]}M\]
and\[
C_{ad}^{\bullet}(K_{f},N)=\mathrm{Hom}_{\mathbf{Z}[G(\mathbf{Q})\times K_{f}]}(C_{\bullet}(D_{\mathbf{A}}),N).\]

\textbf{Proposition 2.2.1. }\emph{There is a canonical isomorphism\[
H^{\ast}(Y(K_{f}),\widetilde{N})=H^{\ast}(C_{ad}^{\bullet}(K_{f},N)).\]
}

\emph{Proof. }Let $C_{\bullet}(D_{\infty})(x_{i})$ denote the complex
of singular chains on $D_{\infty}$, endowed with the natural left
action of $\Gamma(x_{i})$ induced from the left action of $G(\mathbf{Q})$
on $D_{\infty}$; since $D_{\infty}$ is contractible, this is a free
resolution of $\mathbf{Z}$ in the category of $\mathbf{Z}[\Gamma(x_{i})]$-modules.
Let $N(x_{i})$ denote the left $\Gamma(x_{i})$-module whose underlying
module is $N$ but with the action $\gamma\cdot_{x_{i}}n=n|x_{i}^{-1}\gamma^{-1}x_{i}$.
Note that the local system $\widetilde{N}(x_{i})$ obtained by restricting
$\widetilde{N}$ to the connected component $\Gamma(x_{i})\backslash D_{\infty}$
of $Y(K_{f})$ is simply the quotient $\Gamma(x_{i})\backslash\left(D_{\infty}\times N(x_{i})\right)$.
Setting \[
C_{sing}^{\bullet}(K_{f},N)=\oplus_{i}\mathrm{Hom}_{\mathbf{Z}[\Gamma(x_{i})]}(C_{\bullet}(D_{\infty})(x_{i}),N(x_{i})),\]
the map $D_{\infty}\to(D_{\infty},x_{i})\subset D_{\mathbf{A}}$ induces
a morphism $x_{i}^{\ast}=\mathrm{Hom}(C_{\bullet}(D_{\mathbf{A}}),N)\to\mathrm{Hom}(C_{\bullet}(D_{\infty}),N)$,
which in turn induces an isomorphism\[
\oplus_{i}x_{i}^{\ast}:C_{ad}^{\bullet}(K_{f},N)\overset{\sim}{\to}\oplus_{i}\mathrm{Hom}_{\Gamma(x_{i})}(C_{\bullet}(D_{\infty})(x_{i}),N(x_{i})),\]
and passing to cohomology we have\begin{eqnarray*}
H^{\ast}(C_{ad}^{\bullet}(K_{f},N)) & \simeq & \oplus_{i}H^{\ast}(\Gamma(x_{i})\backslash D_{\infty},\widetilde{N}(x_{i}))\\
 & \simeq & H^{\ast}(Y(K_{f}),\widetilde{N})\end{eqnarray*}
as desired. $\square$

For each $x_{i}\in\gamma(K_{f})$, we choose a finite resolution $F_{\bullet}(x_{i})\to\mathbf{Z}\to0$
of $\mathbf{Z}$ by free left $\mathbf{Z}[\Gamma(x_{i})]$-modules
of finite rank as well as a homotopy equivalence $F_{\bullet}(x_{i})\overset{f_{i}}{\underset{g_{i}}{\rightleftarrows}}C_{\bullet}(D_{\infty})(x_{i})$.
We shall refer to the resolution $F_{\bullet}(x_{i})$ as a \emph{Borel-Serre
resolution}; the existence of such resolutions follows from taking
a finite simplicial decomposition of the Borel-Serre compactification
of $\Gamma(x_{i})\backslash D_{\infty}$. Setting \[
C_{\bullet}(K_{f},N)=\oplus_{i}F_{\bullet}(x_{i})\otimes_{\mathbf{Z}[\Gamma(x_{i})]}M(x_{i})\]
and\[
C^{\bullet}(K_{f},N)=\oplus_{i}\mathrm{Hom}_{\mathbf{Z}[\Gamma(x_{i})]}(F_{\bullet}(x_{i}),N(x_{i})),\]
the maps $f_{i},g_{i}$ induce homotopy equivalences\[
C_{\bullet}(K_{f},M)\overset{f_{\ast}}{\underset{g_{\ast}}{\rightleftarrows}}C_{\bullet}^{ad}(K_{f},M)\]
and\[
C^{\bullet}(K_{f},N)\overset{g^{\ast}}{\underset{f^{\ast}}{\rightleftarrows}}C_{ad}^{\bullet}(K_{f},M).\]
We refer to the complexes $C_{\bullet}(K_{f},-)$ and $C^{\bullet}(K_{f},-)$
as \emph{Borel-Serre complexes}, and we refer to these complexes together
with a \emph{fixed }set of homotopy equivalences $\{f_{i},g_{i}\}$
as \emph{augmented Borel-Serre complexes.}

\subsubsection*{Hecke operators}

Suppose now that $R$ is a commutative ring, $\Delta\subset G(\mathbf{A}_{f})$
is a monoid containing $K_{f}$, and $M$ is a left $R[\Delta]$-module.
The complex $C_{\bullet}(D_{\mathbf{A}})\otimes_{\mathbf{Z}[G(\mathbf{Q})]}M$
receives a $\Delta$-action via $\delta\cdot(\sigma\otimes m)=\sigma\delta^{-1}\otimes\delta m$,
and $C_{\bullet}^{ad}(K_{f},M)$ is naturally identified with the
$K_{f}$-coinvariants of this action. Given any double coset $K_{f}\delta K_{f}=\coprod_{j}K_{f}\delta_{j}$,
the action defined on pure tensors by the formula\[
[K_{f}\delta K_{f}]\cdot(\sigma\otimes m)=\sum_{j}\delta_{j}\cdot(\sigma\otimes m)\]
induces a well-defined algebra homomorphism\[
\xi:\mathcal{H}_{R}(K_{f}\backslash\Delta/K_{f})\to\mathrm{End}_{R}(C_{\bullet}^{ad}(K_{f},M)).\]
This action commutes with the boundary maps, and induces the usual
Hecke action defined by correspondences on cohomology. The map \begin{eqnarray*}
\tilde{\xi}:\mathcal{H}_{R}(K_{f}\backslash\Delta/K_{f}) & \to & \mathrm{End}_{\mathbf{K}(R)}(C_{\bullet}(K_{f},M))\\
T & \mapsto & g_{\ast}\circ\xi(T)\circ f_{\ast}\end{eqnarray*}
is well-defined, since $f_{\ast}$ and $g_{\ast}$ are well-defined
up to homotopy equivalence. For any Hecke operator $T$, we will abbreviate
$\tilde{\xi}(T)$ by $\tilde{T}$. Note that any individual lift $\tilde{T}$
is well-defined in $\mathrm{End}_{R}(C_{\bullet}(K_{f},M))$, but
if $T_{1}$ and $T_{2}$ commute in the abstract Hecke algebra, $\tilde{T}_{1}\tilde{T}_{2}$
and $\tilde{T}_{2}\tilde{T}_{1}$ will typically only commute up to
homotopy.

Likewise, if $N$ is a right $\Delta$-module, the complex $\mathrm{Hom}_{\mathbf{Z}[G(\mathbf{Q})]}(C_{\bullet}(D_{\mathbf{A}}),N)$
receives a natural $\Delta$-action via the formula $\phi|\delta=\phi(\sigma\delta^{-1})\delta$,
and $C_{ad}^{\bullet}(K_{f},N)$ is naturally the $K_{f}$-invariants
of this action. The formula \[
\phi|[K_{f}\delta K_{f}]=\sum_{j}\phi|\delta_{j}\]
yields an algebra homomorphism $\xi:\mathcal{H}_{R}(K_{f}\backslash\Delta/K_{f})\to\mathrm{End}_{R}(C_{ad}^{\bullet}(K_{f},N))$
which induces the usual Hecke action on cohomology, and $\tilde{\xi}=f^{\ast}\circ\xi\circ g^{\ast}$
defines an algebra homomorphism $\mathcal{H}_{R}(K_{f}\backslash\Delta/K_{f})\to\mathrm{End}_{\mathbf{K}(R)}(C^{\bullet}(K_{f},M))$.

It is extremely important for us that these Hecke actions are compatible
with the duality isomorphism\[
\mathrm{Hom}_{R}(C_{\bullet}(K_{f},M),P)\simeq C^{\bullet}(K_{f},\mathrm{Hom}_{R}(M,P)),\]
where $P$ is any $R$-module.

\subsection{Slope decompositions of modules and complexes}

Here we review the very general notion of slope decomposition introduced
in \cite{AS}. Let $A$ be a $k$-Banach algebra, and let $M$ be
an $A$-module equipped with an $A$-linear endomorphism $u:M\to M$
(for short, {}``an $A[u]$-module''). Fix a rational number $h\in\mathbf{Q}_{\geq0}$.
We say a polynomial $Q\in A[x]$ is \emph{multiplicative} if the leading
coefficient of $Q$ is a unit in $A$, and that $Q$ has \emph{slope
$\leq h$} if every edge of the Newton polygon of $Q$ has slope $\leq h$.
Write $Q^{\ast}(x)=x^{\deg Q}Q(1/x)$. An element $m\in M$ has slope
$\leq h$ if there is a multiplicative polynomial $Q\in A[T]$ of
slope $\leq h$ such that $ $$Q^{\ast}(u)\cdot m=0$. Let $M_{\leq h}$
be the set of elements of $M$ of slope $\leq h$; according to Proposition
4.6.2 of \emph{loc. cit.}, $M_{\leq h}$ is an $A$-submodule of $M$.

\textbf{Definition 2.3.1. }A \emph{slope-$\leq h$ decomposition }of
$M$ is an $A[u]$-module isomorphism\[
M\simeq M_{\leq h}\oplus M_{>h}\]
such that $M_{\leq h}$ is a finitely generated $A$-module and the
map $Q^{\ast}(u):M_{>h}\to M_{>h}$ is an $A$-module isomorphism
for every multiplicative polynomial $Q\in A[T]$ of slope $\leq h$.

The following proposition summarizes the fundamental results on slope
decompositions.

\textbf{Proposition 2.3.1 (Ash-Stevens): }
\begin{description}
\item [{\emph{a)}}] \emph{Suppose $M$ and $N$ are both $A[u]$-modules
with slope-$\leq h$ decompositions. If $\psi:M\to N$ is a morphism
of $A[u]$-modules, then $\psi(M_{\leq h})\subseteq N_{\leq h}$ and
$\psi(M_{>h})\subseteq N_{>h}$; in particular, a module can have
at most one slope-$\leq h$ decomposition. Furthermore, $\ker\psi$
and $\mathrm{im}\psi$ inherit slope-$\leq h$ decompositions. Given
a short exact sequence \[
0\to M\to N\to L\to0\]
of $A[u]$-modules, if two of the modules admit slope-$\leq h$ decompositions
then so does the third. }
\item [{\emph{b)}}] \emph{If $C^{\bullet}$ is a complex of $A[u]$-modules,
all with slope-$\leq h$ decompositions, then \[
H^{n}(C^{\bullet})\simeq H^{n}(C_{\leq h}^{\bullet})\oplus H^{n}(C_{>h}^{\bullet})\]
is a slope-$\leq h$ decomposition of $H^{n}(C^{\bullet})$.}
\end{description}
\emph{Proof. }This is a rephrasing of (a specific case of) Proposition
4.1.2 of \emph{loc. cit.} $\square$

Suppose now that $A$ is a reduced affinoid algebra, $M$ is an orthonormalizable
Banach $A$-module, and $u$ is a compact operator. Let \[
F(T)=\mathrm{det}(1-uT)|M\in A[[T]]\]
denote the Fredholm determinant for the $u$-action on $M$. We say
$F$ admits a \emph{slope-$\leq h$} \emph{factorization} if we can
write $F(T)=Q(T)\cdot R(T)$ where $Q$ is a multiplicative polynomial
of slope $\leq h$ and $R(T)\in A[[T]]$ is an entire power series
of slope $>h$. Theorem 3.3 of \cite{BuEigen} guarantees that $F$
admits a slope-$\leq h$ factorization if and only if $M$ admits
a slope-$\leq h$ decomposition. Furthermore, given a slope-$\leq h$
factorization $F(T)=Q(T)\cdot R(T)$, we obtain the slope-$\le h$
decomposition of $M$ upon setting $M_{\leq h}=\left\{ m\in M|Q^{\ast}(u)\cdot m=0\right\} $,
and $M_{\leq h}$ in this case is a finite flat $A$-module upon which
$u$ acts invertibly.%
\footnote{Writing $Q^{\ast}(x)=a+x\cdot r(x)$ with $r\in A[x]$ and $a\in A^{\times}$,
$u^{-1}$ on $M_{\leq h}$ is given explicitly by $-a^{-1}r(u)$.%
} Combining this with Theorem 4.5.1 of \cite{AS} and Proposition 2.3.1,
we deduce:

\textbf{Proposition 2.3.2. }\emph{If $C^{\bullet}$ is a bounded complex
of orthonormalizable Banach $A[u]$-modules, and $u$ acts compactly
on the total complex $\oplus_{i}C^{i}$, then for any $x\in\mathrm{Max}(A)$
and any $h\in\mathbf{Q}_{\geq0}$ there is an affinoid subdomain $\mathrm{Max}(A')\subset\mathrm{Max}(A)$
containing $x$ such that the complex $C^{\bullet}\widehat{\otimes}_{A}A'$
of $A'[u]$-modules admits a slope-$\leq h$ decomposition, and $(C^{\bullet}\widehat{\otimes}_{A}A')_{\leq h}$
is a complex of finite flat $A'$-modules.}

\textbf{Proposition 2.3.3. }\emph{If $M$ is an orthonormalizable
Banach $A$-module with a slope-$\leq h$ decomposition, and $A'$
is a Banach $A$-algebra, then $M\widehat{\otimes}_{A}A'$ admits
a slope-$\leq h$ decomposition and in fact\[
(M\widehat{\otimes}_{A}A')_{\leq h}\simeq M_{\leq h}\otimes_{A}A'.\]
}

\textbf{Proposition 2.3.4. }\emph{If $N\in\mathrm{Ban}_{A}$ is finite
and $M\in\mathrm{Ban}_{A}$ is an $A[u]$-module with a slope-$\leq h$
decomposition, the $A[u]$-modules $M\widehat{\otimes}_{A}N$ and
$\mathcal{L}_{A}(M,N)$ inherit slope-$\leq h$ decompositions.}

\emph{Proof. }This is an immediate consequence of the $A$-linearity
of the $u$-action and the fact that $-\widehat{\otimes}_{A}N$ and
$\mathcal{L}_{A}(-,N)$ commute with finite direct sums. $\square$

\subsection{Overconvergent (co)homology}

In this section we establish some foundational results on overconvergent
cohomology. These results likely follow from the formalism introduced
in Chapter 5 of \cite{AS}, but we give different proofs. We use the
notations introduced in §2.1-§2.3.

Fix a connected, reductive group $G/\mathbf{Q}$ with $G/\mathbf{Q}_{p}$
split, and fix a tame level $K^{p}\subset G(\mathbf{A}_{f}^{p})$;
in the presence of overconvergent coefficient modules, we \emph{always}
take our wild level subgroup of $G(\mathbf{Q}_{p})$ to be $I$, which
we drop from the notation, writing e.g. $H_{\ast}(K^{p},\mathbf{\mathcal{A}}_{\Omega})=H_{\ast}(K^{p}I,\mathcal{A}_{\Omega})$.
Fix an augmented Borel-Serre complex $C_{\bullet}(K^{p},-)=C_{\bullet}(K^{p}I,-)$.
Fix an element $t\in\Delta^{+}$, and let $\tilde{U}=\tilde{U}_{t}$
denote the lifting of $U_{t}=[ItI]$ to an endomorphism of the complex
$C_{\bullet}(K^{p},-)$ defined in §2.2. Given a connected admissible
open affinoid subset $\Omega\subset\mathcal{W}$ and any integer $s\geq s[\Omega]$,
the endomorphism $\tilde{U}_{t}\in\mathrm{End}_{A(\Omega)}(C_{\bullet}(K^{p},\mathbf{A}_{\Omega}^{s}))$
is completely continuous; let \[
F_{\Omega}^{s}(T)=\det(1-T\tilde{U}_{t})|C_{\bullet}(K^{p},\mathbf{A}_{\Omega}^{s})\in A(\Omega)[[T]]\]
denote its Fredholm determinant. We say $(U_{t},\Omega,h)$ is a \emph{slope
datum} if $C_{\bullet}(K^{p},\mathbf{A}_{\Omega}^{s})$ admits a slope-$\leq h$
decomposition for the $\tilde{U}_{t}$ action for some $s\geq s[\Omega]$;
we shall see shortly that this agrees with the definition given in
the introduction.

\textbf{Proposition 2.4.1. }\emph{The function $F_{\Omega}^{s}(T)$
is independent of $s$.}

\emph{Proof. }For any integer $s\geq s[\Omega]$ we write $C_{\bullet}^{s}=C_{\bullet}(K^{p},\mathbf{A}_{\Omega}^{s})$
for brevity. By construction, the operator $\tilde{U}_{t}$ factors
into compositions $\rho_{s}\circ\check{U}_{t}$ and $\check{U}_{t}\circ\rho_{s+1}$
where $\check{U}_{t}:C_{\bullet}^{s}\to C_{\bullet}^{s-1}$ is continuous
and $\rho_{s}:C_{\bullet}^{s-1}\to C_{\bullet}^{s}$ is completely
continuous. Now, considering the commutative diagram%
\footnote{The importance of drawing diagrams like this seems to have first been
realized by Hida.%
}\[
\xymatrix{C_{\bullet}^{s}\ar[r]^{\check{U}_{t}}\ar[d]^{\tilde{U}_{t}} & C_{\bullet}^{s-1}\ar[dl]_{\rho_{s}}\ar[d]^{\tilde{U}_{t}}\\
C_{\bullet}^{s}\ar[r]_{\check{U}_{t}} & C_{\bullet}^{s-1}}
\]
we calculate\begin{eqnarray*}
\det(1-T\tilde{U}_{t})|C_{\bullet}^{s} & = & \det(1-T\rho_{s}\circ\check{U}_{t})|C_{\bullet}^{s}\\
 & = & \det(1-T\check{U}_{t}\circ\rho_{s})|C_{\bullet}^{s-1}\\
 & = & \det(1-T\tilde{U}_{t})|C_{\bullet}^{s-1},\end{eqnarray*}
where the second line follows from Lemma 2.7 of \cite{BuEigen}, so
$F_{\Omega}^{s}(T)=F_{\Omega}^{s-1}(T)$ for all $s>s[\Omega]$. $\square$

\textbf{Proposition 2.4.2. }\emph{The slope-$\leq h$ subcomplex $C_{\bullet}(K^{p},\mathbf{A}_{\Omega}^{s})_{\leq h}$,
if it exists, is independent of $s$. If $\Omega'$ is an affinoid
subdomain of $\Omega$, then there is a canonical isomorphism\[
C_{\bullet}(K^{p},\mathbf{A}_{\Omega}^{s})_{\leq h}\otimes_{A(\Omega)}A(\Omega')\simeq C_{\bullet}(K^{p},\mathbf{A}_{\Omega'}^{s})_{\leq h}\]
for any $s\geq s[\Omega]$.}

\emph{Proof. }Since $F_{\Omega}^{s}(T)$ is independent of $s$, we
simply write $F_{\Omega}(T)$. Suppose we are given a slope-$\leq h$
factorization $F_{\Omega}(T)=Q(T)\cdot R(T)$; by the remarks in §2.3,
setting $C_{\bullet}(K^{p},\mathbf{A}_{\Omega}^{s})_{\leq h}=\ker Q^{\ast}(\tilde{U}_{t})$
yields a slope-$\leq h$ decomposition of $C_{\bullet}(K^{p},\mathbf{A}_{\Omega}^{s})$
for any $s\geq s[\Omega]$. By Proposition 2.3.1, the injection $\rho_{s}:C_{\bullet}(K^{p},\mathbf{A}_{\Omega}^{s-1})\hookrightarrow C_{\bullet}(K^{p},\mathbf{A}_{\Omega}^{s})$
gives rise to a canonical injection \[
\rho_{s}:C_{\bullet}(K^{p},\mathbf{A}_{\Omega}^{s-1})_{\leq h}\hookrightarrow C_{\bullet}(K^{p},\mathbf{A}_{\Omega}^{s})_{\leq h}\]
for any $s>s[\Omega]$. The operator $\tilde{U}_{t}$ acts invertibly
on $C_{\bullet}(K^{p},\mathbf{A}_{\Omega}^{s})_{\leq h}$, and its
image factors through $\rho_{s}$, so $\rho_{s}$ is surjective and
hence bijective. This proves the first claim.

For the second claim, by Proposition 2.3.3 we have\begin{eqnarray*}
C_{\bullet}(K^{p},\mathbf{A}_{\Omega}^{s})_{\leq h}\otimes_{A(\Omega)}A(\Omega') & \simeq & \left(C_{\bullet}(K^{p},\mathbf{A}_{\Omega}^{s})\widehat{\otimes}_{A(\Omega)}A(\Omega')\right)_{\leq h}\\
 & \simeq & C_{\bullet}(K^{p},\mathbf{A}_{\Omega'}^{s})_{\leq h},\end{eqnarray*}
so the result now follows from the first claim. $\square$

\textbf{Proposition 2.4.3. }\emph{Given a slope datum $(U_{t},\Omega,h)$
and an affinoid subdomain $\Omega'\subset\Omega$, there is a canonical
isomorphism\[
H_{\ast}(K^{p},\mathbf{A}_{\Omega}^{s})_{\leq h}\otimes_{A(\Omega)}A(\Omega')\simeq H_{\ast}(K^{p},\mathbf{A}_{\Omega'}^{s})_{\leq h}\]
for any $s\geq s[\Omega]$.}

\emph{Proof. }Since $A(\Omega')$ is $A(\Omega)$-flat, the functor
$-\otimes_{A(\Omega)}A(\Omega')$ commutes with taking the homology
of any complex of $A(\Omega)$-modules. Thus we calculate\begin{eqnarray*}
H_{\ast}(K^{p},\mathbf{A}_{\Omega}^{s})_{\leq h}\otimes_{A(\Omega)}A(\Omega) & \simeq & H_{\ast}\left(C_{\bullet}(K^{p},\mathbf{A}_{\Omega}^{s})_{\leq h}\right)\otimes_{A(\Omega)}A(\Omega')\\
 & \simeq & H_{\ast}\left(C_{\bullet}(K^{p},\mathbf{A}_{\Omega}^{s})_{\leq h}\otimes_{A(\Omega)}A(\Omega')\right)\\
 & \simeq & H_{\ast}\left(C_{\bullet}(K^{p},\mathbf{A}_{\Omega'}^{s})_{\leq h}\right)\\
 & \simeq & H_{\ast}(K^{p},\mathbf{A}_{\Omega'}^{s})_{\leq h},\end{eqnarray*}
where the third line follows from Proposition 2.3.3. $\square$

\textbf{Proposition 2.4.4. }\emph{Given a slope datum $(U_{t},\Omega,h)$,
the complex $C_{\bullet}(K^{p},\mathcal{A}_{\Omega})$ and the homology
module $H_{\ast}(K^{p},\mathcal{A}_{\Omega})$ admit slope-$\leq h$
decompositions, and there is an isomorphism \[
H_{\ast}(K^{p},\mathcal{A}_{\Omega})_{\leq h}\simeq H_{\ast}(K^{p},\mathbf{A}_{\Omega}^{s})_{\leq h}\]
for any $s\geq s[\Omega]$. Furthermore, given an affinoid subdomain
$\Omega'\subset\Omega$, there is a canonical isomorphism\[
H_{\ast}(K^{p},\mathcal{A}_{\Omega})_{\leq h}\otimes_{A(\Omega)}A(\Omega')\simeq H_{\ast}(K^{p},\mathcal{A}_{\Omega'})_{\leq h}.\]
}

\emph{Proof. }For any fixed $s\geq s[\Omega]$, we calculate\begin{eqnarray*}
C_{\bullet}(K^{p},\mathcal{A}_{\Omega}) & \simeq & \lim_{\substack{\to\\
s'}
}C_{\bullet}(K^{p},\mathbf{A}_{\Omega}^{s'})\\
 & \simeq & \lim_{\substack{\to\\
s'}
}C_{\bullet}(K^{p},\mathbf{A}_{\Omega}^{s'})_{\leq h}\oplus C_{\bullet}(K^{p},\mathbf{A}_{\Omega}^{s'})_{>h}\\
 & \simeq & C_{\bullet}(K^{p},\mathbf{A}_{\Omega}^{s})_{\leq h}\oplus\lim_{\substack{\to\\
s'}
}C_{\bullet}(K^{p},\mathbf{A}_{\Omega}^{s'})_{>h}\end{eqnarray*}
with the third line following from Proposition 2.4.2. The two summands
in the third line naturally form the components of a slope-$\leq h$
decomposition, so passing to homology yields the first sentence of
the proposition, and the second sentence then follows immediately
from Proposition 2.4.3. $\square$

We're now in a position to prove the subtler cohomology analogue of
Proposition 2.4.4.

\textbf{Proposition 2.4.5. }\emph{Given a slope datum $(U_{t},\Omega,h)$
and a rigid Zariski-closed subset $\Sigma\subset\Omega$, the complex
$C^{\bullet}(K^{p},\mathcal{D}_{\Sigma})$ and the cohomology module
$H^{\ast}(K^{p},\mathcal{D}_{\Sigma})$ admit slope-$\leq h$ decompositions,
and there is an isomorphism \[
H^{\ast}(K^{p},\mathcal{D}_{\Sigma})_{\leq h}\simeq H^{\ast}(K^{p},\mathbf{D}_{\Sigma}^{s})_{\leq h}\]
for any $s\geq s[\Omega]$. Furthermore, given an affinoid subdomain
$\Omega'\subset\Omega$, there are canonical isomorphisms\[
C^{\bullet}(K^{p},\mathcal{D}_{\Omega})_{\leq h}\otimes_{A(\Omega)}A(\Omega')\simeq C^{\bullet}(K^{p},\mathcal{D}_{\Omega'})_{\leq h}\]
and\[
H^{\ast}(K^{p},\mathcal{D}_{\Omega})_{\leq h}\otimes_{A(\Omega)}A(\Omega')\simeq H^{\ast}(K^{p},\mathcal{D}_{\Omega'})_{\leq h}.\]
}

\emph{Proof. }By a topological version of the duality stated in the
final paragraph of §2.2, we have a natural isomorphism\begin{eqnarray*}
C^{\bullet}(K^{p},\mathbf{D}_{\Sigma}^{s}) & = & C^{\bullet}(K^{p},\mathcal{L}_{A(\Omega)}(\mathbf{A}_{\Omega}^{s},A(\Sigma)))\\
 & \simeq & \mathcal{L}_{A(\Omega)}(C_{\bullet}(K^{p},\mathbf{A}_{\Omega}^{s}),A(\Sigma))\end{eqnarray*}
for any $s\geq s[\Omega]$. By assumption, $C_{\bullet}(K^{p},\mathbf{A}_{\Omega}^{s})$
admits a slope-$\leq h$ decomposition, so we calculate\begin{eqnarray*}
C^{\bullet}(K^{p},\mathbf{D}_{\Sigma}^{s}) & \simeq & \mathcal{L}_{A(\Omega)}(C_{\bullet}(K^{p},\mathbf{A}_{\Omega}^{s}),A(\Sigma))\\
 & \simeq & \mathcal{L}_{A(\Omega)}(C_{\bullet}(K^{p},\mathbf{A}_{\Omega}^{s})_{\leq h},A(\Sigma))\\
 &  & \oplus\mathcal{L}_{A(\Omega)}(C_{\bullet}(K^{p},\mathbf{A}_{\Omega}^{s})_{>h},A(\Sigma)).\end{eqnarray*}
By Proposition 2.4.2, passing to the inverse limit over $s$ in this
isomorphism yields a slope-$\leq h$ decomposition of $C^{\bullet}(K^{p},\mathcal{D}_{\Sigma})$
together with a natural isomorphism \[
C^{\bullet}(K^{p},\mathcal{D}_{\Sigma})_{\leq h}\simeq C^{\bullet}(K^{p},\mathbf{D}_{\Sigma}^{s})_{\leq h}\simeq\mathcal{L}_{A(\Omega)}(C^{\bullet}(K^{p},\mathbf{A}_{\Omega}^{s})_{\leq h},A(\Sigma))\]
for any $s\geq s[\Omega]$. This proves the first sentence of the
proposition.

For the second sentence, we first note that since $C_{\bullet}(K^{p},\mathbf{A}_{\Omega}^{s})_{\leq h}$
is a complex of finite $A(\Omega)$-modules, the natural map\[
\mathcal{L}_{A(\Omega)}(C_{\bullet}(K^{p},\mathbf{A}_{\Omega}^{s})_{\leq h},A(\Omega))\to\mathrm{Hom}_{A(\Omega)}(C_{\bullet}(K^{p},\mathbf{A}_{\Omega}^{s})_{\leq h},A(\Omega))\]
is an isomorphism by Lemma 2.2 of \cite{BuEigen}. Next, note that
if $R$ is a commutative ring, $S$ is a flat $R$-algebra, and $M,N$
are $R$-modules with $M$ finitely presented, the natural map $\mathrm{Hom}_{R}(M,N)\otimes_{R}S\to\mathrm{Hom}_{S}(M\otimes_{R}S,N\otimes_{R}S)$
is an isomorphism. With these two facts in hand, we calculate as follows:
\begin{eqnarray*}
C^{\bullet}(K^{p},\mathcal{D}_{\Omega})_{\leq h}\otimes_{A(\Omega)}A(\Omega') & \simeq & \mathrm{Hom}_{A(\Omega)}(C_{\bullet}(K^{p},\mathbf{A}_{\Omega}^{s})_{\leq h},A(\Omega))\otimes_{A(\Omega)}A(\Omega')\\
 & \simeq & \mathrm{Hom}_{A(\Omega')}(C_{\bullet}(K^{p},\mathbf{A}_{\Omega}^{s})_{\leq h}\otimes_{A(\Omega)}A(\Omega'),A(\Omega'))\\
 & \simeq & \mathrm{Hom}_{A(\Omega')}(C_{\bullet}(K^{p},\mathbf{A}_{\Omega'}^{s})_{\leq h},A(\Omega'))\\
 & \simeq & C^{\bullet}(K^{p},\mathcal{D}_{\Omega'})_{\leq h},\end{eqnarray*}
where the third line follows from Proposition 2.3.3. Passing to cohomology,
the result follows as in the proof of Proposition 2.4.3.

Recall from the end of §2.1 the alternate module of distributions
$\tilde{\mathcal{D}}_{\Omega}$. For completeness's sake, we sketch
the following result.

\textbf{Proposition 2.4.6. }\emph{If the complexes $C^{\bullet}(K^{p},\tilde{\mathcal{D}}_{\Omega})$
and $C^{\bullet}(K^{p},\mathcal{D}_{\Omega})$ both admit slope-$\leq h$
decompositions for the $\tilde{U}_{t}$-action, the natural map $\tilde{\mathcal{D}}_{\Omega}\to\mathcal{D}_{\Omega}$
induces an isomorphism\[
H^{\ast}(K^{p},\tilde{\mathcal{D}}_{\Omega})_{\leq h}\simeq H^{\ast}(K^{p},\mathcal{D}_{\Omega})_{\leq h}.\]
}

\emph{Proof. }The morphism $f:\tilde{\mathcal{D}}_{\Omega}\to\mathcal{D}_{\Omega}$
induces a morphism $\phi:C^{\bullet}(K^{p},\tilde{\mathcal{D}}_{\Omega})_{\leq h}\to C^{\bullet}(K^{p},\mathcal{D}_{\Omega})_{\leq h}$.
Choose an arbitrary weight $\lambda\in\Omega(\overline{\mathbf{Q}_{p}})$;
since\[
\mathcal{D}_{\Omega}\otimes_{A(\Omega)}A(\Omega)/\mathfrak{m}_{\lambda}\simeq\tilde{\mathcal{D}}_{\Omega}\otimes_{A(\Omega)}A(\Omega)/\mathfrak{m}_{\lambda}\simeq\mathcal{D}_{\lambda},\]
we calculate\begin{eqnarray*}
C^{\bullet}(K^{p},\tilde{\mathcal{D}}_{\Omega})_{\leq h}\otimes_{A(\Omega)}A(\Omega)/\mathfrak{m}_{\lambda} & \simeq & C^{\bullet}(K^{p},\tilde{\mathcal{D}}_{\Omega}\otimes_{A(\Omega)}A(\Omega)/\mathfrak{m}_{\lambda})_{\leq h}\\
 & \simeq & C^{\bullet}(K^{p},\mathcal{D}_{\lambda})_{\leq h}\\
 & \simeq & C^{\bullet}(K^{p},\mathcal{D}_{\Omega}\otimes_{A(\Omega)}A(\Omega)/\mathfrak{m}_{\lambda})_{\leq h}\\
 & \simeq & C^{\bullet}(K^{p},\mathcal{D}_{\Omega})_{\leq h}\otimes_{A(\Omega)}A(\Omega)/\mathfrak{m}_{\lambda},\end{eqnarray*}
where the first and fourth lines follow by Proposition 2.3.3. Thus
$\phi\,\mathrm{mod}\,\mathfrak{m}_{\lambda}$ is an isomorphism\emph{
}in the category of $A(\Omega)/\mathfrak{m}_{\lambda}$-module complexes.
Writing $C^{\bullet}(\phi)$ for the cone of $\phi$, we have an exact
triangle\[
C^{\bullet}(K^{p},\tilde{\mathcal{D}}_{\Omega})_{\leq h}\overset{\phi}{\to}C^{\bullet}(K^{p},\mathcal{D}_{\Omega})_{\leq h}\to C^{\bullet}(\phi)\to C^{\bullet}(K^{p},\tilde{\mathcal{D}}_{\Omega})_{\leq h}[1]\]
in $\mathbf{D}^{b}(A(\Omega))$. Since $C^{\bullet}(K^{p},\tilde{\mathcal{D}}_{\Omega})_{\leq h}$
and $C^{\bullet}(K^{p},\mathcal{D}_{\Omega})_{\leq h}$ are complexes
of projective $A(\Omega)$-modules, the cone $C^{\bullet}(\phi)$
is a complex of projective modules as well; therefore, the functors
$-\otimes_{A(\Omega)}^{\mathbf{L}}A(\Omega)/\mathfrak{m}$ and $-\otimes_{A(\Omega)}A(\Omega)/\mathfrak{m}$
agree on these three complexes (cf. \cite{Weibelhomalg}, 10.6). In
particular, applying $-\otimes_{A(\Omega)}A(\Omega)/\mathfrak{m}_{\lambda}$
yields an exact triangle\[
C^{\bullet}(K^{p},\mathcal{D}_{\lambda})_{\leq h}\overset{\phi\,\mathrm{mod}\,\mathfrak{m}_{\lambda}}{\longrightarrow}C^{\bullet}(K^{p},\mathcal{D}_{\lambda})_{\leq h}\to C^{\bullet}(\phi)\otimes_{A(\Omega)}A(\Omega)/\mathfrak{m}_{\lambda}\to C^{\bullet}(K^{p},\mathcal{D}_{\lambda})_{\leq h}[1]\]
in $\mathbf{D}^{b}(A(\Omega)/\mathfrak{m}_{\lambda})$. Since $\phi\,\mathrm{mod}\,\mathfrak{m}_{\lambda}$
is an isomorphism\emph{,} $C^{\bullet}(\phi)\otimes_{A(\Omega)}A(\Omega)/\mathfrak{m}_{\lambda}$
is acyclic, which in turn implies that $\mathfrak{m}_{\lambda}$ is
not contained in the support of $H^{\ast}(C^{\bullet}(\phi))$. But
$\mathfrak{m}_{\lambda}$ is an arbitrary maximal ideal in $A(\Omega)$,
so Nakayama's lemma now shows that $H^{\ast}(C^{\bullet}(\phi))$
vanishes identically. Thus $C^{\bullet}(\phi)$ is acyclic and $\phi$
is a quasi-isomorphism as desired. $\square$

We now recall a fundamental theorem of Ash-Stevens and Urban (Theorem
6.4.1 of \cite{AS}, Proposition 4.3.10 of \cite{UrEigen}) relating
overconvergent cohomology classes of small slope to classical automorphic
forms. The possibility of such a result was largely the original \emph{raison
d'etre }of overconvergent cohomology. For any weight $\lambda\in\mathcal{W}(\overline{\mathbf{Q}_{p}})$
we define the \emph{finite slope subspace} $H^{\ast}(K^{p},\mathcal{D}_{\lambda})_{\mathrm{fs}}$
as the intersection $\cap_{n\geq1}\mathrm{im}U_{t}^{n}|H^{\ast}(K^{p},\mathcal{D}_{\lambda})$,
where $U_{t}$ is any controlling operator. This definition is independent
of the choice of $t$. Now suppose $\lambda$ is an arithmetic weight,
and let $k$ be the finite Galois extension of $\mathbf{Q}_{p}$ generated
by the values of $\lambda$. We write $\lambda=\lambda^{\mathrm{alg}}\varepsilon$
where $\lambda^{\mathrm{alg}}\in X_{+}^{\ast}$ and $\varepsilon:T(\mathbf{Z}_{p})\to k^{\times}$
is a finite order character. Given $w\in W$, write $^{w}\lambda^{\mathrm{alg}}$
for the character $w\cdot(\lambda^{\mathrm{alg}}+\rho)-\rho$. Let
$\mathbf{T}_{\lambda,h}$ be the subalgebra of $\mathrm{End}_{k}(H^{\ast}(K^{p},\mathcal{D}_{\lambda})_{\leq h})$
generated by $\mathbf{T}(K^{p})$, and let $\mathbf{T}_{\lambda}$
be the projective limit $\lim_{\infty\leftarrow h}\mathbf{T}_{\lambda,h}$.
If $\phi:\mathbf{T}_{\lambda}\to\overline{\mathbf{Q}_{p}}$ is an
eigenpacket, the semigroup character \begin{eqnarray*}
\mu_{\phi}:\Lambda^{+} & \to & \mathbf{Q}_{>0}\\
t & \mapsto & v_{p}(\phi(U_{t}))\end{eqnarray*}
extends uniquely to a character $\mu_{\phi}:T(\mathbf{Q}_{p})/T(\mathbf{Z}_{p})\to\mathbf{Q}$.

\textbf{Definition 2.4.7. }\emph{Fix a controlling operator $U_{t}$,
$t\in\Lambda$. Given an arithmetic weight $\lambda=\lambda^{\mathrm{alg}}\varepsilon$,
a rational number $h$ is a }small slope for $\lambda$ \emph{if \[
h<\inf_{w\in W\smallsetminus\{\mathrm{id}\}}v_{p}(^{w}\lambda^{\mathrm{alg}}(t))-v_{p}(\lambda^{\mathrm{alg}}(t)).\]
If $\phi:\mathbf{T}_{\lambda}\to\overline{\mathbf{Q}_{p}}$ is an
eigenpacket, $\phi$ is }numerically non-critical \emph{if \[
v_{p}(\mu_{\phi}(t))<\mathrm{inf}_{w\in W\smallsetminus\{id\}}v_{p}(^{w}\lambda^{\mathrm{alg}}(t))-v_{p}(\lambda^{\mathrm{alg}}(t))\]
for all $t=\mu^{\vee}(p^{-1})$ with $\mu^{\vee}\in X_{\ast}(T)$
a positive coroot.}

\textbf{Theorem 2.4.8 (Ash-Stevens, Urban). }\emph{Fix an arithmetic
weight $\lambda=\lambda^{\mathrm{alg}}\varepsilon$.}

\emph{a) Fix a controlling operator $U_{t}$. If $h$ is a small slope
for $\lambda$, there is a natural isomorphism of Hecke modules\[
H^{\ast}(K^{p},\mathcal{D}_{\lambda})_{\leq h}\simeq H^{\ast}(Y(K^{p}I_{1}^{c}),V_{\lambda^{\mathrm{alg}}}(k))_{\leq h}^{T(\mathbf{Z}/p^{c}\mathbf{Z})=\epsilon}\]
for any $c\geq c(\varepsilon)$.}

\emph{b) If $\phi\in H^{\ast}(K^{p},\mathcal{D}_{\lambda})_{\mathrm{fs}}$
is a numerically non-critical eigenclass, the morphism\[
H^{\ast}(K^{p},\mathcal{D}_{\lambda})\to H^{\ast}(Y(K^{p}I_{1}^{c}),V_{\lambda^{\mathrm{alg}}}(k))\]
is a Hecke-equivariant isomorphism on the $\phi$-generalized eigenspace
for any $c\geq c(\varepsilon)$.}

More generally, suppose $\phi:\mathbf{T}(K^{p})\to\overline{\mathbf{Q}_{p}}$
is an eigenpacket associated with a classical algebraic automorphic
representation of weight $\lambda^{\mathrm{alg}}$. We say $\phi$
has \emph{finite slope }if $\phi(U_{t})\neq0$ for some controlling
operator. Choosing $h\geq v_{p}(\phi(U_{t}))$, we define $\phi$
to be \emph{non-critical }if the map\[
H^{\ast}(K^{p},\mathcal{D}_{\lambda})_{\leq h}\to H^{\ast}(Y(K^{p}I_{1}^{c}),V_{\lambda^{\mathrm{alg}}}(k))_{\leq h}\]
induces a Hecke-equivariant isomorphism on the $\phi$-generalized
eigenspaces. According to Theorem 2.4.8, every numerically non-critical
eigenpacket is non-critical. However, all known examples point to
the suspicion that numerically critical cuspidal eigenpackets are
typically non-critical, with the critical cuspidal eigenpackets rising
via functorial maps from smaller groups.

\section{Proofs of the main results}

Fix a choice of a tame level $K^{p}$ and an augmented Borel-Serre
complex $C_{\bullet}(K^{p},-)$.

\subsection{Proof of Theorem 1.1}

\paragraph*{The spectral sequences}

Fix a rigid Zariski closed subset $\Sigma\subset\Omega$. By the isomorphisms
proved in §2.4, it suffices to construct a spectral sequence $\mathrm{Ext}_{A(\Omega)}^{i}(H_{j}(K^{p},\mathbf{A}_{\Omega}^{s})_{\leq h},A(\Sigma))\Rightarrow H^{i+j}(K^{p},\mathbf{D}_{\Sigma}^{s})_{\leq h}$
for some $s\geq s[\Omega]$. From Proposition 2.4.5, we have a natural
isomorphism\[
\mathrm{Hom}_{A(\Omega)}(C_{\bullet}(K^{p},\mathbf{A}_{\Omega}^{s})_{\leq h},A(\Sigma))\simeq C^{\bullet}(K^{p},\mathbf{D}_{\Sigma}^{s})_{\leq h}.\]
Let $\iota:C_{\bullet}(K^{p},\mathbf{A}_{\Omega}^{s})_{\leq h}\hookrightarrow C_{\bullet}(K^{p},\mathbf{A}_{\Omega}^{s})$
denote the canonical inclusion, and $\pi:C_{\bullet}(K^{p},\mathbf{A}_{\Omega}^{s})\twoheadrightarrow C_{\bullet}(K^{p},\mathbf{A}_{\Omega}^{s})_{\leq h}$
the canonical projection. The formula $T\to\pi\circ\tilde{\xi}(T)\circ\iota$
defines an algebra homomorphism $\mathbf{T}(K^{p})\to\mathrm{End}_{\mathbf{K}(A(\Omega))}(C_{\bullet}(K^{p},\mathbf{A}_{\Omega}^{s})_{\leq h})$
which induces the usual Hecke action on $H_{\ast}(K^{p},\mathbf{A}_{\Omega}^{s})_{\leq h}$;
making the analogous definition on $C^{\bullet}(K^{p},\mathbf{D}_{\Sigma}^{s})_{\leq h}$,
the above isomorphism of complexes upgrades to an isomorphism of $\mathbf{T}(K^{p})$-module
complexes in $\mathbf{K}^{b}(A(\Omega))$. Since $C_{\bullet}(K^{p},\mathbf{A}_{\Omega}^{s})_{\leq h}$
is a complex of projective $A(\Omega)$-modules, this isomorphism
in turn induces an isomorphism\[
\mathbf{R}\mathrm{Hom}_{A(\Omega)}(C_{\bullet}(K^{p},\mathbf{A}_{\Omega}^{s})_{\leq h},A(\Sigma))\simeq C^{\bullet}(K^{p},\mathbf{D}_{\Sigma}^{s})_{\leq h}\]
of $\mathbf{T}(K^{p})$-module complexes in $\mathbf{D}^{b}(A(\Omega))$
(cf.\emph{ }Theorem 10.7.4 of \cite{Weibelhomalg}). Passing to cohomology
yields \[
\mathbf{E}\mathrm{xt}_{A(\Omega)}^{i}(C_{\bullet}(K^{p},\mathbf{A}_{\Omega}^{s})_{\leq h},A(\Sigma))\simeq H^{i}(K^{p},\mathbf{D}_{\Sigma}^{s})_{\leq h}.\]
Quite generally, if $C_{\bullet}$ is a chain complex of $R$-modules
equipped with an algebra homomorphism $\varphi:S\to\mathrm{End}_{\mathbf{K}(R)}(C_{\bullet})$
and $N$ is any $R$-module, $\varphi$ induces a natural $S$-module
structure on the homology groups $H_{\ast}(C_{\bullet})$ and the
hyperext groups $\mathbf{E}\mathrm{xt}_{R}^{\ast}(C_{\bullet},N)$,
and the hyperext spectral sequence\[
E_{2}^{i,j}=\mathrm{Ext}_{R}^{i}(H_{j}(A_{\bullet}),N)\Rightarrow\mathbf{E}\mathrm{xt}_{R}^{i+j}(A_{\bullet},N)\]
is a spectral sequence of $S$-modules. The result follows.

For the Tor spectral sequence, the isomorphism\[
C^{\bullet}(K^{p},\mathbf{D}_{\Omega}^{s})_{\leq h}\otimes_{A(\Omega)}A(\Sigma)\simeq C^{\bullet}(K^{p},\mathbf{D}_{\Sigma}^{s})_{\leq h}\]
yields an isomorphism\[
C^{\bullet}(K^{p},\mathbf{D}_{\Omega}^{s})_{\leq h}\otimes_{A(\Omega)}^{\mathbf{L}}A(\Sigma)\simeq C^{\bullet}(K^{p},\mathbf{D}_{\Sigma}^{s})_{\leq h}\]
of $\mathbf{T}(K^{p})$-module complexes in $\mathbf{D}^{b}(A(\Omega))$,
and the result follows analogously from the hypertor spectral sequence\[
\mathrm{Tor}_{-i}^{R}(H^{j}(C^{\bullet}),N)\Rightarrow\mathbf{T}\mathrm{or}_{-i-j}^{R}(C^{\bullet},N).\]

\textbf{Remark 3.1.1. }If $(\Omega,h)$ is a slope datum, $\Sigma_{1}$
is Zariski-closed in $\Omega$, and $\Sigma_{2}$ is\textbf{ }Zariski-closed
in $\Sigma_{1}$, the transitivity of the derived tensor product yields
an isomorphism\begin{eqnarray*}
C^{\bullet}(K^{p},\mathcal{D}_{\Sigma_{2}})_{\leq h} & \simeq & C^{\bullet}(K^{p},\mathcal{D}_{\Omega})_{\leq h}\otimes_{A(\Omega)}^{\mathbf{L}}A(\Sigma_{2})\\
 & \simeq & C^{\bullet}(K^{p},\mathcal{D}_{\Omega})_{\leq h}\otimes_{A(\Omega)}^{\mathbf{L}}A(\Sigma_{1})\otimes_{A(\Sigma_{1})}^{\mathbf{L}}A(\Sigma_{2})\\
 & \simeq & C^{\bullet}(K^{p},\mathcal{D}_{\Sigma_{1}})_{\leq h}\otimes_{A(\Sigma_{1})}^{\mathbf{L}}A(\Sigma_{2})\end{eqnarray*}
which induces a relative version of the Tor spectral sequence, namely\[
E_{2}^{i,j}=\mathrm{Tor}_{A(\Sigma_{1})}^{i}(H^{j}(K^{p},\mathcal{D}_{\Sigma_{1}})_{\leq h},A(\Sigma_{2}))\Rightarrow H^{i+j}(K^{p},\mathcal{D}_{\Sigma_{2}})_{\leq h}.\]
This spectral sequence plays an important role in Newton's Appendix.

\paragraph{The boundary and Borel-Moore/compactly supported spectral sequences}

Recall the complex $C_{\bullet}(K^{p},-)$ was defined by choosing
a triangulation of the Borel-Serre compactification $\overline{Y(K^{p}I)}$
of $Y(K^{p}I)$. The induced triangulation of the boundary yields
a complex $C_{\bullet}^{\partial}(K^{p},-)$ which calculates $H_{\ast}(\partial\overline{Y(K^{p}I)},M)$,
together with a morphism $\phi:C_{\bullet}^{\partial}(K^{p},-)\to C_{\bullet}(K^{p},-)$
inducing the usual morphism on homology. The boundary and Borel-Moore/compactly
supported sequences, and the morphisms between them, follow from beholding
the diagram\[
\xymatrix{\mathbf{R}\mathrm{Hom}_{A(\Omega)}(C_{\bullet}^{\mathrm{BM}}(K^{p},\mathbf{A}_{\Omega}^{s})_{\leq h},A(\Sigma))\ar[r]\ar[d] & C_{c}^{\bullet}(K^{p},\mathbf{D}_{\Sigma}^{s})_{\leq h}\ar[d]\\
\mathbf{R}\mathrm{Hom}_{A(\Omega)}(C_{\bullet}(K^{p},\mathbf{A}_{\Omega}^{s})_{\leq h},A(\Sigma))\ar[r]\ar[d] & C^{\bullet}(K^{p},\mathbf{D}_{\Sigma}^{s})_{\leq h}\ar[d]\\
\mathbf{R}\mathrm{Hom}_{A(\Omega)}(C_{\bullet}^{\partial}(K^{p},\mathbf{A}_{\Omega}^{s})_{\leq h},A(\Sigma))\ar[r] & C_{\partial}^{\bullet}(K^{p},\mathbf{D}_{\Sigma}^{s})_{\leq h}}
\]
in which the horizontal arrows are quasi-isomorphisms, the columns
are exact triangles in $\mathbf{D}^{b}(A(\Omega))$, and the diagram
commutes up to homotopy for the natural action of $\mathbf{T}(K^{p})$.

\subsection{A construction of eigenvarieties}

In this section we briefly sketch the construction of the rigid spaces
$\mathscr{X}(K^{p})$ described in Theorem 1.3. 

\textbf{Step One: the spectral variety and its covering. }Fix a controlling
operator $U=U_{t}$. Given $\Omega\subset\mathcal{W}$ an admissible
open affinoid and $s\geq s[\Omega]$, let $F_{\Omega}(T)=\sum_{n=0}^{\infty}a_{n,\Omega}T^{n}\in A(\Omega)\{\{T\}\}$
be the characteristic power series of $\tilde{U}$ acting on the complex
$C_{\bullet}(K^{p},\mathbf{A}_{\Omega}^{s})$. Suppose $\Omega'\subset\Omega$;
since $F_{\Omega}(T)|_{\Omega'}=F_{\Omega'}(T)$, a simple calculation
gives $a_{n,\Omega}|_{\Omega'}=a_{n,\Omega'}$. By Tate's acyclicity
theorem there exist unique elements $a_{n}\in\mathcal{O}(\mathcal{W})$
such that $F(T)=\sum_{n=0}^{\infty}a_{n}T^{n}\in\mathcal{O}(\mathcal{W})\{\{T\}\}$
restricts to $F_{\Omega}(T)$ on any admissible affinoid open subset
$\Omega\subset\mathcal{W}$. The zero locus of $F(T)$ cuts out a
Fredholm hypersurface $\mathcal{Z}\subset\mathcal{W}\times\mathbf{A}^{1}$,
and we regard $\mathcal{Z}$ as a locally G-ringed space in the usual
way (cf. \cite{Conirred} for a thorough treatment of Fredholm hypersurfaces).
Given a connected affinoid $\Omega\subset\mathcal{W}$ and a rational
$h\in\mathbf{Q}_{\geq0}$, define the admissible open \[
\mathcal{Z}_{\Omega,h}=\mathcal{Z}\cap\left\{ (\omega,z),\omega\in\Omega,z\in\mathbf{A}^{1}\,\mathrm{with}\,|z|\leq p^{h}\right\} ;\]
the coordinate ring is naturally $A(\mathcal{Z}_{\Omega,h})=A(\Omega)\left\langle p^{h}T\right\rangle /(F_{\Omega}(T))$.
The morphism $\mathcal{Z}_{\Omega,h}\to\Omega$ is flat but not necessarily
finite. Let $\mathscr{C}\mathrm{ov}$ be the set of $\mathcal{Z}_{\Omega,h}$'s
such that $\mathcal{Z}_{\Omega,h}\to\Omega$ is \emph{finite flat}.
The significance of the finiteness condition is that $\mathcal{Z}_{\Omega,h}$
is finite if and only if $C_{\bullet}(K^{p},\mathbf{A}_{\Omega}^{s})$
admits a slope-$\leq h$ decomposition: $\mathcal{Z}_{\Omega,h}$
is finite flat if and only if it is disconnected from its complement
in $\mathcal{Z}_{\Omega,\infty}$, if and only if $F_{\Omega}$ admits
a slope-$\leq h$ factorization $F_{\Omega}=Q_{\Omega,h}\cdot R_{\Omega,h}$,
in which case $A(\mathcal{Z}_{\Omega,h})\simeq A(\Omega)[T]/Q_{\Omega,h}(T)$.
Since $C_{\bullet}(K^{p},\mathbf{A}_{\Omega}^{s})_{\leq h}$ is the
kernel of $Q^{\ast}(\tilde{U})$, and $\tilde{U}$ acts invertibly
on $C_{\bullet}(K^{p},\mathbf{A}_{\Omega}^{s})_{\leq h}$, the map
\begin{eqnarray*}
A(\Omega)[T] & \to & \mathrm{End}_{A(\Omega)}(C_{\bullet}(K^{p},\mathbf{A}_{\Omega}^{s})_{\leq h})\\
T & \mapsto & \tilde{U}^{-1}\end{eqnarray*}
puts a canonical $A(\mathcal{Z}_{\Omega,h})$-module structure on
$C_{\bullet}(K^{p},\mathbf{A}_{\Omega}^{s})_{\leq h}$; by the isomorphisms
proven in §2.4, this map induces a canonical $A(\mathcal{Z}_{\Omega,h})$-module
structure on $C_{\bullet}(K^{p},\mathcal{A}_{\Omega})_{\leq h}$ and
on $C^{\bullet}(K^{p},\mathcal{D}_{\Omega})_{\leq h}$. Note that
by Proposition 2.3.2, \emph{every }point $z\in\mathcal{Z}$ is contained
in some $\mathcal{Z}_{\Omega,h}$ with $\mathcal{Z}_{\Omega,h}\in\mathscr{C}\mathrm{ov}$.
By a fundamental theorem of Buzzard (Lemma 4.5 and Theorem 4.6 of \cite{BuEigen}),
the elements of $\mathscr{C}\mathrm{ov}$ form an admissible covering
of $\mathcal{Z}$. 

\textbf{Step Two: the cocycle condition. }If $\mathcal{Z}_{\Omega,h}\in\mathscr{C}\mathrm{ov}$
and $\Omega'\subset\Omega$ with $\Omega'$ connected, a moment's
thought yields $A(\mathcal{Z}_{\Omega',h})\simeq A(\mathcal{Z}_{\Omega,h})\otimes_{A(\Omega)}A(\Omega')$,
so then $\mathcal{Z}_{\Omega',h}\in\mathscr{C}\mathrm{ov}$. Fix $\mathcal{Z}_{\Omega',h'}\subseteq\mathcal{Z}_{\Omega,h}\in\mathscr{C}\mathrm{ov}$
with $\mathcal{Z}_{\Omega',h'}\in\mathscr{C}\mathrm{ov}$; we necessarily
have $\Omega'\subseteq\Omega$, and we may assume $h'\leq h$. Set
$C_{\Omega,h}=C^{\bullet}(K^{p},\mathcal{D}_{\Omega})_{\leq h}$.
We now trace through the following sequence of canonical isomorphisms:\begin{eqnarray*}
C_{\Omega,h}\otimes_{A(\mathcal{Z}_{\Omega,h})}A(\mathcal{Z}_{\Omega',h'}) & \simeq & C_{\Omega,h}\otimes_{A(\mathcal{Z}_{\Omega,h})}A(\mathcal{Z}_{\Omega',h})\otimes_{A(\mathcal{Z}_{\Omega',h})}A(\mathcal{Z}_{\Omega',h'})\\
 & \simeq & \left(C_{\Omega,h}\otimes_{A(\mathcal{Z}_{\Omega,h})}A(\mathcal{Z}_{\Omega,h})\otimes_{A(\Omega)}A(\Omega')\right)\otimes_{A(\mathcal{Z}_{\Omega',h})}A(\mathcal{Z}_{\Omega',h'})\\
 & \simeq & \left(C_{\Omega,h}\otimes_{A(\Omega)}A(\Omega')\right)\otimes_{A(\mathcal{Z}_{\Omega',h})}A(\mathcal{Z}_{\Omega',h'})\\
 & \simeq & C_{\Omega',h}\otimes_{A(\mathcal{Z}_{\Omega',h})}A(\mathcal{Z}_{\Omega',h'})\\
 & \simeq & C_{\Omega',h'}.\end{eqnarray*}
The fourth line here follows from Proposition 2.4.5.

\textbf{Step Three: coherent $\mathcal{O}$-modules. }Given $\mathcal{Z}_{\Omega',h'}\subseteq\mathcal{Z}_{\Omega,h}\in\mathscr{C}\mathrm{ov}$
with $\mathcal{Z}_{\Omega',h'}\in\mathscr{C}\mathrm{ov}$, $\mathcal{Z}_{\Omega',h'}$
is an affinoid subdomain of $\mathcal{Z}_{\Omega,h}$, so $A(\mathcal{Z}_{\Omega',h'})$
is $A(\mathcal{Z}_{\Omega,h})$-flat. In particular, the functor $-\otimes_{A(\mathcal{Z}_{\Omega,h})}A(\mathcal{Z}_{\Omega',h'})$
commutes with taking cohomology of any complex of $A(\mathcal{Z}_{\Omega,h})$-modules.
Hence, taking cohomology in the isomorphism of step three yields a
canonical isomorphism\[
H^{\ast}(K^{p},\mathcal{D}_{\Omega})_{\leq h}\otimes_{A(\mathcal{Z}_{\Omega,h})}A(\mathcal{Z}_{\Omega',h'})\simeq H^{\ast}(K^{p},\mathcal{D}_{\Omega'})_{\leq h'}.\]
This is exactly what we need in order to verify that the assignments
\[
\mathcal{Z}_{\Omega,h}\mapsto H^{\ast}(K^{p},\mathcal{D}_{\Omega})_{\leq h},\,\mathcal{Z}_{\Omega,h}\in\mathscr{C}\mathrm{ov}\]
glue together into a coherent $\mathcal{O}_{\mathcal{Z}}$-module
sheaf $\mathcal{M}^{\ast}=\oplus\mathcal{M}^{n}$ such that $\mathcal{M}^{n}(\mathcal{Z}_{\Omega,h})\simeq H^{n}(K^{p},\mathcal{D}_{\Omega})_{\leq h}$
for $\mathcal{Z}_{\Omega,h}\in\mathscr{C}\mathrm{ov}$. Now, for $\mathcal{Z}_{\Omega,h}\in\mathscr{C}\mathrm{ov}$
let $\mathbf{T}_{\Omega,h}$ be the commutative subalgebra of $\mathrm{End}_{A(\mathcal{Z}_{\Omega,h})}(H^{\ast}(K^{p},\mathcal{D}_{\Omega})_{\leq h})$
generated by $\mathbf{T}(K^{p})\otimes_{\mathbf{Q}_{p}}A(\Omega)$,
with \[
\psi_{\Omega,h}:\mathbf{T}(K^{p})\otimes_{\mathbf{Q}_{p}}A(\Omega)\to\mathbf{T}_{\Omega,h}\]
the structure map. For any Noetherian ring $A$, finite $A$-module
$M$, and flat $A$-algebra $B$, there is a canonical isomorphism
$\mathrm{End}_{A}(M)\otimes_{A}B\simeq\mathrm{End}_{B}(M\otimes_{A}B)$.
In particular, we immediately obtain canonical isomorphisms\[
\mathbf{T}_{\Omega,h}\otimes_{A(\mathcal{Z}_{\Omega,h})}A(\mathcal{Z}_{\Omega',h'})\simeq\mathbf{T}_{\Omega',h'},\]
whereby the $A(\mathcal{Z}_{\Omega,h})$-algebras $\mathbf{T}_{\Omega,h}$
glue together into a coherent sheaf $\mathcal{T}$ of $\mathcal{O}_{\mathcal{Z}}$-algebras.
By an obvious rigid version of the {}``relative Spec'' construction
(see e.g. §2.2 of \cite{Conrelative}) the affinoid rigid spaces $\mathscr{X}_{\Omega,h}=\mathrm{Sp}\mathbf{T}_{\Omega,h}$
glue together into a rigid space $\mathscr{X}$, and the natural morphisms
\[
w:\mathscr{X}_{\Omega,h}\to\mathcal{Z}_{\Omega,h}\to\Omega\]
glue into a morphism $w:\mathscr{X}\to\mathcal{W}$. Since $\mathcal{Z}$
is separated and $\mathscr{X}$ is finite over $\mathcal{Z}$, $\mathscr{X}$
is separated. Each $H^{\ast}(K^{p},\mathcal{D}_{\Omega})_{\leq h}$
is naturally a finite $\mathbf{T}_{\Omega,h}$-module; since\begin{eqnarray*}
H^{\ast}(K^{p},\mathcal{D}_{\Omega})_{\leq h}\otimes_{\mathbf{T}_{\Omega,h}}\mathbf{T}_{\Omega',h'} & = & H^{\ast}(K^{p},\mathcal{D}_{\Omega})_{\leq h}\otimes_{\mathbf{T}_{\Omega,h}}\mathbf{T}_{\Omega,h}\otimes_{A(\mathcal{Z}_{\Omega,h})}A(\mathcal{Z}_{\Omega',h'})\\
 & \simeq & H^{\ast}(K^{p},\mathcal{D}_{\Omega'})_{\leq h'},\end{eqnarray*}
the cohomology groups $H^{\ast}(K^{p},\mathcal{D}_{\Omega})_{\leq h}$
glue together into a sheaf of coherent $\mathcal{O}_{\mathscr{X}}$-modules.
For any fixed element $T\in\mathbf{T}(K^{p})$, the sections $\psi_{\Omega,h}(T\otimes1)\in\mathcal{O}(\mathscr{X}_{\Omega,h})$
glue to a unique global section $\psi(T)\in\mathcal{O}(\mathscr{X})$,
and $\psi$ is easily seen to be an algebra homomorphism. Power-boundedness
of $\psi(T)$ follows immediately from the fact that any element $\gamma\in\Delta$
acts on $\mathcal{D}_{\Omega}$ with operator norm at most one for
the family of norms defining the Fréchet topology on $\mathcal{D}_{\Omega}$.
At this point, claims i., iii., and iv. of Theorem 1.5 are immediate
consequences of our construction. If $x\in\mathscr{X}(\overline{\mathbf{Q}_{p}})$,
choose some $\mathscr{X}_{\Omega,h}$ with $w(x)\in\Omega$ and $h>v_{p}(\psi(x)(U_{t}))$;
writing $\mathfrak{M}_{x}$ for the maximal ideal of $\mathbf{T}_{\Omega,h}$
determined by $x$, claim$ $ ii. follows immediately from setting
$\Sigma=w(x)$ in the Tor spectral sequence and then localizing the
spectral sequence at $\mathfrak{M}_{x}$. Claim v. is straightforward,
though slightly tedious, and we point the reader to \cite{Hthesis}
for details.

\subsection{The support of overconvergent cohomology}

Let $R$ be a Noetherian ring, and let $M$ be a finite $R$-module.
We say $M$ has \emph{full support }if $\mathrm{Supp}(M)=\mathrm{Spec}(R)$,
and that $M$ is \emph{torsion }if $\mathrm{ann}(M)\neq0$. We shall
repeatedly use the following basic result.

\textbf{Proposition 3.3.1. }\emph{If $\mathrm{Spec}(R)$ is reduced
and irreducible, the following are equivalent:}

\emph{i) $M$ is faithful (i.e. $\mathrm{ann}(M)=0$),}

\emph{ii) $M$ has full support,}

\emph{iii)} \emph{$M$ has nonempty open support,}

\emph{iv) $\mathrm{Hom}_{R}(M,R)\neq0$,}

\emph{v) $M\otimes_{R}K\ne0$, $K=\mathrm{Frac}(R)$.}

\emph{Proof. }Since $M$ is finite, $\mathrm{Supp}(M)$ is the underlying
topological space of $\mathrm{Spec}(R/\mathrm{ann}(M))$, so i) obviously
implies ii). If $\mathrm{Spec}(R/\mathrm{ann}(M))=\mathrm{Spec}(R)$
as topological spaces, then $\mathrm{ann}(M)\subset\surd(0)=(0)$
since $R$ is reduced, so ii) implies i). The set $\mathrm{Supp}(M)=\mathrm{Spec}(R/\mathrm{ann}(M))$
is \emph{a priori }closed; since $\mathrm{Spec}(R)$ is irreducible
by assumption, the only nonempty simultaneously open and closed subset
of $\mathrm{Spec}(R)$ is all of $\mathrm{Spec}(R)$, so ii) and iii)
are equivalent. By finiteness, $M$ has full support if and only if
$(0)$ is an associated prime of $M$, if and only if there is an
injection $R\hookrightarrow M$; tensoring with $K$ implies the equivalence
of ii) and v). Finally, $\mathrm{Hom}_{R}(M,R)\otimes_{R}K\simeq\mathrm{Hom}_{K}(M\otimes_{R}K,K)$,
so $M\otimes_{R}K\neq0$ if and only if $\mathrm{Hom}_{R}(M,R)\neq0$,
whence iv) and v) are equivalent. $\square$

\emph{Proof of Theorem 1.2.i}. (I'm very grateful to Jack Thorne for
suggesting this proof.) Tensoring the Ext spectral sequence with $K(\Omega)=\mathrm{Frac}(A(\Omega))$,
it degenerates to isomorphisms \[
\mathrm{Hom}_{K(\Omega)}(H_{i}(K^{p},\mathcal{A}_{\Omega})_{\leq h}\otimes_{A(\Omega)}K(\Omega),K(\Omega))\simeq H^{i}(K^{p},\mathcal{D}_{\Omega})_{\leq h}\otimes_{A(\Omega)}K(\Omega),\]
so the claim is immediate from the preceding proposition. $\square$

\emph{Proof of Theorem 1.2.ii. }We give the proof in two steps, with
the first step naturally breaking into two cases. In the first step,
we prove the result assuming $\Omega$ contains an arithmetic weight.
In the second step, we eliminate this assumption via analytic continuation.

\paragraph*{Step One, Case One: G doesn't have a discrete series.}

Let $\mathcal{W}^{\mathrm{sd}}$ be the rigid Zariski closure in $\mathcal{W}$
of the arithmetic weights whose algebraic parts are the highest weights
of irreducible $G$-representations with nonvanishing $(\mathfrak{g},K_{\infty})$-cohomology.
A simple calculation using §II.6 of \cite{BWcohom} shows that $\mathcal{W}^{\mathrm{sd}}$
is the union of its countable set of irreducible components, each
of dimension $r(G)$. An arithmetic weight is \emph{non-self-dual
}if $\lambda\notin\mathcal{W}^{\mathrm{sd}}$.

Now, by assumption $\Omega$ contains an arithmetic weight, so $\Omega$
automatically contains a Zariski dense $ $set $\mathcal{N}_{h}\subset\Omega\smallsetminus\Omega\cap\mathcal{W}^{\mathrm{sd}}$
of non-self-dual arithmetic weights for which $h$ is a small slope.
By Theorem 2.4.8 together with Matsushima's formula, $H^{\ast}(K^{p},\mathcal{D}_{\lambda})_{\leq h}$
vanishes identically for any $\lambda\in\mathcal{N}_{h}$. For any
fixed $\lambda\in\mathcal{N}_{h}$, suppose $\mathfrak{m}_{\lambda}\in\mathrm{Supp}_{\Omega}H^{\ast}(K^{p},\mathcal{D}_{\Omega})_{\leq h}$;
let $d$ be the largest integer with $\mathfrak{m}_{\lambda}\in\mathrm{Supp}_{\Omega}H^{d}(K^{p},\mathcal{D}_{\Omega})_{\leq h}$.
Taking $\Sigma=\lambda$ in the Tor spectral sequence gives\[
E_{2}^{i,j}=\mathrm{Tor}_{-i}^{A(\Omega)}(H^{j}(K^{p},\mathcal{D}_{\Omega})_{\leq h},A(\Omega)/\mathfrak{m}_{\lambda})\Rightarrow H^{i+j}(K^{p},\mathcal{D}_{\lambda})_{\leq h}.\]
The entry $E_{2}^{0,d}=H^{d}(K^{p},\mathcal{D}_{\Omega})_{\leq h}\otimes_{A(\Omega)}A(\Omega)/\mathfrak{m}_{\lambda}$
is nonzero by Nakayama's lemma, and is stable since every row of the
$E_{2}$-page above the \emph{d}th row vanishes by assumption. In
particular, $E_{2}^{0,d}$ contributes a nonzero summand to the grading
on $H^{d}(K^{p},\mathcal{D}_{\lambda})_{\leq h}$ - but this module
is zero, contradicting our assumption that $\mathfrak{m}_{\lambda}\in\mathrm{Supp}_{\Omega}H^{\ast}(K^{p},\mathcal{D}_{\Omega})$.
Therefore, $H^{\ast}(K^{p},\mathcal{D}_{\Omega})_{\leq h}$ does \emph{not
}have full support, so is not a faithful $A(\Omega)$-module.

\paragraph*{Step One, Case Two: G has a discrete series.}

The idea is the same as Case One, but with $\mathcal{N}_{h}$ replaced
by $\mathcal{R}_{h}$, the set of arithmetic weights with \emph{regular
}algebraic part for which $h$ is a small slope. For these weights,
Proposition 2.4.8 and Matsushima's formula together with known results
on $(\mathfrak{g},K_{\infty})$-cohomology (see e.g. Sections 4-5
of \cite{LSEiscohom}) implies that $H^{i}(K^{p},\mathcal{D}_{\lambda})_{\leq h},\,\lambda\in\mathcal{R}_{h}$\emph{
}vanishes for $i\neq d_{G}=\frac{1}{2}\mathrm{dim}G(\mathbf{R})/Z_{\infty}K_{\infty}$.
The Tor spectral sequence with $\Sigma=\lambda\in\mathcal{R}_{h}$
then shows that $\mathcal{R}_{h}$ doesn't meet $\mathrm{Supp}_{\Omega}H^{i}(K^{p},\mathcal{D}_{\Omega})_{\leq h}$
for any $i>d_{G}$. The Ext spectral sequence with $\Sigma=\lambda\in\mathcal{R}_{h}$
then shows that $\mathcal{R}_{h}$ doesn't meet $\mathrm{Supp}_{\Omega}H_{i}(K^{p},\mathcal{A}_{\Omega})_{\leq h}$
for any $i<d_{G}$, whence the Ext spectral sequence with $\Sigma=\Omega$
shows that $\mathcal{R}_{h}$ doesn't meet $\mathrm{Supp}_{\Omega}H^{i}(K^{p},\mathcal{D}_{\Omega})_{\leq h}$
for any $i<d_{G}$. The result follows.

\paragraph*{Step Two.}

We maintain the notation of §3.2. As in that subsection, $H^{n}(K^{p},\mathcal{D}_{\Omega})_{\leq h}$
glues together over the affinoids $\mathcal{Z}_{\Omega,h}\in\mathscr{C}\mathrm{ov}$
into a coherent $\mathcal{O_{Z}}$-module sheaf $\mathcal{M}^{n}$,
and in particular, the support of $\mathcal{M}^{n}$ is a closed analytic
subset of $\mathcal{Z}$. Let $\pi:\mathcal{Z}\to\mathcal{W}$ denote
the natural projection. For any $\mathcal{Z}_{\Omega,h}\in\mathscr{C}\mathrm{ov}$,
we have \[
\pi_{\ast}\mathrm{Supp}_{\mathcal{Z}_{\Omega,h}}\mathcal{M}^{n}(\mathcal{Z}_{\Omega,h})=\mathrm{Supp}_{\Omega}H^{n}(K^{p},\mathcal{D}_{\Omega})_{\leq h}.\]
Suppose $\mathrm{Supp}_{\Omega}H^{n}(K^{p},\mathcal{D}_{\Omega})_{\leq h}=\Omega$
for some $\mathcal{Z}_{\Omega,h}\in\mathscr{C}\mathrm{ov}$. This
implies that $\mathrm{Supp}_{\mathcal{Z}_{\Omega,h}}\mathcal{M}^{n}(\mathcal{Z}_{\Omega,h})$
contains a closed subset of dimension equal to $\mathrm{dim}\mathcal{Z}$,
so contains an irreducible component of $\mathcal{Z}_{\Omega,h}$.
Any irreducible component of $\mathcal{Z}_{\Omega,h}$ is an admissible
open affinoid in $\mathcal{Z}$, so this in turn implies that $\mathrm{Supp}_{\mathcal{Z}}\mathcal{M}^{n}$
contains an affinoid open. Since $\mathrm{Supp}_{\mathcal{Z}}\mathcal{M}^{n}$
is a priori closed, Corollary 2.2.6 of \cite{Conirred} implies that
$\mathrm{Supp}_{\mathcal{Z}}\mathcal{M}^{n}$ contains an entire irreducible
component of $\mathcal{Z}$, say $\mathcal{Z}_{0}$. The irreducible
component $\mathcal{Z}_{0}$ corresponds, by Theorem 4.3.2 of \cite{Conirred},
to a nonconstant irreducible Fredholm series \[
F_{0}(T)=1+\sum_{j=1}^{\infty}\omega_{j}T^{j},\,\omega_{j}\in\mathcal{O}(\mathcal{W})\]
dividing $F(T)$. I claim the image of $\mathcal{Z}_{0}$ under $\pi$
is Zariski-open in $\mathcal{W}$. Indeed, by Lemma 1.3.2 of \cite{CMeigencurve},
the fiber $\mathcal{Z}_{0}\cap\pi^{-1}(\lambda)$ is empty for a given
$\lambda\in\mathcal{W}$ if and only if $\omega_{j}\in\mathfrak{m}_{\lambda}$
for all $j$, if and only if $\mathscr{I}=(\omega_{1},\omega_{2},\omega_{3},\dots)\subset\mathfrak{m}_{\lambda}$.
The ideal $\mathscr{I}\subset\mathcal{O}(\mathcal{W})$ is naturally
identified with the global sections of a coherent ideal sheaf over
$\mathcal{W}$, which cuts out a closed analytic subset $V(\mathscr{I})$
in the usual way; the complement of $V(\mathscr{I})$ is precisely
$\pi_{\ast}\mathcal{Z}_{0}$. Fix an arithmetic weight $\lambda_{0}\in\pi_{\ast}\mathcal{Z}_{0}$.
For some sufficiently large $h_{0}$ and some affinoid $\Omega_{0}$
containing $\lambda_{0}$, $\mathcal{Z}_{\Omega_{0},h_{0}}$ will
contain $\mathcal{Z}_{\Omega_{0},h_{0}}\cap\mathcal{Z}_{0}$ as a
nonempty union of irreducible components, and the latter intersection
will be finite flat over $\Omega_{0}$. Since $\mathcal{M}^{n}(\mathcal{Z}_{\Omega_{0},h_{0}})\simeq H^{n}(K^{p},\mathcal{D}_{\Omega_{0}})_{\leq h_{0}}$,
we deduce that $\mathrm{Supp}_{\Omega_{0}}H^{n}(K^{p},\mathcal{D}_{\Omega_{0}})_{\leq h_{0}}=\Omega_{0}$,
whence $H^{n}(K^{p},\mathcal{D}_{\Omega_{0}})_{\leq h_{0}}$ is faithful,
so by Step One $G^{\mathrm{der}}(\mathbf{R})$ has a discrete series
and $n=\frac{1}{2}\mathrm{dim}G(\mathbf{R})/Z_{\infty}K_{\infty}$.

\subsection{Some cases of Urban's conjecture}

In this subsection we prove Theorem 1.6.

By the basic properties of irreducible components together with the
construction given in §3.2, it suffices to work locally over a fixed
$\mathcal{Z}_{\Omega,h}\in\mathscr{C}\mathrm{ov}$. Suppose $\phi:\mathbf{T}_{\Omega,h}\to\overline{\mathbf{Q}_{p}}$
is an eigenpacket corresponding to a cuspidal non-critical regular
classical point $x\in\mathscr{X}_{\Omega,h}$. Set $\mathfrak{M}=\ker\phi$,
and let $\mathfrak{m}=\mathfrak{m}_{\lambda}$ be the contraction
of $\mathfrak{M}$ to $A(\Omega)$. Let $\mathscr{P}\subset\mathbf{T}_{\Omega,h}$
be any minimal prime contained in $\mathfrak{M}$, and let $\wp$
be its contraction to a prime in $A(\Omega)$. The ring $\mathbf{T}_{\Omega,h}/\mathscr{P}$
is a finite integral extension of $A(\Omega)/\wp$, so both rings
have the same dimension. In particular, Conjecture 1.5 follows from
the equality $\mathrm{ht}\wp=l(G)$; this latter equality is what
we prove.

\textbf{Proposition 3.4.1. }\emph{The largest degrees for which $\phi$
occurs in $H^{\ast}(K^{p},\mathcal{D}_{\Omega})_{\leq h}$ and $H^{\ast}(K^{p},\mathcal{D}_{\lambda})_{\leq h}$
coincide, and the smallest degrees for which $\phi$ occurs in $H^{\ast}(K^{p},\mathcal{D}_{\lambda})_{\leq h}$
and $H_{\ast}(K^{p},\mathcal{A}_{\Omega})_{\leq h}$ coincide. Finally,
the smallest degree for which $\phi$ occurs in $H^{\ast}(K^{p},\mathcal{D}_{\Omega})_{\leq h}$
is greater than or equal to the smallest degree for which $\phi$
occurs in $H_{\ast}(K^{p},\mathcal{A}_{\Omega})_{\leq h}$.}

\emph{Proof. }For the first claim, localize the Tor spectral sequence
at $\mathfrak{M}$, with $\Sigma=\lambda$. If $\phi$ occurs in $H^{i}(K^{p},\mathcal{D}_{\lambda})_{\leq h}$
then it occurs in a subquotient of $\mathrm{Tor}_{j}^{A(\Omega)}(H^{i+j}(K^{p},\mathcal{D}_{\Omega})_{\leq h},A(\Omega)/\mathfrak{m}_{\lambda})$
for some $j\geq0$. On the other hand, if $d$ is the largest degree
for which $\phi$ occurs in $H^{d}(K^{p},\mathcal{D}_{\Omega})_{\leq h}$,
the entry $E_{2}^{0,d}$ of the Tor spectral sequence is stable and
nonzero after localizing at $\mathfrak{M}$, and it contributes to
the grading on $H^{d}(K^{p},\mathcal{D}_{\lambda})_{\leq h,\mathfrak{M}}$.
The second and third claims follow from an analogous treatment of
the Ext spectral sequence. $\square$

First we treat the case where $l(G)=0$, so $G^{\mathrm{der}}(\mathbf{R})$
has a discrete series. By the noncriticality of $\phi$ together with
the nonexistence of CAP forms at regular weights and the results recalled
in §3.3, the only degree for which $\phi$ occurs in $H^{i}(K^{p},\mathcal{D}_{\lambda})_{\leq h}$
is the middle degree $d=\frac{1}{2}\mathrm{dim}G(\mathbf{R})/Z_{\infty}K_{\infty}$,
so Proposition 3.4.1 implies that the only degree for which $\phi$
occurs in $H^{\ast}(K^{p},\mathcal{D}_{\Omega})_{\leq h}$ is the
middle degree as well. The Tor spectral sequence localized at $\mathfrak{M}$
now degenerates, and yields\[
\mathrm{Tor}_{i}^{A(\Omega)_{\mathfrak{m}}}(H^{d}(K^{p},\mathcal{D}_{\Omega})_{\leq h,\mathfrak{M}},A(\Omega)/\mathfrak{m})=0\,\mathrm{for\, all}\, i\geq1.\]
By Proposition A.3, $H^{d}(K^{p},\mathcal{D}_{\Omega})_{\leq h,\mathfrak{M}}$
is a \emph{free }module over $A(\Omega)_{\mathfrak{m}}$, so $\mathrm{Ass}_{A(\Omega)_{\mathfrak{m}}}(H^{d}(K^{p},\mathcal{D}_{\Omega})_{\leq h,\mathfrak{M}})=\{(0)\}$.
Let us write $A=A(\Omega)_{\mathfrak{m}}$ and $\mathbf{T}=(\mathbf{T}_{\Omega,h})_{\mathfrak{M}}$
for brevity, so there is a diagram \[
A\to\mathbf{T}\hookrightarrow\mathrm{End}_{A}\left(H^{d}(K^{p},\mathcal{D}_{\Omega})_{\leq h,\mathfrak{M}}\right).\]
I claim the structure map $A\to\mathbf{T}$ is a \emph{flat }local
homomorphism. To see this, note by the isomorphism\[
H^{d}(K^{p},\mathcal{D}_{\Omega})_{\leq h,\mathfrak{M}}\otimes_{A}A/\mathfrak{m}\simeq H^{d}(K^{p},\mathcal{D}_{\lambda})_{\leq h,\mathfrak{M}}\]
that the residue ring $\mathbf{T}/\mathfrak{m}\mathbf{T}$ is simply
the subalgebra of $\mathrm{End}_{A/\mathfrak{m}A}\left(H^{d}(K^{p},\mathcal{D}_{\lambda})_{\leq h,\mathfrak{M}}\right)$
generated by $\mathbf{T}(K^{p})$. By the non-criticality of $x$
together with the semisimplicity of the spherical Hecke algebra, $\mathbf{T}/\mathfrak{m}\mathbf{T}$
is a \emph{field }and $H^{d}(K^{p},\mathcal{D}_{\lambda})_{\leq h,\mathfrak{M}}$
is a free $\mathbf{T}/\mathfrak{m}\mathbf{T}$-module of finite rank.
Choosing lifts of a minimal set of $\mathbf{T}/\mathfrak{m}\mathbf{T}$-module
generators of $H^{d}(K^{p},\mathcal{D}_{\lambda})_{\leq h,\mathfrak{M}}$
determines a surjection \[
\psi:\mathbf{T}^{r}\twoheadrightarrow H^{d}(K^{p},\mathcal{D}_{\Omega})_{\leq h,\mathfrak{M}}\]
of $A$-modules. Applying $-\otimes_{A}A/\mathfrak{m}A$ to the sequence\[
0\to\ker\psi\to\mathbf{T}^{r}\to H^{d}(K^{p},\mathcal{D}_{\Omega})_{\leq h,\mathfrak{M}}\to0\]
and using the freeness of the final term over $A$, we easily see
that $(\ker\psi)\otimes_{A}A/\mathfrak{m}A=0$, so $\psi$ is an isomorphism.
Since $H^{d}(K^{p},\mathcal{D}_{\Omega})_{\leq h,\mathfrak{M}}$ is
free over $A$ and over $\mathbf{T}$, it follows that $\mathbf{T}$
is a free $A$-module. The rings $A$ and $\mathbf{T}/\mathfrak{m}\mathbf{T}$
are regular, so $\mathbf{T}$ is regular by Theorem 23.7.ii of \cite{Matcrt}.
Since $\mathbf{T}/\mathfrak{m}\mathbf{T}$ is a field, $w$ is étale
at $x$.

Now we turn to the case $l(G)\geq1$. First, there is an affinoid
open $\mathscr{Y}\subset\mathscr{X}_{\Omega,h}$ containing $x$,
and meeting every component of $\mathscr{X}_{\Omega,h}$ containing
$x$, such that every regular classical non-critical point in $\mathscr{Y}$
is cuspidal. Indeed, at regular weights there are no CAP representations,
so the Hecke action splits the cuspidal and Eisenstein subspaces of
$H^{\ast}(K^{p},V_{\lambda^{\mathrm{alg}}})$, whence the top horizontal
arrow in the diagram\[
\xymatrix{H_{c}^{i}(K^{p},V_{\lambda^{\mathrm{alg}}})_{\leq h}\ar[r] & H^{i}(K^{p},V_{\lambda^{\mathrm{alg}}})_{\leq h}\\
H_{c}^{i}(K^{p},\mathcal{D}_{\lambda})_{\le h}\ar[u]\ar[r] & H^{i}(K^{p},\mathcal{D}_{\lambda})_{\le h}\ar[u]}
\]
becomes an isomorphism after localizing at $\mathfrak{M}$. The vertical
arrows are isomorphisms by the noncriticality assumption, so the bottom
arrow becomes an isomorphism after localization at $\mathfrak{M}$.
Localizing the sequence\[
\dots\to H_{c}^{i}(K^{p},\mathcal{D}_{\lambda})\to H^{i}(K^{p},\mathcal{D}_{\lambda})\to H_{\partial}^{i}(K^{p},\mathcal{D}_{\lambda})\to H_{c}^{i+1}(K^{p},\mathcal{D}_{\lambda})\to\dots\]
at $\mathfrak{M}$ then shows that $\phi$ does not occur in $H_{\partial}^{\ast}(K^{p},\mathcal{D}_{\lambda})_{\leq h}$,
so by the boundary spectral sequence $\phi$ does not occur in $H_{\partial}^{\ast}(K^{p},\mathcal{D}_{\Omega})_{\leq h}$.
Since $\mathrm{Supp}_{\mathbf{T}_{\Omega,h}}H_{\partial}^{\ast}(K^{p},\mathcal{D}_{\Omega})_{\leq h}$
is closed in $\mathscr{X}_{\Omega,h}$ and does not meet $x$, the
existence of a suitable $\mathscr{Y}$ now follows easily. Shrinking
$\Omega$ and $\mathscr{Y}$ as necessary, we may assume that $A(\mathscr{Y})$
is finite over $A(\Omega)$, and thus $\mathscr{M}^{\ast}(\mathscr{Y})=H^{\ast}(K^{p},\mathcal{D}_{\Omega})_{\leq h}\otimes_{\mathbf{T}_{\Omega,h}}A(\mathscr{Y})$
is finite over $A(\Omega)$ as well. Exactly as in the proof\emph{
}of Theorem 1.2, the Tor spectral sequence shows that $\mathrm{Supp}_{\Omega}\mathscr{M}^{\ast}(\mathscr{Y})$
doesn't contain any regular non-self-dual weights for which $h$ is
a small slope, so $\mathscr{M}^{\ast}(\mathscr{Y})$ and $A(\mathscr{Y})$
are torsion $A(\Omega)$-modules.

Finally, suppose $l(G)=1$. Set $d=q(G)$, so $\phi$ occurs in $H^{\ast}(K^{p},\mathcal{D}_{\lambda})_{\leq h}\simeq H^{\ast}(K^{p},V_{\lambda})_{\leq h}$
only in degrees $d$ and $d+1$. By the argument of the previous paragraph
and the faithful flatness of $\mathcal{O}_{\Omega,\lambda}$ over
$A(\Omega)_{\mathfrak{m}}$, $H^{\ast}(K^{p},\mathcal{D}_{\Omega})_{\leq h,\mathfrak{M}}$
is a torsion $A(\Omega)_{\mathfrak{m}}$-module. Taking $\Sigma=\lambda$
in the Ext spectral sequence and localizing at $\mathfrak{M}$, Proposition
3.4.1 yields\[
H^{d}(K^{p},\mathcal{D}_{\Omega})_{\leq h,\mathfrak{M}}\simeq\mathrm{Hom}_{A(\Omega)_{\mathfrak{m}}}(H_{d}(K^{p},\mathcal{A}_{\Omega})_{\leq h,\mathfrak{M}},A(\Omega)_{\mathfrak{m}}).\]
Since the left-hand term is a torsion $A(\Omega)_{\mathfrak{m}}$-module,
Proposition 3.3.1 implies that both modules vanish identically. Proposition
3.4.1 now shows that $d+1$ is the only degree for which $\phi$ occurs
in $H^{\ast}(K^{p},\mathcal{D}_{\Omega})_{\leq h}$. Taking $\Sigma=\lambda$
in the Tor spectral sequence and localizing at $\mathfrak{M}$, the
only nonvanishing entries are $E_{2}^{0,d+1}$ and $E_{2}^{-1,d+1}$.
In particular, $\mathrm{Tor}_{i}^{A(\Omega)_{\mathfrak{m}}}(H^{d+1}(K^{p},\mathcal{D}_{\Omega})_{\leq h,\mathfrak{M}},A(\Omega)/\mathfrak{m})=0$
for all $i\geq2$, so $H^{d+1}(K^{p},\mathcal{D}_{\Omega})_{\leq h,\mathfrak{M}}$
has projective dimension at most one by Proposition A.3. Summarizing,
we've shown that $H^{i}(K^{p},\mathcal{D}_{\Omega})_{\leq h,\mathfrak{M}}$
vanishes in degrees $\neq q(G)+1$, and that $H^{q(G)+1}(K^{p},\mathcal{D}_{\Omega})_{\le h,\mathfrak{M}}$
is a torsion $A(\Omega)_{\mathfrak{m}}$-module of projective dimension
one, so the theorem now follows from Propositions A.4 and A.6. $\square$

As an exotic example of a group with $l(G)=1$, let $F$ be a totally
real number field with $\Sigma=\left\{ \sigma:F\hookrightarrow\mathbf{R}\right\} $,
and choose a quadratic form $Q$ in an even number of variables over
$F$ such that $\sigma(\mathrm{disc}(Q))<0$ for exactly one $\sigma\in\Sigma$.
Then $G=\mathrm{Res}_{F/\mathbf{Q}}\mathrm{SO}_{Q}$ has $l(G)=1$.

\section{The imaginary quadratic eigencurve}

In this section we sketch the proof of Theorem 1.8 in the nonsplit
case. The split case requires an analysis of boundary cohomology;
we defer an analysis of this case to \cite{Hbs}, since that paper
will develop tools for systematically calculating overconvergent cohomology
over the Borel-Serre boundary.

Fix an imaginary quadratic field $F/\mathbf{Q}$ and a prime $p$
split in $F$, say $p=\mathfrak{p}_{1}\mathfrak{p}_{2}$. Fix a quaternion
algebra $D/F$, and let $\mathfrak{d}\subset\mathcal{O}_{F}$ be the
product of the primes where $D$ ramifies. We assume $D$ splits at
every prime over $p$. Set $G=\mathrm{Res}_{F/\mathbf{Q}}D^{\times}$;
if $R$ is an associative $\mathbf{Q}$-algebra, we shall identify
$G(R)\simeq(D\otimes_{\mathbf{Q}}R)^{\times}$ without particular
comment. In particular, there is an isomorphism \[
\iota_{1}\times\iota_{2}:G(\mathbf{Q}_{p})=\mathrm{GL}_{2}(F\otimes_{\mathbf{Q}}\mathbf{Q}_{p})\simeq\mathrm{GL}_{2}(\mathbf{Q}_{p})\times\mathrm{GL}_{2}(\mathbf{Q}_{p}).\]
The weight space for $G$ is simply the product $\mathcal{W}^{(1)}\times\mathcal{W}^{(2)}$
of two copies of the weight space $\mathcal{W}$ for $\mathrm{GL}_{2}/\mathbf{Q}_{p}$.
We define the \emph{null space }$\mathcal{W}_{0}$ for $G$ as the
product of null subspaces $\mathcal{W}_{0}^{(1)}\times\mathcal{W}_{0}^{(2)}$,
so $\mathcal{W}_{0}$ is two-dimensional. In particular, a weight
$\lambda\in\mathcal{W}_{0}(L)$ is simply a pair of characters $\lambda_{1},\lambda_{2}:\mathbf{Z}_{p}^{\times}\to L^{\times}$
with $\lambda:\mathrm{diag}(x,1)\mapsto\lambda_{1}(\iota_{1}(x))\lambda_{2}(\iota_{2}(x)).$

Given a finite place $v$, we write $\mathfrak{p}_{v}$ for the associated
prime ideal of $\mathcal{O}_{F}$ and $\mathcal{O}_{v}$ for the $v$-adic
completion of $\mathcal{O}_{F}$. Given two ideals $\mathfrak{m},\mathfrak{n}\subset\mathcal{O}_{F}$
prime to $\mathfrak{d}$, we define an open compact subgroup $K(\mathfrak{n},\mathfrak{m})=\prod_{v}K(\mathfrak{n},\mathfrak{m})_{v}$
of $G(\mathbf{A}_{f})$ as follows:
\begin{itemize}
\item $K(\mathfrak{n},\mathfrak{m})_{v}\simeq\mathrm{GL}_{2}(\mathcal{O}_{v})$
if $\mathfrak{p}_{v}\nmid\mathfrak{mnd}$
\item $K(\mathfrak{n},\mathfrak{m})_{v}\simeq\mathcal{O}_{D_{v}}^{\times}$
if $\mathfrak{p}_{v}|\mathfrak{d}$, where $\mathcal{O}_{D_{v}}$
is a choice of maximal order in $D_{v}=D\otimes_{F}F_{v}$
\item $K(\mathfrak{n},\mathfrak{m})_{v}\simeq\left\{ \left(\begin{array}{cc}
a & b\\
c & d\end{array}\right)\in\mathrm{GL}_{2}(\mathcal{O}_{v}),\, c\in(\mathfrak{m}\cap\mathfrak{n})\mathcal{O}_{v},\, d-1\in\mathfrak{n}\mathcal{O}_{v}\right\} .$
\end{itemize}
This is analogous to the classical subgroup $\Gamma_{1}(N)\cap\Gamma_{0}(M)$
of $\mathrm{GL}_{2}(\mathbf{Z})$. The projection of $K(\mathfrak{m},\mathfrak{n})$
onto $G(\mathbf{Q}_{p})$ is contained in the Iwahori if and only
if $\mathfrak{m}\cap\mathfrak{n}\subseteq(p)$. We shall be interested
in the cohomology of local systems on the hyperbolic three-manifold
$Y(K(\mathfrak{n},(p)))$. For brevity, we set $H_{\ast}(\mathfrak{n},\mathcal{A}_{\Omega})=H_{\ast}(Y(K(\mathfrak{n},(p))),\mathcal{A}_{\Omega})$,
etc.

For $i=1,2$, fix a uniformizer $\varpi_{i}$ of $\mathcal{O}_{\mathfrak{p}_{i}}$,
and let $U_{i}$ be the Hecke operator associated with the matrix
$\left(\begin{array}{cc}
1\\
 & \varpi_{i}\end{array}\right)$. We set $U_{p}=U_{1}U_{2}$; this is a controlling operator. Theorem
1.8 follows from the construction in §3.2 together with the following
theorem.

\textbf{Theorem 4.1. }\emph{Let $F/\mathbf{Q}$ be an imaginary quadratic
field in which $p$ splits. Fix a quaternion division algebra $D/F$
split at the primes over $p$, and let $G/\mathbf{Q}$ be the inner
form of $\mathrm{Res}_{F/\mathbf{Q}}\mathrm{GL}_{2}$ associated with
$D$. Choose a slope datum $(U_{t},\Omega,h)$ with $\Omega$ contained
in the two-dimensional null space $\mathcal{W}_{0}\subset\mathcal{W}$
defined in }§\emph{4.3. Then $H_{i}(K^{p},\mathcal{A}_{\Omega})_{\leq h}=0$
unless $i=1$ and $H^{i}(K^{p},\mathcal{D}_{\Omega})_{\leq h}=0$
unless $i=2$ or $i=3$. The Ext spectral sequence degenerates, for
$\Sigma=\Omega$, to isomorphisms\begin{eqnarray*}
H^{2}(K^{p},\mathcal{D}_{\Omega})_{\leq h} & \simeq & \mathrm{Ext}_{A(\Omega)}^{1}(H_{1}(K^{p},\mathcal{A}_{\Omega})_{\leq h},A(\Omega)),\\
H^{3}(K^{p},\mathcal{D}_{\Omega})_{\leq h} & \simeq & \mathrm{Ext}_{A(\Omega)}^{2}(H_{1}(K^{p},\mathcal{A}_{\Omega})_{\leq h},A(\Omega)).\end{eqnarray*}
For any weight $\lambda\in\Omega$, the Tor spectral sequence yields
a long exact sequence \[
0\to H^{2}(K^{p},\mathcal{D}_{\Omega})_{\leq h}\otimes_{A(\Omega)}A(\Omega)/\mathfrak{m}_{\lambda}\to H^{2}(K^{p},\mathcal{D}_{\lambda})_{\leq h}\to\mathrm{Tor}_{A(\Omega)}^{1}(H^{3}(K^{p},\mathcal{D}_{\Omega})_{\leq h},A(\Omega)/\mathfrak{m}_{\lambda})\to0,\]
and the third term vanishes unless $\lambda$ the trivial weight.
Finally, $H^{2}(K^{p},\mathcal{D}_{\Omega})_{\leq h}$, if nonzero,
is a Cohen-Macaulay $A(\Omega)$-module of projective dimension and
grade one, and all of its associated primes have height one.}

\textbf{Lemma 4.3.1. }\emph{For any slope datum, we have $H^{0}(\mathfrak{n},\mathcal{D}_{\lambda})_{\leq h}=0$
and $H^{0}(\mathfrak{n},\mathcal{D}_{\Omega})_{\leq h}=0$.}

\emph{Proof. }Let $\mathbf{D}^{s}$ denote the $k$-Banach dual of
the space $\mathbf{A}^{s}$ of continuous $k$-valued functions $f:\mathbf{Z}_{p}\to k$
which are analytic on each coset of $p^{s}\mathbf{Z}_{p}$. Let $b\in\mathbf{Z}_{p}$
act on $\mathbf{A}^{s}$ by translation, i.e. $(t_{b}f)(x)=f(x+b)$,
and let $t_{b}^{\ast}$ denote the dual action on $\mathbf{D}^{s}$.
By strong approximation it's enough to show that for any $b\neq0$,
the only $t_{b}^{\ast}$-fixed element of $\mathbf{D}^{s}$ is the
zero distribution. By Theorem 1.7.8 of {[}Col{]}, the Amice transform
\[
\mu\mapsto A_{\mu}(T)=\int_{\mathbf{Z}_{p}}(1+T)^{x}d\mu(x)\]
defines an isometric isomorphism of $\mathbf{D}^{s}$ onto the ring
\[
\mathcal{R}^{s}=\left\{ \sum_{n=0}^{\infty}b_{n}T^{n},\, b_{n}\in k\,\mathrm{with}\,|b_{n}|\ll p^{v_{p}(\left\lfloor p^{-s}n\right\rfloor !)}\,\mathrm{as}\, n\to\infty\right\} ,\]
where $\mathcal{R}^{s}$ is equipped with the norm $|\rho|=\mathrm{sup}_{n}|b_{n}(\rho)|p^{-v_{p}(\left\lfloor p^{-s}n\right\rfloor !)}$.
The ring $\mathcal{R}^{s}$ is an integral domain and injects densely
into the Frechet algebra of power series which converge on the open
disk $|T|<p^{-\frac{1}{p^{s}(p-1)}}$. A simple computation gives
$A_{t_{b}^{\ast}\mu}(T)=(1+T)^{b}A_{\mu}(T)$, so $t_{b}^{\ast}\mu=\mu$
implies $\left((1+T)^{b}-1\right)A_{\mu}(T)=0$, whence $A_{\mu}(T)=0$
and $\mu=0$ as desired. $\square$

\emph{Proof of Theorem 4.1. }Choose a slope datum $(\Omega,h)$ with
$\Omega\subset\mathcal{W}_{0}$. By Theorem 1.2, $H_{\ast}(\mathfrak{n},\mathcal{A}_{\Omega})_{\leq h}$
and $H^{\ast}(\mathfrak{n},\mathcal{D}_{\Omega})_{\leq h}$ are torsion
$A(\Omega)$-modules. We examine the Tor spectral sequence for $\Sigma=\lambda\in\Omega$
an arbitrary weight; the sequence reads\[
E_{2}^{i,j}=\mathrm{Tor}_{-i}^{A(\Omega)}(H^{j}(\mathfrak{n},\mathcal{D}_{\Omega})_{\leq h},A(\Omega)/\mathfrak{m}_{\lambda})\Rightarrow H^{i+j}(\mathfrak{n},\mathcal{D}_{\lambda})_{\leq h}.\]
As the ring $A(\Omega)$ is regular of dimension two, nonzero entries
on the $E_{2}$ page may only occur with coordinates $-2\leq i\leq0,0\leq j\leq3$,
and the sequence stabilizes at its $E_{3}$ page. By Lemma 4.3.1,
$H^{0}(\mathfrak{n},\mathcal{D}_{\Omega})_{\leq h}=0$ and so the
entries $E_{2}^{i,0}$ along the zeroth row vanish identically. This
implies that the entries $E_{2}^{-2,1}$ and $E_{2}^{-1,1}$ are stable;
they contribute to gradings on $H^{j}(\mathfrak{n},\mathcal{D}_{\lambda})_{\leq h}$
for $j\in\{-1,0\}$, which both vanish, so they vanish and hence $\mathrm{Tor}_{A(\Omega)}^{i}(H^{1}(\mathfrak{n},\mathcal{D}_{\Omega})_{\leq h},A(\Omega)/\mathfrak{m}_{\lambda})=0$
for any $\lambda\in\Omega$ and any $i>0$. By Proposition A.3, $H^{1}(\mathfrak{n},\mathcal{D}_{\Omega})_{\leq h}$
is either zero or projective, but Theorem 1.2 rules out projectivity,
so $H^{1}(\mathfrak{n},\mathcal{D}_{\Omega})_{\leq h}=0$. This in
turn implies that the entries $E_{2}^{-2,2}$ and $E_{2}^{-1,2}$
are stable; $E_{2}^{-2,2}$ contributes to the grading on $H^{0}(\mathfrak{n},\mathcal{D}_{\Omega})_{\leq h}=0$,
so it vanishes. Thus $\mathrm{Tor}_{i}^{A(\Omega)}(H^{2}(\mathfrak{n},\mathcal{D}_{\Omega})_{\leq h},A(\Omega)/\mathfrak{m}_{\lambda})=0$
for any $\lambda\in\Omega$ and any $i>1$, so $H^{2}(\mathfrak{n},\mathcal{D}_{\Omega})_{\leq h}$
has projective dimension at most one by Proposition A.3; since it's
already a torsion module, its grade and projective dimension are both
at least one. Thus $H^{2}(\mathfrak{n},\mathcal{D}_{\Omega})_{\leq h}$
is a perfect $A(\Omega)$-module of projective dimension one, so the
final sentence of Theorem 4.1 follows from Proposition A.6.

An analysis as in Lemma 3.9 of \cite{BeCrit} shows that $H^{3}(\mathfrak{n},\mathcal{D}_{\Omega})$
is zero unless $\Omega$ contains the trivial weight $\lambda_{0}$,
in which case there is an $A(\Omega)$-module isomorphism $H^{3}(\mathfrak{n},\mathcal{D}_{\Omega})\simeq(A(\Omega)/\mathfrak{m}_{\lambda_{0}})^{|\mathrm{Cl}(F)|}$.
A simple calculation shows $H_{3}(\mathfrak{n},\mathcal{A}_{\Omega})=0$,
so the remainder of the theorem will follow if we can show $H_{2}(\mathfrak{n},\mathcal{A}_{\Omega})_{\leq h}=0$.
Considering the Ext spectral sequence with $\Omega=\lambda$, the
entry $E_{2}^{2,2}$ is stable and necessarily vanishes for any $\lambda\in\Omega$,
so $\mathrm{projdim}_{A(\Omega)}(H_{2}(\mathfrak{n},\mathcal{A}_{\Omega})_{\leq h})\leq1$.
The entry $E_{2}^{1,2}$ is stable for any $\lambda$, and contributes
to the grading on $H^{3}(\mathfrak{n},\mathcal{D}_{\Omega})_{\leq h}$,
which has zero dimensional support, so $E_{2}^{1,2}$ must vanish
for $\lambda$ outside a zero-dimensional analytic set. This rules
out the possibility that $H_{2}(\mathfrak{n},\mathcal{A}_{\Omega})_{\leq h}$
has projective dimension one, and Theorem 1.2 rules out projectivity.
$\square$

\subsubsection*{Two conjectures}

Fix a character $\lambda\in\mathcal{W}_{0}$ with value field $k$.
Note that $\lambda$ is simply a pair of characters $\lambda_{1},\lambda_{2}$
with $\lambda_{i}:\mathbf{Z}_{p}^{\times}\to k^{\times}$. Set $\delta(\lambda_{i})=\frac{\log\lambda_{i}(1+p)}{\log(1+p)}\in k$.

For each prime $\mathfrak{l}\nmid\mathfrak{nd}$, choose a uniformizer
$\varpi_{\mathfrak{l}}$ of $\mathcal{O}_{v_{\mathfrak{l}}}$ and
let $T_{\mathfrak{l}}$ be the double coset operator associated with
the matrix $\left(\begin{array}{cc}
1\\
 & \varpi_{\mathfrak{l}}\end{array}\right).$ For $\mathfrak{l}\nmid\mathfrak{nd}p$, let $S_{\mathfrak{l}}$ be
the double coset operator of $\left(\begin{array}{cc}
\varpi_{\mathfrak{l}}\\
 & \varpi_{\mathfrak{l}}\end{array}\right)$. We define the \emph{$\mathfrak{n}$-}new\emph{ }subspace $H^{2}(\mathfrak{n},\mathcal{D}_{\lambda})^{\mathfrak{n}-\mathrm{new}}$\emph{
}of $H^{2}(\mathfrak{n},\mathcal{D}_{\lambda})$ by intersecting the
kernels of the various level-lowering maps into $H^{2}(\mathfrak{n'},\mathcal{D}_{\lambda})^{\oplus2}$
for varying $\mathfrak{n}'|\mathfrak{n}$ in the usual way. Write
$\mathbf{T}_{\lambda}^{\mathrm{new}}(\mathfrak{n})$ for subalgebra
of $ $$\mathrm{End}_{k}(H^{2}(\mathfrak{n},\mathcal{D}_{\lambda})_{<\infty})$
$ $generated by $T_{\mathfrak{l}}$ for all $\mathfrak{l}\nmid\mathfrak{nd}$
and by $S_{\mathfrak{l}}$ for all $\mathfrak{l}\nmid\mathfrak{nd}p$.

\textbf{Conjecture. }\emph{Given a $k$-algebra homomorphism $\pi:\mathbf{T}_{\lambda}^{\mathrm{new}}(\mathfrak{n})\to\mathbf{C}_{p}$,
there is a continuous semisimple Galois representation $\rho_{\pi}:G_{F}\to\mathrm{GL}_{2}(\mathbf{C}_{p})$,
unique up to conjugation, satisfying the following properties:}

\emph{i. For any $\mathfrak{l}\nmid\mathfrak{n}p$, $\mathrm{tr}\rho_{\pi}(\mathrm{Frob}_{\mathfrak{l}})=\pi(T_{\mathfrak{l}})$
and $\mathrm{det}\rho(\mathrm{Frob}_{\mathfrak{l}})=\mathrm{Nm}(\mathfrak{l})\pi(S_{\mathfrak{l}})$.}

\emph{ii. The prime-to-$p$ Artin conductor of $\rho_{\pi}$ divides
$\mathfrak{nd}$.}

\emph{iii. The Hodge-Tate-Sen weights of $\rho_{\pi}$ at the decomposition
groups over $p$ are $\left\{ 0,\delta(\lambda_{1})+1\right\} $ and
$\left\{ 0,\delta(\lambda_{2})+1\right\} .$}

\emph{iv. The representation $\rho_{\pi}|D_{\mathfrak{p}}$ is trianguline
for all $\mathfrak{p}|p$.}

Next, recall that our group $G$ is the inner form of $\mathrm{Res}_{F/\mathbf{Q}}\mathrm{GL}_{2}$
associated with the quaternion algebra of discriminant $\mathfrak{d}$,
where $\mathfrak{d}$ is a squarefree ideal prime to $p$ and divisible
by an even number of prime ideals. Given a slope datum, let $\mathbf{T}_{\Omega,h}^{\mathfrak{d}}(\mathfrak{n})$
denote the subalgebra of $\mathrm{End}_{A(\Omega)}(H^{2}(\mathfrak{n},\mathcal{D}_{\Omega})_{\leq h})$
generated by $T_{\mathfrak{l}}$ for all $\mathfrak{l}\nmid\mathfrak{nd}$
and by $S_{\mathfrak{l}}$ for all $\mathfrak{l}\nmid\mathfrak{nd}p$.
We wish to formulate a Jacquet-Langlands conjecture relating the algebra
$\mathbf{T}_{\Omega,h}^{\mathfrak{d}}(\mathfrak{n})$ with an analogous
algebra defined using the split form of $G$. 

Let $G'=\mathrm{Res}_{F/\mathbf{Q}}\mathrm{GL}_{2}$ be the split
form, and let $K(\mathfrak{n},\mathfrak{d})$ denote the subgroup
of $G^{'}(\mathbf{A}_{f})$ defined as above (so the second bulleted
condition is vacuous). Let $H_{!}^{2}(K(\mathfrak{n},p\mathfrak{d}),\mathcal{D}_{\Omega})_{\leq h}^{\mathfrak{d}-\mathrm{new}}$
denote the intersection of the parabolic cohomology with the kernels
of the various level-lowering maps, and let $\mathbf{T}_{\Omega,h}(\mathfrak{n},\mathfrak{d})_{!}^{\mathfrak{d}-\mathrm{new}}$
denote the subalgebra of \[
\mathrm{End}_{A(\Omega)}\left(H_{!}^{2}(K(\mathfrak{n},p\mathfrak{d}),\mathcal{D}_{\Omega})_{\leq h}^{\mathfrak{d}-\mathrm{new}}\right)\]
generated by $T_{\mathfrak{l}}$ for all $\mathfrak{l}\nmid\mathfrak{nd}$
and by $S_{\mathfrak{l}}$ for all $\mathfrak{l}\nmid\mathfrak{nd}p$.

\textbf{Conjecture. }\emph{There is an $A(\Omega)$-algebra isomorphism\[
\mathbf{T}_{\Omega,h}^{\mathfrak{d}}(\mathfrak{n})\simeq\mathbf{T}_{\Omega,h}(\mathfrak{n},\mathfrak{d})_{!}^{\mathfrak{d}-\mathrm{new}}.\]
}

\appendix

\section{Some commutative algebra}

In this appendix we collect some results relating the projective dimension
of a module $M$ and its localizations, the nonvanishing of certain
Tor and Ext groups, and the heights of the associated primes of $M$.
We also briefly recall the definition of a perfect module, and explain
their basic properties. These results are presumably well-known to
experts, but they are not given in our basic reference \cite{Matcrt}.

Throughout this subsection, $R$ is a commutative Noetherian ring
and $M$ is a finite $R$-module. Our notations follow \cite{Matcrt},
with one addition: we write $\mathrm{mSupp}(M)$ for the set of maximal
ideals in $\mathrm{Supp}(M)$.

\textbf{Proposition A.1. }\emph{There is an equivalence\[
\mathrm{projdim}_{R}(M)\geq n\Leftrightarrow\mathrm{Ext}_{R}^{n}(M,N)\neq0\,\mathrm{for\, some}\, N\in\mathrm{Mod}_{R}.\]
}See e.g. p. 280 of \cite{Matcrt} for a proof.

\textbf{Proposition A.2. }\emph{The equality\[
\mathrm{projdim}_{R}(M)=\mathrm{sup}_{\mathfrak{m}\in\mathrm{mSupp}(M)}\mathrm{projdim}_{R_{\mathfrak{m}}}(M_{\mathfrak{m}})\]
holds.}

\emph{Proof.} Any projective resolution of $M$ localizes to a projective
resolution of $M_{\mathfrak{m}}$, so $\mathrm{projdim}_{R_{\mathfrak{m}}}(M_{\mathfrak{m}})\leq\mathrm{projdim}_{R}(M)$
for all $\mathfrak{m}$. On the other hand, if $\mathrm{projdim}_{R}(M)\geq n$,
then $\mathrm{Ext}_{R}^{n}(M,N)\neq0$ for some $N$, so $\mathrm{Ext}_{R}^{n}(M,N)_{\mathfrak{m}}\neq0$
for some $\mathfrak{m}$; but $\mathrm{Ext}_{R}^{n}(M,N)_{\mathfrak{m}}\simeq\mathrm{Ext}_{R_{\mathfrak{m}}}^{n}(M_{\mathfrak{m}},N_{\mathfrak{m}})$,
so $\mathrm{projdim}_{R_{\mathfrak{m}}}(M_{\mathfrak{m}})\geq n$
for some $\mathfrak{m}$ by Proposition A.1. $\square$

\textbf{Proposition A.3. }\emph{For $M$ any finite $R$-module, the
equality\[
\mathrm{projdim}_{R}(M)=\mathrm{sup}_{\mathfrak{m}\in\mathrm{mSupp}(M)}\mathrm{sup}\left\{ i|\mathrm{Tor}_{i}^{R}(M,R/\mathfrak{m})\neq0\right\} \]
holds. If furthermore $\mathrm{projdim}_{R}(M)<\infty$ then the equality
\[
\mathrm{projdim}_{R}(M)=\mathrm{sup}\left\{ i|\mathrm{Ext}_{R}^{i}(M,R)\neq0\right\} \]
holds as well.}

\emph{Proof. }The module $\mathrm{Tor}_{i}^{R}(M,R/\mathfrak{m})$
is a finite-dimensional $R/\mathfrak{m}$-vector space, so localization
at $\mathfrak{m}$ leaves it unchanged, yielding\begin{eqnarray*}
\mathrm{Tor}_{i}^{R}(M,R/\mathfrak{m}) & \simeq & \mathrm{Tor}_{i}^{R}(M,R/\mathfrak{m})_{\mathfrak{m}}\\
 & \simeq & \mathrm{Tor}_{i}^{R_{\mathfrak{m}}}(M_{\mathfrak{m}},R_{\mathfrak{m}}/\mathfrak{m}).\end{eqnarray*}
Since the equality $\mathrm{projdim}_{S}(N)=\mathrm{sup}\left\{ i|\mathrm{Tor}_{i}^{S}(N,S/\mathfrak{m}_{S})\neq0\right\} $
holds for any local ring $S$ and any finite $S$-module $N$ (see
e.g. Lemma 19.1.ii of \cite{Matcrt}), the first claim now follows
from Proposition A.2.

For the second claim, we first note that if $S$ is a local ring and
$N$ is a finite $S$-module with $\mathrm{projdim}_{S}(N)<\infty$,
then $\mathrm{projdim}_{S}(N)=\mathrm{sup}\{i|\mathrm{Ext}_{S}^{i}(N,S)\ne0\}$
by Lemma 19.1.iii of \cite{Matcrt}. Hence by Proposition A.2 we have\begin{eqnarray*}
\mathrm{projdim}_{R}(M) & = & \mathrm{sup}_{\mathfrak{m}\in\mathrm{mSupp}(M)}\mathrm{sup}\left\{ i|\mathrm{Ext}_{R_{\mathfrak{m}}}^{i}(M_{\mathfrak{m}},R_{\mathfrak{m}})\neq0\right\} \\
 & = & \mathrm{sup}\{i|\mathrm{Ext}_{R}^{i}(M,R)_{\mathfrak{m}}\neq0\,\mathrm{for\, some}\,\mathfrak{m}\}\\
 & = & \mathrm{sup}\{i|\mathrm{Ext}_{R}^{i}(M,R)\neq0\},\end{eqnarray*}
as desired. $\square$

\textbf{Proposition A.4. }\emph{If $R$ is a Cohen-Macaulay ring,
$M$ is a finite $R$-module of finite projective dimension, and $\mathfrak{p}$
is an associated prime of $M$, then $\mathrm{ht}\mathfrak{p}=\mathrm{projdim}_{R_{\mathfrak{p}}}(M_{\mathfrak{p}})$.
In particular, $\mathrm{ht}\mathfrak{p}\leq\mathrm{projdim}_{R}(M)$.}

\emph{Proof. }Supposing $\mathfrak{p}$ is an associated prime of
$M$, there is an injection $R/\mathfrak{p}\hookrightarrow M$; this
localizes to an injection $R_{\mathfrak{p}}/\mathfrak{p}\hookrightarrow M_{\mathfrak{p}}$,
so $\mathrm{depth}_{R_{\mathfrak{p}}}(M_{\mathfrak{p}})=0$. Now we
compute\begin{eqnarray*}
\mathrm{ht}\mathfrak{p} & = & \mathrm{dim}(R_{\mathfrak{p}})\\
 & = & \mathrm{depth}_{R_{\mathfrak{p}}}(R_{\mathfrak{p}})\;(\mathrm{by\, the\, CM\, assumption})\\
 & = & \mathrm{depth}_{R_{\mathfrak{p}}}(M_{\mathfrak{p}})+\mathrm{projdim}_{R_{\mathfrak{p}}}(M_{\mathfrak{p}})\;(\mathrm{by\, the\, Auslander}-\mathrm{Buchsbaum\, formula)}\\
 & = & \mathrm{projdim}_{R_{\mathfrak{p}}}(M_{\mathfrak{p}}),\end{eqnarray*}
whence the result. $\square$

Now we single out an especially nice class of modules, which are equidimensional
in essentially every sense of the word. Recall the \emph{grade }of
a module $M$, written $\mathrm{grade}_{R}(M)$, is the $\mathrm{ann}_{R}(M)$-depth
of $R$; by Theorems 16.6 and 16.7 of \cite{Matcrt}, \[
\mathrm{grade}_{R}(M)=\mathrm{inf}\{i|\mathrm{Ext}_{R}^{i}(M,R)\neq0\},\]
so quite generally $\mathrm{grade}_{R}(M)\leq\mathrm{projdim}_{R}(M)$.

\textbf{Definition A.5. }\emph{A finite $R$-module $M$ is }perfect
\emph{if $\mathrm{grade}_{R}(M)=\mathrm{projdim}_{R}(M)<\infty$.}

\textbf{Proposition A.6. }\emph{Let $R$ be a Noetherian ring, and
let $M$ be a perfect $R$-module, with $\mathrm{grade}_{R}(M)=\mathrm{projdim}_{R}(M)=d$.
Then for any $\mathfrak{p}\in\mathrm{Supp}(M)$ we have $\mathrm{grade}_{R_{\mathfrak{p}}}(M_{\mathfrak{p}})=\mathrm{projdim}_{R_{\mathfrak{p}}}(M_{\mathfrak{p}})=d$.
If furthermore $R$ is Cohen-Macaulay, then $M$ is Cohen-Macaulay
as well, and every associated prime of $M$ has height $d$.}

\emph{Proof. }The grade of a module can only increase under localization
(as evidenced by the Ext definition above), while the projective dimension
can only decrease; on the other hand, $\mathrm{grade}_{R}(M)\leq\mathrm{projdim}_{R}(M)$
for any finite module over any Noetherian ring. This proves the first
claim.

For the second claim, Theorems 16.6 and 17.4.i of \cite{Matcrt} combine
to yield the formula\[
\mathrm{dim}(M_{\mathfrak{p}})+\mathrm{grade}_{R_{\mathfrak{p}}}(M_{\mathfrak{p}})=\mathrm{dim}(R_{\mathfrak{p}})\]
for any $\mathfrak{p}\in\mathrm{Supp(}M)$. The Auslander-Buchsbaum
formula reads\[
\mathrm{depth}_{R_{\mathfrak{p}}}(M_{\mathfrak{p}})+\mathrm{projdim}_{R_{\mathfrak{p}}}(M_{\mathfrak{p}})=\mathrm{depth}_{R_{\mathfrak{p}}}(R_{\mathfrak{p}}).\]
But $\mathrm{dim}(R_{\mathfrak{p}})=\mathrm{depth}_{R_{\mathfrak{p}}}(R_{\mathfrak{p}})$
by the Cohen-Macaulay assumption, and $\mathrm{grade}_{R_{\mathfrak{p}}}(M_{\mathfrak{p}})=\mathrm{projdim}_{R_{\mathfrak{p}}}(M_{\mathfrak{p}})$
by the first claim. Hence $\mathrm{depth}_{R_{\mathfrak{p}}}(M_{\mathfrak{p}})=\mathrm{dim}(M_{\mathfrak{p}})$
as desired. The assertion regarding associated primes is immediate
from the first claim and Proposition A.4. $\square$

\section{The dimension of regular cuspidal components}

\begin{center}
{\large by James Newton}%
\footnote{Trinity College, Cambridge, CB2 1TQ; jjmn2@cam.ac.uk%
}
\par\end{center}{\large \par}

In this appendix we use the results of the above article to give some
additional evidence for Conjecture 1.5. In the notation and terminology
of Section 1 above, we prove

\textbf{Proposition B.1. }\emph{Any irreducible component of $\mathscr{X}(K^{p})$
containing a cuspidal non-critical regular classical point has dimension
at least $\mathrm{dim}(\mathcal{W})-l(G)$.}

Note that Proposition 5.7.4 of \cite{UrEigen} implies that at least
one of these components has dimension at least $\mathrm{dim}(\mathcal{W})-l(G)$.
This is stated without proof in that reference, and is due to G. Stevens
and E. Urban. We learned the idea of the proof of this result from
E. Urban --- in this appendix we adapt that idea and make essential
use of Theorem 1.1 (in particular the \textquoteleft{}Tor spectral
sequence\textquoteright{}) to provide a fairly simple proof of Proposition
B.1.

We place ourselves in the setting of Section 1. In particular, $G$
is a reductive group over $\mathbf{Q}$, which is split over $\mathbf{Q}_{p}$.
Fix an open compact subgroup $K^{p}\subset G(\mathbf{A}_{f}^{p})$
and a slope datum $(U_{t},\Omega,h)$. Set $q=q(G)$, $l=l(G)$ (these
quantities are defined in the paragraph before the statement of Conjecture
1.5), and suppose that $\mathfrak{M}$ is a maximal ideal of $\mathbf{T}_{\Omega,h}(K^{p})$
corresponding to a cuspidal non-critical regular classical point of
$\mathscr{X}(K^{p})$. Denote by $\mathfrak{m}$ the contraction of
$\mathfrak{M}$ to $A(\Omega)$. Let $\mathscr{P}$ be a minimal prime
of $\mathbf{T}_{\Omega,h}(K^{p})$ contained in $\mathfrak{M}$. Since
$H^{\ast}(K^{p},\mathcal{D}_{\Omega})_{\leq h}$ is a finite faithful
$\mathbf{T}_{\Omega,h}(K^{p})$-module, minimal primes of $\mathbf{T}_{\Omega,h}(K^{p})$
are in bijection with minimal elements of \[
\mathrm{Supp}_{\mathbf{T}_{\Omega,h}(K^{p})}(H^{\ast}(K^{p},\mathcal{D}_{\Omega})_{\leq h});\]
by Theorem 6.5 of \cite{Matcrt}, minimal elements of the latter set
are in bijection with minimal elements of \[
\mathrm{Ass}_{\mathbf{T}_{\Omega,h}(K^{p})}(H^{*}(K^{p},\mathcal{D}_{\Omega})_{\le h}).\]

\textbf{Definition B.2.} \emph{Denote by $r$ the minimal index $i$
such that $\mathscr{P}$ is in the support of $H^{i}(K^{p},\mathcal{D}_{\Omega})_{\le h,\mathfrak{M}}$. }

Let $\wp$ denote the contraction of $\mathscr{P}$ to a prime of
$A(\Omega)_{\mathfrak{m}}$; in particular, $\wp$ is an associated
prime of $H^{r}(K^{p},\mathcal{D}_{\Omega})_{\le h,\mathfrak{M}}$.
The ring $A(\Omega)_{\mathfrak{m}}$ is a regular local ring. The
localisation $A(\Omega)_{\wp}$ is therefore a regular local ring,
with maximal ideal $\wp A(\Omega)_{\wp}$. We let $(x_{1},...,x_{d})$
denote a regular sequence generating $\wp A(\Omega)_{\wp}$. After
multiplying the $x_{i}$ by units in $A(\Omega)_{\wp}$, we may assume
that the $x_{i}$ are in $A(\Omega)$. Note that $(x_{1},...,x_{d})A(\Omega)_{\mathfrak{m}}$
may be a proper submodule of $\wp$. Nevertheless, we have \[
d=\mathrm{dim}(A(\Omega)_{\wp})=\mathrm{ht}(\wp).\]
 We will show that $d\le l$.

Denote by $A_{i}$ the quotient $A(\Omega)_{\wp}/(x_{1},...,x_{i})A(\Omega)_{\wp}$
and denote by $\Sigma_{i}$ the Zariski closed subspace of $\Omega$
defined by the ideal $(x_{1},...,x_{i})A(\Omega)$. The affinoids
$\Sigma_{i}$ may be non-reduced. Note that $A_{i}=A(\Sigma_{i})_{\wp}$
and $A(\Sigma_{i+1})=A(\Sigma_{i})/x_{i+1}A(\Sigma_{i})$.

\textbf{Lemma B.3. }\emph{The space \[
H^{r-d}(K^{p},\mathcal{D}_{\Sigma_{d}})_{\le h,\mathscr{P}}\]
is non-zero.} 

\emph{Proof. }By induction, it suffices to prove the following: let
$i$ be an integer satisfying $0\le i\le d-1$. Suppose \[
H^{r-i}(K^{p},\mathcal{D}_{\Sigma_{i}})_{\le h,\mathscr{P}}\]
 is a non-zero $A_{i}$-module, with $\wp A_{i}$ an associated prime,
and \[
H^{t}(K^{p},\mathcal{D}_{\Sigma_{i}})_{\le h,\mathscr{P}}=0\]
 for every $t<r-i$. Then \[
H^{r-i-1}(K^{p},\mathcal{D}_{\Sigma_{i+1}})_{\le h,\mathscr{P}}\]
 is a non-zero $A_{i+1}$-module, with $\wp A_{i+1}$ an associated
prime, and \[
H^{t}(K^{p},\mathcal{D}_{\Sigma_{i+1}})_{\le h,\mathscr{P}}=0\]
 for every $t<r-i-1$.

Note that the hypothesis of this claim holds for $i=0$, by the minimality
of $r$. Suppose the hypothesis is satisfied for $i$. It will suffice
to show that 
\begin{itemize}
\item $H^{t}(K^{p},\mathcal{D}_{\Sigma_{i+1}})_{\le h,\mathscr{P}}=0$ for
every $t<r-i-1$ 
\item there is an isomorphism of non-zero $A_{i}$-modules \[
\iota:\mathrm{Tor}_{1}^{A_{i}}(H^{r-i}(K^{p},\mathcal{D}_{\Sigma_{i}})_{\le h,\mathscr{P}},A_{i}/x_{i+1}A_{i})\cong H^{r-i-1}(K^{p},\mathcal{D}_{\Sigma_{i+1}})_{\le h,\mathscr{P}}.\]
 
\end{itemize}
Indeed, the left hand side (which we denote by $T$) of the isomorphism
$\iota$ is given by the $x_{i+1}$-torsion in $H^{r-i}(K^{p},\mathcal{D}_{\Sigma_{i}})_{\le h,\mathscr{P}}$,
so a non-zero $A_{i}$-submodule of $H^{r-i}(K^{p},\mathcal{D}_{\Sigma_{i}})_{\le h,\mathscr{P}}$
with annihilator $\wp A_{i}$ immediately gives a non-zero $A_{i+1}$-submodule
of $T$ with annihilator $\wp A_{i+1}$.

Both the claimed facts are shown by studying the localisation at $\mathscr{P}$
of the spectral sequence \[
E_{2}^{s,t}:\mathrm{Tor}_{-s}^{A(\Sigma_{i})}(H^{t}(K^{p},\mathcal{D}_{\Sigma_{i}})_{\le h},A(\Sigma_{i+1}))\Rightarrow H^{s+t}(K^{p},\mathcal{D}_{\Sigma_{i+1}})_{\le h}\]
(cf. Remark 3.1.1). After localisation at $\mathscr{P}$, the spectral
sequence degenerates at $E_{2}$. This is because we have a free resolution
\[
0\rightarrow A_{i}\overset{\times x_{i+1}}{\rightarrow}A_{i}\rightarrow A_{i+1}\rightarrow0\]
of $A_{i+1}$ as an $A_{i}$-module (we use the fact that $x_{i+1}$
is not a zero-divisor in $A_{i}$), so $(E_{2}^{s,t})_{\mathscr{P}}$
vanishes whenever $s\notin\{-1,0\}$. Moreover, since \[
H^{t}(K^{p},\mathcal{D}_{\Sigma_{i}})_{\le h,\mathscr{P}}=0\]
for every $t<r-i$, we know that $(E_{2}^{s,t})_{\mathscr{P}}$ vanishes
for $t<r-i$. The existence of the isomorphism $\iota$ and the desired
vanishing of $H^{t}(K^{p},\mathcal{D}_{\Sigma_{i+1}})_{\le h,\mathscr{P}}$
are therefore demonstrated by the spectral sequence, since the only
non-zero term $(E_{2}^{s,t})_{\mathscr{P}}$ contributing to $(E_{\infty}^{r-i-1})_{\mathscr{P}}$
is given by $s=-1,t=r-i$, whilst $(E_{2}^{s,t})_{\mathscr{P}}=0$
for all $(s,t)$ with $s+t<r-i-1$. $\square$

\textbf{Corollary B.4.}\emph{ We have an inequality $r-d\ge q$. Since
$r\le q+l$ we obtain $d\le l$. In particular $\wp$ has height $\le l$,
so the irreducible component of $\mathbf{T}_{\Omega,h}(K^{p})$ corresponding
to $\mathscr{P}$ has dimension $\ge$ $\mathrm{dim}(\Omega)-l$.} 

\emph{Proof. }It follows from Proposition 3.4.1 (with $\Omega$ replaced
by $\Sigma_{d}$) that \[
H^{i}(K^{p},\mathcal{D}_{\Sigma_{d}})_{\le h,\mathfrak{M}}\]
is zero for $i<q$. Our Lemma therefore implies that $r-d\ge q$.
The conclusion on dimensions follows from the observation made in
Section 3.4 that $\mathbf{T}_{\Omega,h}(K^{p})/\mathscr{P}$ has the
same dimension as $A(\Omega)/\wp$. $\square$ 

Proposition B.1 follows immediately from the Corollary. We have also
shown that if $d=l$ (which would be the case if Conjecture 1.5 is
true), then $r=q+l$.

\bibliographystyle{alpha}

\end{document}